\RequirePackage{tikz-cd}
\tikzset{Rightarrow/.style={double equal sign distance,>={Implies},->},
triple/.style={-,preaction={draw,Rightarrow}},
quadruple/.style={preaction={draw,Rightarrow,shorten >=0pt},shorten >=1pt,-,double,double
distance=0.2pt}}

\documentclass[a4paper, 12pt]{article} 
\usepackage{amsmath} %odgovoran za \newcommand{\Bi}{\operatorname{Im}}
\usepackage{theorem}
\usepackage{amssymb}
\usepackage{amsfonts}
\usepackage{latexsym}
\usepackage{amscd}
\usepackage{diagramb}
\usepackage{graphicx}
\usepackage{enumerate}
\usepackage{cancel}
\usepackage{pxfonts}  %za south-west-arrow
\usepackage{float}
\DeclareFontFamily{U}{mathx}{}
\DeclareFontShape{U}{mathx}{m}{n}{<-> mathx10}{}
\DeclareSymbolFont{mathx}{U}{mathx}{m}{n}
\DeclareMathAccent{\widehat}{0}{mathx}{"70}
\DeclareMathAccent{\widecheck}{0}{mathx}{"71}

\usepackage[hyperindex,breaklinks]{hyperref}

\usepackage{stmaryrd}

\input xy
\xyoption{all}

\usepackage{xspace,colortbl}
\usepackage{color}

\usepackage{arydshln}

\DeclareMathAlphabet{\mathpzc}{OT1}{pzc}{m}{it}

    \newcommand{\axiom}[1]{\textbf{\hypertarget{#1}{(#1)}}}
    \newcommand{\axiomref}[1]{\hyperlink{#1}{(#1)}}

\theoremstyle{plain}
\theoremheaderfont{\normalfont\bfseries}

\newtheorem{thm}{Theorem}[section]

\newtheorem{cor}[thm]{Corollary}
\newtheorem{defn}[thm]{Definition}

\theorembodyfont{\rmfamily}

\newtheorem{rem}[thm]{Remark}
\newtheorem{prop}[thm]{Proposition}
\newtheorem{ex}[thm]{Example}

\newcommand{\qed}{\hfill\quad\fbox{\rule[0mm]{0,0cm}{0,0mm}}  \par\bigskip}

\definecolor{rojo}{rgb}{1,0,0}
\definecolor{grey}{rgb}{0.7,0.7,0.7}

\newcommand{\x}{\mbox{-}}
\newcommand{\s}{\hspace{0,06cm}}
\newcommand{\w}{\hspace{-0,06cm}}

\newcommand{\comp}{\circ}
\newcommand{\iso}{\cong}
\newcommand{\ot}{\otimes}
\newcommand{\crta}{\overline}

\newcommand{\Aa}{{\mathbb A}}
\newcommand{\Bb}{{\mathbb B}}
\newcommand{\Cc}{{\mathbb C}}
\newcommand{\Dd}{{\mathbb D}}
\newcommand{\Ee}{{\mathbb E}}
\newcommand{\Mm}{{\mathbb M}}

\newcommand{\Kk}{{\mathbb K}}
\newcommand{\HH}{{\mathcal H}}

\newcommand{\F}{{\mathcal F}}

\newcommand{\Ps}{\operatorname{\mathbb P s}}

\newcommand{\C}{{\mathcal C}}
\newcommand{\V}{{\mathcal V}}
\newcommand{\K}{{\mathcal K}}
\newcommand{\Ll}{{\mathcal L}}

\newcommand{\M}{{\mathcal M}}

\def\ul{\underline}

\newcommand{\MonDbl}{\operatorname {MonDbl}}
\newcommand{\MMonDbl}{\operatorname {\mathbb M onDbl}}

\newcommand{\Id}{\operatorname {Id}}

\newcommand{\Fun}{\operatorname {Fun}}

\newcommand{\MonBicat}{\operatorname {MonBicat}}
\newcommand{\Para}{\operatorname {Para}}
\newcommand{\coPara}{\operatorname {coPara}}
\newcommand{\Lax}{\operatorname{\mathbb L ax}}
\newcommand{\Pseudo}{\operatorname{\mathbb P seudo}}

\newcommand{\Epsilon}{\varepsilon}

\def\Dd{{\mathbb D}}
\newcommand{\rtr}{\triangleright}
\newcommand{\ltr}{\triangleleft}

\usepackage{yfonts}
\usepackage{upgreek}

\newcommand{\eqlabel}[1]{\label{eq:#1}}
\newcommand{\equref}[1]{(\ref{eq:#1})}
\newcommand{\thlabel}[1]{\label{th:#1}}
\newcommand{\thref}[1]{Theorem~\ref{th:#1}}
\newcommand{\delabel}[1]{\label{de:#1}}
\newcommand{\deref}[1]{Definition~\ref{de:#1}}
\newcommand{\prlabel}[1]{\label{pr:#1}}
\newcommand{\prref}[1]{Proposition~\ref{pr:#1}}
\newcommand{\colabel}[1]{\label{co:#1}}
\newcommand{\coref}[1]{Corollary~\ref{co:#1}}
\newcommand{\rmlabel}[1]{\label{rm:#1}}
\newcommand{\rmref}[1]{Remark~\ref{rm:#1}}
\newcommand{\selabel}[1]{\label{se:#1}}
\newcommand{\seref}[1]{Section~\ref{se:#1}}
\newcommand{\sslabel}[1]{\label{ss:#1}}
\newcommand{\ssref}[1]{Subsection~\ref{ss:#1}}
\newcommand{\ssslabel}[1]{\label{ss:#1}}
\newcommand{\sssref}[1]{Subsection~\ref{ss:#1}}

\newcommand*{\threefrac}[3]{%
  \begin{array}{@{\,}c@{\,}}%
    #1\\
    \hline
    #2\\
    \hline
    #3%
  \end{array}%
}

\newcommand*{\fourfrac}[4]{%
 \begin{array}{@{\,}c@{\,}c@{\,}}%
    #1\\
    \hline
    #2\\
    \hline
    #3\\
    \hline
		#4%
  \end{array}%
	}

\newdir{:=}{{}}
\newcommand{\fit}[3]{\ar@{:=}@/{#3}/[#1] |{\Downarrow #2} }

\begin{document}

\title{Para construction for double categories }

\author{Bojana Femi\'c \vspace{6pt} \\
{\small Mathematical Institite of  \vspace{-2pt}}\\
{\small Serbian Academy of Sciences and Arts } \vspace{-2pt}\\
{\small Kneza Mihaila 36,} \vspace{-2pt}\\
{\small 11 000 Belgrade, Serbia} \vspace{-2pt}\\
{\small femicenelsur@gmail.com}}

\date{}
\maketitle

\begin{abstract}
We introduce Para and coPara double categories for double categories. They rely on a horizontal action $\crta\ot$ of a horizontally monoidal double category 
$\Mm$ on a double category $\Dd$. We prove a series of properties, 
most importantly, we characterize monoidality of $\coPara_\Mm(\Dd)$ in the way that it extends monoidality of $\Dd$ as: lax monoidality of the action $\crta\ot$;  
bistrong commutativity of $\crta\ot$; and purely central premonoidality of $\coPara_\Mm(\Dd)$. 
\end{abstract}

% {\em 2020 Mathematics Subject Classification: 18N10}. \\

%\begin{keywords}

\tableofcontents

\bigskip

\section{Introduction}

Para construction is a categorical framework with applications in computer science. It serves as a ``category'' of parametric morphisms. 
Apart from taking an input $A$, agents take additional data that can be understood as ``parameters'' $P$ to deliver an output $B$. 
This way they are morphisms of the shape $P\crta\ot A\to B$. The objects $P$ live in a monoidal category $\M$, while $A,B$ live in an {\em actegory} $\C$. 
Such parametric morphisms determine objects of a category $\Para(\C)$. %Dually, one can consider coparametric morphisms $A\to P\crta\ot B$ in $\coPara(\C)$. 

The Para construction is studied in the works of David Spivak, David Jaz Myers, Bruno Gavranovi\'c, Mateo Cappucci (see for example 
\cite{FST}, \cite{CGHF}, \cite{CJM}), among others. 
Some of the applications of $\Para(\C)$ are in deep learning, where it was first described in \cite{FST}, open games, and cybernetics. Earlier appearances 
of the 1-categorical notion %version of Para already appears e.g.
can be found in \cite{HT}, with predecessor in \cite{P}, and bicategorical in \cite[Thm. 13.2]{B}.

For $\M$ and $\C$ as above, that $\C$ is an $\M$-actegory (or an {\em $\M$-module category}) it 
means that there is an action functor $\crta\ot: \M\times\C\to\C$ and two natural transformations determining the associativity and unitality of the action 
expressed via a classical pentagonal and a triangular diagram. In this setting there is a %bicategory $\Para_\M(\C)$, the Para construction, whose objects are 
double category $\Para_\M(\C)$, the Para construction, whose objects are that of $\C$, horizontal 1-cells $X\to Y$ are pairs $(P,f)$ where 
$P$ is an object of $\M$ and $f:P\crta\ot X\to Y$ a morphism in $\C$, vertical 1-cells are morphisms in $\C$, and 2-cells 
for horizontal 1-cells $(P,f):X'\to Y'$ and $(Q,g):X\to Y$ and vertical 
1-cells $u:X\to X', v:Y\to Y'$ are morphisms $m:Q\to P$ in $\M$ such that $f(m\crta\ot u)=vg$. 

%In the 1-categorical setting, i
When $\M$ is the trivial category $*$, the action functor becomes an endofunctor $T:\C\to\C$ and the rest of the action data 
deliver $T$ a structure of a monad on $\C$ (so that the morphisms $\mu_A, \eta_A$ are invertible). 
(The Para construction actually can be defined using non-invertible associativity and unitality of the action, the so called 
{\em (co)lax actions}.) In this case, the 1-cells of $\Para_*(\C)$ are precisely the 1-cells of the coKleisli category of $T$, and the 2-cells are trivial. Thus the 
Para construction generalizes the coKleisli category. 

A (co)lax action $\crta\ot: \M\times\C\to\C$ equivalently corresponds to a (co)lax monoidal functor $\M\to\Fun(\C,\C)$, and such functors 
(more precisely, colax monoidal functors into a functor category) are called {\em graded comonads}. The Para construction is the 
Kleisli double category of a graded comonad.

Dually to $\Para_\M(\C)$ there is a double category $\coPara_\M(\C)$ whose horizontal 1-cells have the form $f:X\to P\crta\ot Y$ 
(playing the role of ``coparametric morphisms'') 
and the 2-cells go in the opposite direction. In particular, $\coPara_\M(\C)$ generalizes the Kleisli category of a monad $T$ on $\C$. 

\smallskip

In the present article we introduce the Para and coPara double category for double categories as a generalization of the 
Para construction for 1-categories. In \cite{Fem3} we studied the Kleisli double category $\Kk l(T)$ for (double monads $T$ on) double categories, as well as  
actions of monoidal double categories $\Mm$ on double categories $\Dd$. 
Given that the coPara construction presents a generalization of the Kleisli category, our motivation is to extend our results on $\Kk l(T)$ from \cite{Fem3} 
to the setting of coPara for double categories, facilitating in this way a wider range of applications both of coPara (and Para) in computer science, 
and possibly of the Kleisli construction in pure category theory.  

\smallskip

We prove a variety of properties related to the double category $\coPara_\Mm(\Dd)$, with double categories $\Mm$ and $\Dd$ as above. We will list them below 
in the description of the composition of the paper. We treat structures both in the vertical and in the horizontal direction, and we study when vertical structures 
``lift'' to horizontal ones, assuming existence of certain companions. Our main result is \thref{main eq}, which we cite here: 

\bigskip

\noindent {\bf Theorem 1.}  %\begin{thm} \thlabel{main eq}
Let $\Mm$ and $\Dd$ be horizontally monoidal double categories with a horizontal lax action  $(F, \beta,\nu)$ %$F:\Mm\times\Dd\to\Dd$ 
and assume that $\Mm$ is braided with a double braiding $\Phi$. The following are equivalent: 
\begin{enumerate}
\item there is a structure $(F, F^2,F^0)$ of a lax monoidal horizontal lax action, \vspace{-0,2cm}
\item there is a structure $(F, s,t, q,c)$ of a bistrong and commutative horizontal lax action, \vspace{-0,7cm}
\item there is a purely central horizontally premonoidal structure on $\coPara_\Mm(\Dd)$ with trivial 2-cells $(u,U)$, \vspace{-0,2cm}
\item there is a horizontally monoidal structure on $\coPara_\Mm(\Dd)$ extending that on $\Dd$.
\end{enumerate}
%\end{thm}

The above named properties: lax monoidality, bistrength and commutativity  of an action, premonoidality and monoidality, and double braiding, 
we all study both in their vertical and horizontal versions. 
The equivalence of 1., 2. and 4. is our extension to double categories and to coPara of \cite[Theorems 7.3-7.5]{HF2} for bicategorical monads. 
The equivalence of 3. and 4. relies on our \thref{premon-mon-hor}, which is a horizontal version of our vertical result \cite[Theorem 7.10]{Fem2}, 
stating that there is a 1-1 correspondence between the so called {\em purely central} premonoidal structures on a double category $\Dd$ and its 
monoidal structures. The proof of Theorem 1. is technically pretty involved, for this reason we restrain ourselves to only stress the 
main points of the proof relying on the bicategorical results. 
The implication $(1\Rightarrow 4)$ truncated to 1-categories recovers \cite[Corollary 5.21]{CJM}. 

\smallskip

We now expose the structure of the paper. In Section 2 we recall the preliminaries on double categories, their functors, transformations and modifications, 
we recall the notion of companions, as well as our results from \cite{Fem3} on lifting of vertical transformations and identity vertical modifictaions 
to their horizontal counterparts. In Section 3 we recall the notion of a binoidal structure for double categories from \cite{Fem3} and introduce 
{\em horizontally premonoidal} double categories (vertically premonoidal double categories we studied in \cite{Fem3}). 
%, now we record when vertical notion lifts to the horizontal one. 
In Section 4 we speak about locally cubical bicategories from \cite{GG}: we use this setting to fix ideas on two different ways of considering monoidal 
structures on double categories. For explicit definitions of monoidal (lax) double functors, transformations and modifications in what we will 
consider to be vertical sense, we refer the reader to \cite{GGV}. We modify these into what we will consider to be their horizontal versions, so that 
so defined monoidal structures generalize the corresponding underlying bicategorical ones. 
We relate the two groups of structures via a morphism of locally cubical bicategories. We use the term (1- and 2-dimensional) {\em cocycles} 
to simplify the description of the axioms defining monoidality structures on all levels of cells. 
In Section 5 we recall the notion of vertical and horizontal action 
of double categories from \cite{Fem3}, we introduce vertical and horizontal notions of a {\em double braiding}, and of vertically and horizontally 
{\em lax monoidal} action. The latter generalizes the notion of monoidal (vertical and horizontal) double monad from \cite{GGV}. 
In \prref{ex-braid-prem} we provide an example of a (horizontally) premonoidal double category stemming from a (horizontally) braided 
monoidal double category. 

In Section 6 we introduce (left and right, vertical and horizontal) strengths on an action of double categories, as well as bistrengths and commutativity 
of an action. Our first more relevant result is \prref{lax monoidal-comm h} (equivalence of lax monoidality and bistrong commutativity of an action) %$\crta\ot$
followed by its vertical analogue. Section 7 is the last section: here we introduce the double categories Para and coPara for a given action 
$\crta\ot:\Mm\times\Dd\to\Dd$. We introduce the notion of the canonical embedding $P:\Dd\to\coPara_\Mm(\Dd)$. In Subsection 7.1 we prove that a left 
strength induces a horizontal icon and a left action $\rtr:\Dd\times\coPara_\Mm(\Dd)\to\coPara_\Mm(\Dd)$, and that a bistrong lax action makes 
$\coPara_\Mm(\Dd)$ premonoidal (all in horizontal sense). In Subsection 7.2 we recall our result from \cite{Fem3} that relates by equivalence purely 
central premonoidal 
double categories with monoidal double categories (in the vertical sense), and prove its horizontal version in \thref{premon-mon-hor}. In \thref{iff} 
we prove 1-1 correspondence between horizontal strengths and extensions of horizontal actions of $\Dd$ on itself, and in \thref{main eq} we prove the 
above highlighted Theorem 1. We record which results from Section 7 have their vertical versions.

\section{Lifting of vertical to horizontal structures}

In the literature there are a couple of versions of horizontal and vertical transformations between (weak) double functors. 
In order to be coherent with our terminology used in author's previous articles, in this section we record the definitions that we will be using throughout. 
We also recall the definition of modifications, as well as of companions and conjoints, and we record few important results 
from \cite{Fem3} on how vertical transformations and identity vertical modifications lift to their horizontal counterparts.

\subsection{Basics on double categories and notation}

A double category is a category internal to the category of categories. A weaker notion is a pseudodouble category, which is a pseudocategory internal to the 2-category of categories. In \cite[Section~7.5]{GP:L} it is proved that any pseudodouble category is double equivalent to a double category. To simplify the writing and computation we will deal with double categories. 
For a double category $\Dd$ we denote by $\Dd_0$ the category of objects, and by $\Dd_1$ the category of morphisms. 
For horizontal 1-cells we will say shortly 1h-cells, for vertical 1-cells we will say 1v-cells, and squares we will just call 2-cells. 
In our convention the horizontal direction is weak and the vertical direction is strict. 
The underlying horizontal 2-category of $\Dd$ will be denoted by $\HH(\Dd)$. It consists of objects, 1h-cells and those 2-cells 
whose 1v-cells are identities. 
Composition of 1h-cells $A\stackrel{f}{\to}A'\stackrel{f'}{\to}A''$ we will write as $f'f$ or $[f\,\, \vert\,\, f']$ and vertical composition 
of 1v-cells $A\stackrel{u}{\to}\tilde A\stackrel{v}{\to}\tilde{\tilde A}$ as a fraction $\frac{u}{v}$. Similarly, to simplify the writing or avoid the use of too many pasting diagrams, we will denote horizontal composition of 2-cells $\alpha$ and $\beta$, 
where $\alpha$ acts first, by $[\alpha\vert\beta]$, and vertical composition of 2-cells $\alpha$ and $\gamma$, where $\alpha$ acts first, by $\frac{\alpha}{\gamma}$. In the images below we illustrate the mnemotechnical value of such notation 
$$[\alpha\vert \beta]=
\scalebox{0.8}{
\bfig
\putmorphism(0,200)(1,0)[A`B`f]{450}1a
\putmorphism(460,200)(1,0)[`C`f']{450}1a
\putmorphism(0,220)(0,-1)[\phantom{Y_2}``u]{380}1l
\putmorphism(450,220)(0,-1)[\phantom{Y_2}``u']{380}1r
\putmorphism(890,220)(0,-1)[\phantom{Y_2}``u'']{380}1r
\put(150,0){\fbox{$\alpha$}}
\put(600,0){\fbox{$\beta$}}
\putmorphism(0,-150)(1,0)[\tilde A`\tilde B`g]{450}1b
\putmorphism(460,-150)(1,0)[`\tilde C, `g']{450}1b
\efig}
\qquad\quad
\frac{\alpha}{\gamma}=
\scalebox{0.8}{
\bfig
\putmorphism(0,350)(1,0)[A`B`f]{450}1a
\putmorphism(0,370)(0,-1)[\phantom{Y_2}``u]{380}1l
\putmorphism(450,370)(0,-1)[\phantom{Y_2}``u']{380}1r
\put(150,150){\fbox{$\alpha$}}
\put(150,-210){\fbox{$\gamma$}}
\putmorphism(0,10)(0,-1)[\phantom{Y_2}``v]{380}1l
\putmorphism(450,10)(0,-1)[\phantom{Y_2}``v']{380}1r
\putmorphism(0,0)(1,0)[\tilde A`\tilde B`g]{450}1b
\putmorphism(0,-360)(1,0)[\tilde{\tilde A}`\tilde{\tilde B}.`h]{450}1b
\efig}
$$
For 2-cells of the form of $a$ as below we will say that they are {\em horizontally globular}, while for those of the form of $b$ 
we will say that they are {\em vertically globular}. 
$$ 
\scalebox{0.9}{
\bfig
\putmorphism(-150,170)(1,0)[A`B`f]{460}1a
\putmorphism(-150,-180)(1,0)[A`B`g]{460}1a
\putmorphism(-150,170)(0,-1)[\phantom{Y_2}``=]{350}1l
\putmorphism(310,170)(0,-1)[\phantom{Y_2}``=]{350}1r
\put(20,10){\fbox{$a$}}
\efig}
\qquad\qquad
\scalebox{0.9}{
\bfig
\putmorphism(-150,170)(1,0)[A`A`=]{460}1a
\putmorphism(-150,-180)(1,0)[\tilde A`\tilde B`=]{460}1a
\putmorphism(-150,170)(0,-1)[\phantom{Y_2}``u]{350}1l
\putmorphism(310,170)(0,-1)[\phantom{Y_2}``v]{350}1r
\put(20,0){\fbox{$b$}}
\efig}
$$ 
Strict double functors $F:\Dd\to\Ee$ are functors internal in the category of categories. They consist of functors $F_0:\Dd_0\to\Ee_0$ 
(on objects and 1v-cells) and $F_1:\Dd_1\to\Ee_1$ (on 1h-cells and 2-cells). From weak double functors we will use lax and pseudo ones, 
they are weak in the horizontal direction. For more on double categories we recommend \cite{GP:L, MG}. We pass to the definitions of 
transformations of (weak) double functors.

\begin{defn} \delabel{hor nat tr}
A {\em horizontal oplax transformation} $\alpha$ between lax double functors $F,G\colon \Aa\to\Bb$ consists of the following:
\begin{enumerate}
\item for every 0-cell $A$ in $\Aa$ a 1h-cell $\alpha(A)\colon F(A)\to G(A)$ in $\Bb$,
\item for every 1v-cell $u\colon A\to A'$ in $\Aa$ a 2-cell in $\Bb$:
$$
\scalebox{0.86}{
\bfig
\putmorphism(-150,50)(1,0)[F(A)`G(A)`\alpha(A)]{560}1a
\putmorphism(-150,-320)(1,0)[F(A')`G(A')`\alpha(A')]{600}1a
\putmorphism(-180,50)(0,-1)[\phantom{Y_2}``F(u)]{370}1l
\putmorphism(410,50)(0,-1)[\phantom{Y_2}``G(u)]{370}1r
\put(30,-110){\fbox{$\alpha^u$}}
\efig}
$$
\item for every 1h-cell $f\colon A\to B$ in $\Aa$ a 2-cell in $\Bb$:
$$
\scalebox{0.86}{
\bfig
 \putmorphism(-170,500)(1,0)[F(A)`F(B)`F(f)]{540}1a
 \putmorphism(360,500)(1,0)[\phantom{F(f)}`G(B) `\alpha(B)]{560}1a
 \putmorphism(-170,120)(1,0)[F(A)`G(A)`\alpha(A)]{540}1a
 \putmorphism(360,120)(1,0)[\phantom{G(B)}`G(A) `G(f)]{560}1a
\putmorphism(-180,500)(0,-1)[\phantom{Y_2}``=]{380}1r
\putmorphism(940,500)(0,-1)[\phantom{Y_2}``=]{380}1r
\put(280,310){\fbox{$\delta_{\alpha,f}$}}
\efig}
$$
\end{enumerate}
so that the following are satisfied: 
\begin{itemize}
\item (coherence with compositors for $\delta_{\alpha,-}$): for any composable 1h-cells $f$ and $g$ in $\Aa$ the 2-cell 
$\delta_{\alpha,gf}$ satisfies: \\
{\em \axiom{h.o.t.-1}} 
$$\scalebox{0.82}{
\bfig
  \putmorphism(-750,200)(1,0)[F(A)`\phantom{F(B)}`F(f)]{600}1a
\putmorphism(-130,200)(1,0)[F(A)`F(C)`F(g)]{580}1a
\putmorphism(-730,200)(0,-1)[\phantom{Y_2}``=]{400}1r
\putmorphism(420,200)(0,-1)[\phantom{Y_2}``=]{400}1r
\putmorphism(-730,-210)(0,-1)[\phantom{Y_2}``=]{400}1r
\putmorphism(1030,-210)(0,-1)[\phantom{Y_2}``=]{400}1r
 \putmorphism(450,-210)(1,0)[F(C)`G(C) `\alpha(C)]{580}1a
  \putmorphism(-750,-210)(1,0)[F(A)`\phantom{F(B)}`F(gf)]{1200}1a
\put(-270,20){\fbox{$F_{gf}$}}
 \putmorphism(-750,-615)(1,0)[F(A)`G(A)`\alpha(A)]{620}1a
 \putmorphism(-120,-615)(1,0)[\phantom{F(B)}`G(C) `G(gf)]{1170}1a
\put(-250,-420){\fbox{$\delta_{\alpha,gf}$}}
\efig}= 
\scalebox{0.82}{
\bfig
 \putmorphism(450,150)(1,0)[F(B)`F(C) `F(g)]{680}1a
 \putmorphism(1120,150)(1,0)[\phantom{F(B)}`G(C) `\alpha(C)]{600}1a
\put(1000,-100){\fbox{$\delta_{\alpha,g}$}}

  \putmorphism(-150,-300)(1,0)[F(A)` F(B) `F(f)]{600}1a
\putmorphism(450,-300)(1,0)[\phantom{F(A)}` G(B) `\alpha(B)]{680}1a
 \putmorphism(1120,-300)(1,0)[\phantom{F(A)}`G(C) ` G(g)]{620}1a

\putmorphism(450,150)(0,-1)[\phantom{Y_2}``=]{450}1l
\putmorphism(1710,150)(0,-1)[\phantom{Y_2}``=]{450}1r

 \putmorphism(-150,-750)(1,0)[F(A)`G(A)`\alpha(A)]{600}1a
 \putmorphism(450,-750)(1,0)[\phantom{F(B)}`G(B) `G(f)]{680}1a
 \putmorphism(1120,-750)(1,0)[\phantom{F(B)}`G(C) `G(g)]{620}1a

\putmorphism(-180,-300)(0,-1)[\phantom{Y_2}``=]{450}1r
\putmorphism(1040,-300)(0,-1)[\phantom{Y_2}``=]{450}1r
\put(350,-540){\fbox{$\delta_{\alpha,f}$}}
\put(1000,-960){\fbox{$G_{gf}$}}

 \putmorphism(450,-1200)(1,0)[G(A)` G(C) `G(gf)]{1300}1a

\putmorphism(450,-750)(0,-1)[\phantom{Y_2}``=]{450}1l
\putmorphism(1750,-750)(0,-1)[\phantom{Y_2}``=]{450}1r
\efig}
$$ 
(coherence with unitors for $\delta_{\alpha,-}$): for any object $A\in\Aa$: 
$$\text{{\em \axiom{h.o.t.-2}}}  \qquad\quad
\scalebox{0.86}{
\bfig
 \putmorphism(-150,420)(1,0)[F(A)`F(A)`=]{500}1a
\putmorphism(-180,420)(0,-1)[\phantom{Y_2}``=]{370}1l
\putmorphism(320,420)(0,-1)[\phantom{Y_2}``=]{370}1r
 \putmorphism(-150,50)(1,0)[F(A)`F(A)`F(1_A)]{500}1a
 \put(-80,250){\fbox{$F_A$}} %(0,-180)
\putmorphism(330,50)(1,0)[\phantom{F(A)}`G(A) `\alpha(A)]{560}1a
 \putmorphism(-170,-350)(1,0)[F(A)`G(A)`\alpha(A)]{520}1a
 \putmorphism(350,-350)(1,0)[\phantom{F(A)}`G(A) `G(1_A)]{560}1a

\putmorphism(-180,50)(0,-1)[\phantom{Y_2}``=]{400}1l
\putmorphism(910,50)(0,-1)[\phantom{Y_2}``=]{400}1r
\put(240,-150){\fbox{$\delta_{\alpha,1_A}$}}
\efig}
\quad=\quad
\scalebox{0.86}{
\bfig
 \putmorphism(-150,420)(1,0)[F(A)`G(A)`\alpha(A)]{500}1a
\putmorphism(-180,420)(0,-1)[\phantom{Y_2}``=]{370}1l
\putmorphism(320,420)(0,-1)[\phantom{Y_2}``=]{370}1r
  \put(-100,230){\fbox{$\Id_{\alpha(A)}$}} %(0,-180)
\putmorphism(-150,50)(1,0)[F(A)` \phantom{Y_2} `\alpha(A)]{450}1a

\putmorphism(350,50)(1,0)[G(A)` G(A) `=]{470}1a
\putmorphism(330,-300)(1,0)[G(A)` G(A) `G(1_A)]{480}1b
\putmorphism(330,50)(0,-1)[\phantom{Y_2}``=]{350}1l
\putmorphism(800,50)(0,-1)[\phantom{Y_2}``=]{350}1r
\put(470,-150){\fbox{$G_A$}}
\efig}
$$

\item (coherence with vertical composition and identity for $\alpha^\bullet$): for any composable 1v-cells $u$ and $v$ in $\Aa$:
$$\text{{\em \axiom{h.o.t.-3}}} \label{h.o.t.-3} \qquad\alpha^{\frac{u}{v}}=\frac{\alpha^u}{\alpha^v}\quad\qquad\text{ and}\quad\qquad
\text{{\em \axiom{h.o.t.-4}}} \label{h.o.t.-4} \qquad\alpha^{1^A}=\Id_{\alpha(A)};$$

\item (oplax naturality of 2-cells):
for every 2-cell in $\Aa$
$\scalebox{0.86}{
\bfig
\putmorphism(-150,50)(1,0)[A` B`f]{400}1a
\putmorphism(-150,-270)(1,0)[A'`B' `g]{400}1b
\putmorphism(-170,50)(0,-1)[\phantom{Y_2}``u]{320}1l
\putmorphism(250,50)(0,-1)[\phantom{Y_2}``v]{320}1r
\put(0,-140){\fbox{$a$}}
\efig}$ 
the following identity in $\Bb$ must hold:\\
$\text{{\em \axiom{h.o.t.-5}}}$ 
$$
\scalebox{0.86}{
\bfig
\putmorphism(-150,500)(1,0)[F(A)`F(B)`F(f)]{600}1a
 \putmorphism(450,500)(1,0)[\phantom{F(A)}`G(B) `\alpha(B)]{640}1a

 \putmorphism(-150,50)(1,0)[F(A')`F(B')`F(g)]{600}1a
 \putmorphism(450,50)(1,0)[\phantom{F(A)}`G(B') `\alpha(B')]{640}1a

\putmorphism(-180,500)(0,-1)[\phantom{Y_2}``F(u)]{450}1l
\putmorphism(450,500)(0,-1)[\phantom{Y_2}``]{450}1r
\putmorphism(300,500)(0,-1)[\phantom{Y_2}``F(v)]{450}0r
\putmorphism(1100,500)(0,-1)[\phantom{Y_2}``G(v)]{450}1r
\put(0,260){\fbox{$F(a)$}}
\put(700,270){\fbox{$\alpha^v$}}

\putmorphism(-150,-400)(1,0)[F(A')`G(A') `\alpha(A')]{640}1a
 \putmorphism(450,-400)(1,0)[\phantom{A'\ot B'}` G(B') `G(g)]{680}1a

\putmorphism(-180,50)(0,-1)[\phantom{Y_2}``=]{450}1l
\putmorphism(1120,50)(0,-1)[\phantom{Y_3}``=]{450}1r
\put(320,-200){\fbox{$\delta_{\alpha,g}$}}

\efig}
\quad=\quad
\scalebox{0.86}{
\bfig
\putmorphism(-150,500)(1,0)[F(A)`F(B)`F(f)]{600}1a
 \putmorphism(450,500)(1,0)[\phantom{F(A)}`G(B) `\alpha(B)]{680}1a
 \putmorphism(-150,50)(1,0)[F(A)`G(A)`\alpha(A)]{600}1a
 \putmorphism(450,50)(1,0)[\phantom{F(A)}`G(B) `G(f)]{680}1a

\putmorphism(-180,500)(0,-1)[\phantom{Y_2}``=]{450}1r
\putmorphism(1100,500)(0,-1)[\phantom{Y_2}``=]{450}1r
\put(350,260){\fbox{$\delta_{\alpha,f}$}}
\put(650,-180){\fbox{$G(a)$}}

\putmorphism(-150,-400)(1,0)[F(A')`G(A') `\alpha(A')]{640}1a
 \putmorphism(490,-400)(1,0)[\phantom{F(A')}` G(B'). `G(g)]{640}1a

\putmorphism(-180,50)(0,-1)[\phantom{Y_2}``F(u)]{450}1l
\putmorphism(450,50)(0,-1)[\phantom{Y_2}``]{450}1l
\putmorphism(610,50)(0,-1)[\phantom{Y_2}``G(u)]{450}0l %470
\putmorphism(1120,50)(0,-1)[\phantom{Y_3}``G(v)]{450}1r
\put(40,-180){\fbox{$\alpha^u$}} %(0,-180)
\efig}
$$
\end{itemize}
A {\em horizontal lax transformation} is defined by using the opposite direction of the 2-cells $\delta_{\alpha,f}$ in item 3 
and accommodating the axioms correspondingly. \\
A {\em horizontal pseudonatural transformation} is a horizontal oplax transformation for which the 2-cells $\delta_{\alpha,f}$ are isomorphisms. 
A {\em horizontal strict transformation} is a horizontal oplax transformation for which the 2-cells $\delta_{\alpha,f}$ are identities. 
\end{defn}

\begin{defn} \delabel{vlt}
A {\em vertical lax transformation} $\alpha_0$ between lax double functors $F,G\colon \Aa\to\Bb$ consists of: 
\begin{enumerate}
\item a 1v-cell $\alpha_0(A)\colon F(A)\to G(A)$ in $\Bb$ for every 0-cell $A$ in $\Aa$; 
\item %\indent \hspace{1cm}  (I) \hspace{5,7cm} (II) \vspace{-0,5cm}
for every 1h-cell $f\colon A\to B$ in $\Aa$ a 2-cell in $\Bb$:
$$ 
\scalebox{0.9}{
\bfig
\putmorphism(-150,180)(1,0)[F(A)`F(B)`F(f)]{560}1a
\putmorphism(-150,-190)(1,0)[G(A)`G(B)`G(f)]{600}1a
\putmorphism(-150,180)(0,-1)[\phantom{Y_2}``\alpha_0(A)]{370}1l
\putmorphism(410,180)(0,-1)[\phantom{Y_2}``\alpha_0(B)]{370}1r
\put(-30,20){\fbox{$(\alpha_0)_f$}}
\efig}
$$
\item for every 1v-cell $u\colon A\to A'$ in $\Aa$ a 2-cell in $\Bb$:
$$
\scalebox{0.9}{
\bfig
 \putmorphism(-90,500)(1,0)[F(A)`F(A) `=]{540}1a
\putmorphism(450,500)(0,-1)[\phantom{Y_2}`F(\tilde A) `F(u)]{400}1r
\putmorphism(-90,-300)(1,0)[G(\tilde A)`G(\tilde A) `=]{540}1a
\putmorphism(450,100)(0,-1)[\phantom{Y_2}``\alpha_0(\tilde A)]{400}1r
\putmorphism(-100,100)(0,-1)[\phantom{Y_2}``G(u)]{400}1l
\putmorphism(-100,500)(0,-1)[\phantom{Y_2}`G(A) `\alpha_0(A)]{400}1l
\put(80,70){\fbox{$\alpha_0^u$}}
\efig}
$$
\end{enumerate}
which need to satisfy: 
\begin{itemize}
\item (coherence with compositors for $(\alpha_0)_\bullet$): for any composable 1h-cells $f$ and $g$ in $\Aa$: \\
$$\text{{\em \axiom{v.l.t.\x 1}}}  \quad
\scalebox{0.84}{
\bfig
\putmorphism(-150,500)(1,0)[F(A)`F(B)`F(f)]{600}1a
 \putmorphism(450,500)(1,0)[\phantom{F(A)}`F(C) `F(g)]{600}1a

 \putmorphism(-150,50)(1,0)[G(A)`G(B)`G(f)]{600}1a
 \putmorphism(450,50)(1,0)[\phantom{F(A)}`G(C) `G(g)]{600}1a

\putmorphism(-180,500)(0,-1)[\phantom{Y_2}``\alpha_0(A)]{450}1l
\putmorphism(450,500)(0,-1)[\phantom{Y_2}``]{450}1r
\putmorphism(300,500)(0,-1)[\phantom{Y_2}``\alpha_0(B)]{450}0r
\putmorphism(1060,500)(0,-1)[\phantom{Y_2}``\alpha_0(C)]{450}1r
\put(0,260){\fbox{$(\alpha_0)_f$}}
\put(660,270){\fbox{$(\alpha_0)_g$}}

\putmorphism(-150,-400)(1,0)[G(A)`G(C) `G(gf)]{1220}1a

\putmorphism(-180,50)(0,-1)[\phantom{Y_2}``=]{450}1l
\putmorphism(1060,50)(0,-1)[\phantom{Y_3}``=]{450}1r
\put(320,-170){\fbox{$G_{gf}$}}

\efig}
=\quad\hspace{-0,2cm}
\scalebox{0.84}{
\bfig
\putmorphism(-160,500)(1,0)[F(A)`F(B)`F(f)]{580}1a
 \putmorphism(420,500)(1,0)[\phantom{F(A)}`F(C) `F(g)]{620}1a

 \putmorphism(-160,50)(1,0)[F(A)`F(C)`F(gf)]{1260}1a

\putmorphism(-180,500)(0,-1)[\phantom{Y_2}``=]{450}1r
\putmorphism(1060,500)(0,-1)[\phantom{Y_2}``=]{450}1l
\put(330,260){\fbox{$F_{gf}$}}
\put(300,-180){\fbox{$(\alpha_0)_{gf}$}}

\putmorphism(-180,50)(0,-1)[\phantom{Y_2}``\alpha_0(A)]{450}1r
\putmorphism(1060,50)(0,-1)[\phantom{Y_3}``\alpha_0(C)]{450}1l
\putmorphism(-150,-400)(1,0)[G(A)`G(C) `G(gf)]{1220}1a
\efig}
$$

(coherence with unitors for $(\alpha_0)_\bullet$): for any object $A$ in $\Aa$:
$$\text{{\em \axiom{v.l.t.\x 2}}}  \qquad\qquad %(\alpha_0)_{1_a}=\Id^{\alpha_0(A)}  \hspace{4cm}$$
\scalebox{0.86}{
\bfig
 \putmorphism(-170,420)(1,0)[F(A)`F(A)`=]{500}1a
\putmorphism(-180,420)(0,-1)[\phantom{Y_2}``=]{370}1l
\putmorphism(280,420)(0,-1)[\phantom{Y_2}``=]{370}1r
  \put(-40,250){\fbox{$F_A$}} 
\putmorphism(-170,50)(1,0)[F(A)` F(A) `F(1_A)]{450}1a

\putmorphism(-170,50)(0,-1)[\phantom{Y_2}``\alpha_0(A)]{350}1l
\putmorphism(280,50)(0,-1)[\phantom{Y_2}``\alpha_0(A)]{350}1r
\putmorphism(-170,-300)(1,0)[G(A)` G(A) `G(1_A)]{460}1b
\put(-100,-140){\fbox{$(\alpha_0)_{1_A}$}}
\efig}\quad
=
\scalebox{0.86}{
\bfig
 \putmorphism(-170,420)(1,0)[F(A)`F(A)`=]{500}1a
\putmorphism(-180,420)(0,-1)[\phantom{Y_2}``\alpha_0(A)]{370}1l
\putmorphism(280,420)(0,-1)[\phantom{Y_2}``\alpha_0(A)]{370}1r
  \put(-120,220){\fbox{$\Id^{\alpha_0(A)}$}} %(0,-180)
\putmorphism(-170,50)(1,0)[G(A)` G(A) `=]{450}1a

\putmorphism(-170,50)(0,-1)[\phantom{Y_2}``=]{350}1l
\putmorphism(280,50)(0,-1)[\phantom{Y_2}``=]{350}1r
\putmorphism(-170,-300)(1,0)[G(A)` G(A) `G(1_A)]{460}1b
\put(-40,-140){\fbox{$G_A$}}
\efig}
$$
\item (coherence with vertical composition for $\alpha_0^\bullet$): for any composable 1v-cells $u$ and $v$ in $\Aa$: 
$$ \text{{\em \axiom{v.l.t.\x 3}}}   \qquad\quad 
\alpha_0^{\frac{u}{u'}}=
\scalebox{0.86}{
\bfig 
 \putmorphism(-150,500)(1,0)[F(A)`F(A) `=]{500}1a
\putmorphism(-130,500)(0,-1)[\phantom{Y_2}`G(A) `\alpha_0(A)]{400}1l
\put(0,250){\fbox{$\alpha_0^u$}}
\putmorphism(-150,-300)(1,0)[G(\tilde A)`G(\tilde A) `=]{460}1a
\putmorphism(-130,110)(0,-1)[\phantom{Y_2}``G(u)]{400}1l
\putmorphism(380,500)(0,-1)[\phantom{Y_2}` `F(u)]{400}1r
\putmorphism(380,100)(0,-1)[F(\tilde A)` `\alpha_0(\tilde A)]{400}1l
\putmorphism(480,110)(1,0)[`F(\tilde A)`=]{460}1a
\putmorphism(390,-700)(1,0)[\phantom{G(A)}`G(\tilde{\tilde A})`=]{570}1a
\putmorphism(920,100)(0,-1)[\phantom{(B, \tilde A')}``F(u')]{400}1r
\putmorphism(920,-300)(0,-1)[F(\tilde{\tilde A})`` \alpha_0(\tilde{\tilde A})]{400}1r
\putmorphism(400,-300)(0,-1)[\phantom{(B, \tilde A)}`G(\tilde{\tilde A}) `G(u')]{400}1l
\put(530,-190){\fbox{$\alpha_0^{u'}$}}
\efig}
$$
(coherence with vertical identity for $\alpha^\bullet$): for any object $A$ in $\Aa$: \\
$$\text{{\em \axiom{v.l.t.\x 4}}}  \qquad\qquad \alpha_0^{1^A}=\Id^{\alpha_0(A)}  \hspace{4cm}$$
\item (lax naturality of 2-cells):
for every 2-cell in $\Aa$
$\scalebox{0.86}{
\bfig
\putmorphism(-150,50)(1,0)[A` B`f]{400}1a
\putmorphism(-150,-270)(1,0)[\tilde A` \tilde B `g]{400}1b
\putmorphism(-150,50)(0,-1)[\phantom{Y_2}``u]{320}1l
\putmorphism(250,50)(0,-1)[\phantom{Y_2}``v]{320}1r
\put(0,-140){\fbox{$a$}}
\efig}$ 
the following identity in $\Bb$ must hold:\\
$$\text{{\em \axiom{v.l.t.\x 5}}}  \quad
\scalebox{0.84}{
\bfig
 \putmorphism(-130,500)(1,0)[F(A)`F(A) `=]{560}1a
 \putmorphism(530,500)(1,0)[` `F(f)]{400}1a
\putmorphism(-130,500)(0,-1)[\phantom{Y_2}`G(A) `\alpha_0(A)]{450}1l
\put(40,50){\fbox{$\alpha_0^u$}}
\putmorphism(-130,-400)(1,0)[G(\tilde A)` `=]{470}1a
\putmorphism(-150,50)(0,-1)[\phantom{Y_2}``G(u)]{450}1l
\putmorphism(450,50)(0,-1)[\phantom{Y_2}`G(\tilde A)`\alpha_0(\tilde A)]{450}1l
\putmorphism(450,500)(0,-1)[\phantom{Y_2}`F(\tilde A) `F(u)]{450}1l
\put(600,260){\fbox{$F(a)$}}
\putmorphism(420,50)(1,0)[\phantom{(B, \tilde A)}``F(g)]{500}1a
\putmorphism(1030,50)(0,-1)[\phantom{(B, A')}`G(\tilde B)`]{450}1r
\putmorphism(1010,50)(0,-1)[\phantom{(B, A')}``\alpha_0(\tilde B)]{450}0r
\putmorphism(1030,500)(0,-1)[F(B)`F(\tilde B)`]{450}1r
\putmorphism(1010,500)(0,-1)[``F(v)]{450}0r
\putmorphism(420,-400)(1,0)[\phantom{(B, \tilde A)}``G(g)]{500}1a
\put(600,-170){\fbox{$(\alpha_0)_g$}}
\efig}
\quad\hspace{-0,1cm}=\quad\hspace{-0,14cm}
\scalebox{0.84}{
\bfig
 \putmorphism(-130,500)(1,0)[F(A)`F(B) `F(f)]{560}1a
 \putmorphism(550,500)(1,0)[` `=]{380}1a
\putmorphism(-140,500)(0,-1)[\phantom{Y_2}`G(A) `]{450}1l
\putmorphism(-120,500)(0,-1)[\phantom{Y_2}` `\alpha_0(A)]{450}0l
\put(650,50){\fbox{$\alpha_0^v$}}
\putmorphism(-130,-400)(1,0)[G(\tilde A)` `G(g)]{470}1a
\putmorphism(-140,50)(0,-1)[\phantom{Y_2}``]{450}1l
\putmorphism(-120,50)(0,-1)[\phantom{Y_2}``G(u)]{450}0l
\putmorphism(450,50)(0,-1)[\phantom{Y_2}`G(\tilde B)`G(v)]{450}1r
\putmorphism(450,500)(0,-1)[\phantom{Y_2}`G(B) `\alpha_0(B)]{450}1r
\put(40,280){\fbox{$(\alpha_0)_f$}}
\putmorphism(-150,50)(1,0)[\phantom{(B, \tilde A)}``G(f)]{500}1a
\putmorphism(1030,50)(0,-1)[\phantom{(B, A')}`G(\tilde B).`]{450}1l
\putmorphism(1050,50)(0,-1)[``\alpha_0(\tilde B)]{450}0l
\putmorphism(1030,500)(0,-1)[F(B)`F(\tilde B)`]{450}1l
\putmorphism(1050,500)(0,-1)[``F(v)]{450}0l
\putmorphism(430,-400)(1,0)[\phantom{(B, \tilde A)}``=]{470}1b
\put(0,-170){\fbox{$G(a)$}}
\efig}
$$
\end{itemize}

A {\em vertical oplax transformation} is defined by using the opposite direction of the 2-cells $\alpha_0^u$ in item 3 
and accommodating the axioms correspondingly. 

A vertical lax transformation is called a {\em vertical pseudonatural transformation} if the globular 2-cells $\alpha_0^u$ are horizontally invertible for all 1v-cells $u$ in $\Aa$. 

It is called a {\em vertical strict transformation} if all the globular 2-cells $\alpha_0^u$ are identities. 
A vertical strict transformation is said to be 
{\em invertible} if the 2-cells $(\alpha_0)_f$ for all 1h-cells $f$ in $\Aa$ are vertically invertible (this includes the condition that 
the 1v-cells $\alpha_0(A)$ are invertible for all $A\in\Aa$). 
\end{defn}

\subsection{Companions, liftings and modifications}

Recall that a {\em companion} for a 1v-cell $u:A\to B$ is a 1h-cell  $\hat u:A\to B$ together with two 2-cells $\Epsilon$ and $\eta$ 
as below satisfying $[\eta\vert\Epsilon]=\Id_{\hat u}$ and $\frac{\eta}{\Epsilon}=\Id_{u}$. On the other hand, 
a {\em conjoint} for a 1v-cell $u:A\to B$ is a 1h-cell  $\check{u}:B\to A$ together with two 2-cells $\Epsilon^*$ and $\eta^*$ 
as below satisfying $[\Epsilon^*\vert\eta^*]=\Id_{\check{u}}$ and $\frac{\eta^*}{\Epsilon^*}=\Id^{u}$.
$$ 
\scalebox{0.9}{
\bfig
\putmorphism(-150,170)(1,0)[A`B`\hat{u}]{400}1a
\putmorphism(-150,-160)(1,0)[B`B`=]{400}1a
\putmorphism(-150,170)(0,-1)[\phantom{Y_2}``u]{330}1l
\putmorphism(250,170)(0,-1)[\phantom{Y_2}``=]{330}1r
\put(0,0){\fbox{$\Epsilon$}}
\efig}
\qquad%\qquad
\scalebox{0.9}{
\bfig
\putmorphism(-150,170)(1,0)[A`A`=]{400}1a
\putmorphism(-150,-160)(1,0)[A`B`\hat{u}]{400}1b
\putmorphism(-150,170)(0,-1)[\phantom{Y_2}``=]{330}1l
\putmorphism(250,170)(0,-1)[\phantom{Y_2}``u]{330}1r
\put(0,10){\fbox{$\eta$}}
\efig}
\qquad\qquad
\scalebox{0.9}{
\bfig
\putmorphism(-150,170)(1,0)[A`A`=]{400}1a
\putmorphism(-150,-160)(1,0)[B`A`\check{u}]{400}1b
\putmorphism(-150,170)(0,-1)[\phantom{Y_2}``u]{330}1l
\putmorphism(250,170)(0,-1)[\phantom{Y_2}``=]{330}1r
\put(-10,-10){\fbox{$\eta^*$}}
\efig}
\qquad%\qquad
\scalebox{0.9}{
\bfig
\putmorphism(-150,170)(1,0)[B`A`\check{u}]{400}1a
\putmorphism(-150,-160)(1,0)[B`B`=]{400}1b
\putmorphism(-150,170)(0,-1)[\phantom{Y_2}``=]{330}1l
\putmorphism(250,170)(0,-1)[\phantom{Y_2}``u]{330}1r
\put(-10,-10){\fbox{$\Epsilon^*$}}
\efig}
$$
Companions are unique up to a unique globular isomorphism. For more properties of companions and conjoints 
see \cite[Section 1.2]{GP:Adj}, \cite[Section 3]{Shul}. We only recall the following few fact about them. 

\begin{cor} \colabel{conj}
Let $u$ be an invertible 1v-cell. If $u^{-1}$ has a companion $\widehat{(u^{-1})}$, then $u$ has a conjoint given by $\check u=\widehat{(u^{-1})}$ with $\Epsilon^*_u=\frac{\Epsilon_{u^{-1}}}{\Id^u}$ and $\eta^*_u=\frac{\Id^u}{\eta_{u^{-1}}}$. 
\end{cor}

We will need existence of companions and conjoints only for 1v-components of invertible vertical strict transformations.  
In this case by \coref{conj} it is sufficient to require the existence of companions. 
The announced lifting properties of vertical to horizontal structures relies on the following result.

\medskip

\begin{prop} \cite[Proposition 2.14]{Fem3} \prlabel{lifting 1v to equiv}
Let $\alpha_0:F\Rightarrow G$ be a vertical pseudonatural transformation between pseudodouble functors acting between double categories 
$\Dd\to\Ee$. Suppose that the 2-cell components $(\alpha_0)_f$ of $\alpha_0$ are vertically invertible for every 1h-cell $f:A\to B$. 
\begin{enumerate}
\item The following data define a horizontal pseudonatural transformation $\alpha_1:F\Rightarrow G$:
\begin{enumerate} [a)]
\item for all 1v-cell components $\alpha_0(A)$ of $\alpha_0$ a fixed choice of companions (and conjoints) in $\Ee$, 
(we denote their companions by $\alpha_1(A)$ for every 0-cell $A$ of $\Dd$, the corresponding 2-cells by $\Epsilon^\alpha_A$ and 
$\eta^\alpha_A$, and also by $\Epsilon^*_A$ and $\eta^*_A$ the 2-cells related to conjoints of the inverse of $\alpha_0(A)$); 
\item the 2-cells 
$$\delta_{\alpha_1,f}=\quad
\scalebox{0.86}{
\bfig
 \putmorphism(-150,250)(1,0)[F(A)`\phantom{F(A)} `=]{500}1a
\put(-100,30){\fbox{$\eta^\alpha_A$}}
\putmorphism(380,250)(0,-1)[\phantom{Y_2}` `]{450}1l
\putmorphism(410,250)(0,-1)[\phantom{Y_2}` `\alpha_0(A)]{450}0l
\putmorphism(950,250)(0,-1)[\phantom{Y_2}` `]{450}1r
\putmorphism(930,250)(0,-1)[\phantom{Y_2}` `\alpha_0(B)]{450}0r
\putmorphism(350,250)(1,0)[F(A)`F(B)`F(f)]{600}1a
 \putmorphism(950,250)(1,0)[\phantom{F(A)}`G(B) `\alpha_1(B)]{600}1a
 \putmorphism(470,-200)(1,0)[`G(B)`G(f)]{500}1b
 \putmorphism(1060,-200)(1,0)[`G(B)`=]{500}1b
\putmorphism(1570,250)(0,-1)[\phantom{Y_2}``=]{450}1r
\put(530,10){\fbox{$(\alpha_0)_f$}}
\put(1280,30){\fbox{$\Epsilon^\alpha_B$}}
\putmorphism(-150,-200)(1,0)[F(A)`G(A) `\alpha_1(A)]{520}1a
\putmorphism(-150,250)(0,-1)[``=]{450}1l
\efig}
$$
in $\Ee$ for every 1h-cell $f:A\to B$ in $\Dd$;
\item the 2-cells 
$$(\alpha_1)^u=\quad
\scalebox{0.86}{
\bfig
 \putmorphism(-40,50)(1,0)[``=]{320}1b

\put(550,30){\fbox{$\delta_{\alpha_0,u}$}}
\putmorphism(-150,-400)(1,0)[F(\tilde A)`G(\tilde A) `\alpha_1(\tilde A)]{520}1a
\putmorphism(-150,50)(0,-1)[F(A')``=]{450}1l
\putmorphism(380,500)(0,-1)[\phantom{Y_2}` `F(u)]{450}1l
\putmorphism(950,500)(0,-1)[\phantom{Y_2}`G(A) `\alpha_0(A)]{450}1r
\putmorphism(380,50)(0,-1)[F(\tilde A)` `\alpha_0(\tilde A)]{450}1r
\putmorphism(350,500)(1,0)[F(A)`F(A)`=]{600}1a
 \putmorphism(950,500)(1,0)[\phantom{F(A)}`G(A) `\alpha_1(A)]{650}1a
 \putmorphism(1060,50)(1,0)[`G(A)`=]{500}1b

\putmorphism(1570,500)(0,-1)[\phantom{Y_2}``=]{450}1r
\put(1250,260){\fbox{$\Epsilon^\alpha_A$}}
\putmorphism(480,-400)(1,0)[`G(\tilde A) `=]{500}1a
\putmorphism(950,50)(0,-1)[\phantom{Y_2}``G(u)]{450}1r
\put(0,-160){\fbox{$\eta^\alpha_{\tilde A}$}}
\efig}
$$
in $\Ee$ for every 1v-cell $u:A\to \tilde A$ in $\Dd$;
\item the inverses of the 2-cells $\delta_{\alpha_1,f}$ are given by 
$$\delta_{\alpha_1,f}^{-1}=\quad
\scalebox{0.86}{
\bfig
 \putmorphism(-150,250)(1,0)[F(A)`\phantom{F(A)} `\alpha_1(A)]{500}1a
\put(-120,30){\fbox{$\Epsilon^*_A$}}
\putmorphism(400,250)(0,-1)[\phantom{Y_2}` `]{450}1l
\putmorphism(430,250)(0,-1)[\phantom{Y_2}` `\alpha_0(A)^{-1}]{450}0l
\putmorphism(950,250)(0,-1)[\phantom{Y_2}` `]{450}1r
\putmorphism(930,250)(0,-1)[\phantom{Y_2}` `\alpha_0(B)^{-1}]{450}0r
\putmorphism(350,250)(1,0)[G(A)`G(B)`G(f)]{600}1a
 \putmorphism(950,250)(1,0)[\phantom{F(A)}`G(B) `=]{600}1a
 \putmorphism(470,-200)(1,0)[`F(B)`F(f)]{500}1b
 \putmorphism(1060,-200)(1,0)[`G(B).`\alpha_1(B)]{500}1b
\putmorphism(1570,250)(0,-1)[\phantom{Y_2}``=]{450}1r
\put(530,10){\fbox{$(\alpha_0)_f^{-1}$}}
\put(1300,30){\fbox{$\eta^*_B$}} 
\putmorphism(-150,-200)(1,0)[F(A)`F(A) `=]{520}1a
\putmorphism(-150,250)(0,-1)[``=]{450}1l
\efig}
$$
\end{enumerate}
\item The 1h-cells $\alpha_1(A)$ are adjoint equivalence 1-cells in $\HH(\Ee)$. 
\item The data from 1. define a horizontal equivalence $\alpha_1:F\Rightarrow G:\Dd\to\Ee$ and a  
pseudonatural equivalence $\HH(\alpha_1):\HH(F)\Rightarrow \HH(G):\HH(\Dd)\to\HH(\Ee)$. 
\end{enumerate}
\end{prop}

\begin{defn} 
Those 1v-cells that have companions we will call {\em companion-liftable} (or shortly {\em liftable}) 1v-cells. \\
Those vertical transformations all of whose 1v-cell components are liftable we will call {\em companion-liftable} (or shortly {\em liftable}) 
vertical transformations. \\
A horizontal pseudonatural equivalence $\alpha_1$ obtained in \prref{lifting 1v to equiv} from an invertible liftable vertical strict transformation $\alpha_0$ we will call a {\em companion-lift} of $\alpha_0$. 
\end{defn}

We recall the definition of a modification, our choice in it that the horizontal transformations are oplax and the vertical transformations 
are lax is fully arbitrary. 

\begin{defn} \delabel{modif-hv}
A modification $\Theta$ between two horizontal oplax transformations $\alpha$ and $\beta$ and two vertical lax transformations $\alpha_0$ and $\beta_0$ 
depicted below on the left, where the lax double functors $F, G, F\s', G'$ act between $\Aa\to\Bb$, is given 
by a collection of 2-cells in $\Bb$ depicted below on the right:
$$%%\begin{equation} \eqlabel{modification cells}
\scalebox{0.86}{
\bfig
\putmorphism(-150,50)(1,0)[F` G`\alpha]{400}1a
\putmorphism(-150,-270)(1,0)[F'`G' `\beta]{400}1b
\putmorphism(-170,50)(0,-1)[\phantom{Y_2}``\alpha_0]{320}1l
\putmorphism(250,50)(0,-1)[\phantom{Y_2}``\beta_0]{320}1r
\put(-30,-140){\fbox{$\Theta$}}
\efig}
\qquad\qquad
\scalebox{0.86}{
\bfig
\putmorphism(-180,50)(1,0)[F(A)` G(A)`\alpha(A)]{550}1a
\putmorphism(-180,-270)(1,0)[F\s'(A)`G'(A) `\beta(A)]{550}1b
\putmorphism(-170,50)(0,-1)[\phantom{Y_2}``\alpha_0(A)]{320}1l
\putmorphism(350,50)(0,-1)[\phantom{Y_2}``\beta_0(A)]{320}1r
\put(0,-140){\fbox{$\Theta_A$}}
\efig}
$$%\end{equation}
which satisfy the following rules: 

\medskip

\noindent {\em \axiom{m.ho-vl.-1}}  for every 1h-cell $f$, we have  
$$
\scalebox{0.86}{
\bfig
\putmorphism(-150,500)(1,0)[F(A)`F(B)`F(f)]{600}1a
 \putmorphism(450,500)(1,0)[\phantom{F(A)}`G(B) `\alpha(B)]{620}1a

 \putmorphism(-150,50)(1,0)[F\s'(A)`F\s'(B)`F\s'(f)]{600}1a
 \putmorphism(450,50)(1,0)[\phantom{F(A)}`G'(B) `\beta(B)]{620}1a

\putmorphism(-180,500)(0,-1)[\phantom{Y_2}``\alpha_0(A)]{450}1l
\putmorphism(450,500)(0,-1)[\phantom{Y_2}``]{450}1r
\putmorphism(300,500)(0,-1)[\phantom{Y_2}``\alpha_0(B)]{450}0r
\putmorphism(1080,500)(0,-1)[\phantom{Y_2}``\beta_0(B)]{450}1r
\put(0,280){\fbox{$(\alpha_0)_f$}}
\put(670,280){\fbox{$\Theta_B$}}

\putmorphism(-150,-400)(1,0)[F\s'(A)`G'(A) `\beta(A)]{600}1a
\putmorphism(510,-400)(1,0)[\phantom{Y_2}`G'(B) `G'(f)]{580}1a

\putmorphism(-180,50)(0,-1)[\phantom{Y_2}``=]{450}1l
\putmorphism(1080,50)(0,-1)[\phantom{Y_3}``=]{450}1r
\put(320,-180){\fbox{$\delta_{\beta,f}$}}

\efig}
\quad=\quad
\scalebox{0.86}{
\bfig
\putmorphism(-150,500)(1,0)[F(A)`F(B)`F(f)]{600}1a
 \putmorphism(450,500)(1,0)[\phantom{F(A)}`G(B) `\alpha(B)]{620}1a
\putmorphism(-150,50)(1,0)[F(A)`G(A) `\alpha(A)]{600}1a
\putmorphism(510,50)(1,0)[\phantom{Y_2}`G(B) `G(f)]{580}1a

\putmorphism(-180,500)(0,-1)[\phantom{Y_2}``=]{450}1r
\putmorphism(1080,500)(0,-1)[\phantom{Y_2}``=]{450}1r
\put(350,280){\fbox{$\delta_{\alpha,f}$}}

\putmorphism(-180,50)(0,-1)[\phantom{Y_2}``\alpha_0(A)]{450}1l
\putmorphism(1080,50)(0,-1)[\phantom{Y_3}``\beta_0(B)]{450}1r
\put(20,-180){\fbox{$\Theta_A$}}
\put(670,-180){\fbox{$(\beta_0)_f$}}

\putmorphism(450,50)(0,-1)[\phantom{Y_2}``]{450}1r
\putmorphism(300,50)(0,-1)[\phantom{Y_2}``\beta_0(A)]{450}0r

\putmorphism(-150,-400)(1,0)[F\s'(A)`G'(A) `\beta(A)]{600}1a
\putmorphism(510,-400)(1,0)[\phantom{Y_2}`G'(B) `G'(f)]{580}1a
\efig}
$$
and

\noindent {\em \axiom{m.ho-vl.-2}}   for every 1v-cell $u$, we have
$$
\scalebox{0.86}{
\bfig
 \putmorphism(-150,500)(1,0)[F(A)`F(A) `=]{600}1a
 \putmorphism(550,500)(1,0)[` `\alpha(A)]{400}1a
\putmorphism(-180,500)(0,-1)[\phantom{Y_2}`F\s'(A) `\alpha_0(A)]{450}1l
\put(30,50){\fbox{$\alpha_0^u$}}
\putmorphism(-150,-400)(1,0)[F\s'(\tilde A)` `=]{480}1a
\putmorphism(-180,50)(0,-1)[\phantom{Y_2}``F\s'(u)]{450}1l
\putmorphism(450,50)(0,-1)[\phantom{Y_2}`F\s'(\tilde A)` \alpha_0(\tilde A)]{450}1l
\putmorphism(450,500)(0,-1)[\phantom{Y_2}`F(\tilde A) `F(u)]{450}1l
\put(660,280){\fbox{$\alpha^u$}}
\putmorphism(450,50)(1,0)[\phantom{(B, \tilde A)}``\alpha(\tilde A)]{500}1a
\putmorphism(1070,50)(0,-1)[\phantom{(B, A')}`G'(\tilde A)`\beta_0(\tilde A)]{450}1r
\putmorphism(1070,500)(0,-1)[G(A)`G(\tilde A)`G(u)]{450}1r
\putmorphism(450,-400)(1,0)[\phantom{(B, \tilde A)}``\beta(\tilde A)]{500}1a
\put(640,-170){\fbox{$ \Theta_{\tilde A}$ } } % ???
\efig}\quad=\quad
\scalebox{0.86}{
\bfig
 \putmorphism(-150,500)(1,0)[F(A)`G(A) `\alpha(A)]{600}1a
 \putmorphism(450,500)(1,0)[\phantom{(B,A)}` `=]{460}1a
\putmorphism(-180,500)(0,-1)[\phantom{Y_2}`F\s'(A) `\alpha_0(A)]{450}1l
\put(650,50){\fbox{$\beta_0^u$}}
\putmorphism(-180,-400)(1,0)[F\s'(\tilde A)` `\beta(\tilde A)]{500}1a
\putmorphism(-180,50)(0,-1)[\phantom{Y_2}``F\s'(u)]{450}1l
\putmorphism(450,50)(0,-1)[\phantom{Y_2}`G'(\tilde A)`G'(u)]{450}1r
\putmorphism(450,500)(0,-1)[\phantom{Y_2}`G'(A) `\beta_0(A)]{450}1r
\put(0,260){\fbox{$\Theta_A$}}
\putmorphism(-180,50)(1,0)[\phantom{(B, \tilde A)}``\beta(A)]{500}1a
\putmorphism(1030,50)(0,-1)[\phantom{(B, A')}` G'(\tilde A). ` \beta_0(\tilde A)]{450}1r
\putmorphism(1030,500)(0,-1)[G(A)`G(\tilde A)` G(u)]{450}1r
\putmorphism(430,-400)(1,0)[\phantom{(B, \tilde A)}``=]{480}1b
\put(70,-170){\fbox{$\beta^u$}}
\efig}
$$
\end{defn}

Horizontal and vertical restrictions of modifications in \deref{modif-hv} yield:
\begin{itemize} 
\item {\em horizontal modifications} or {\em modifications between horizontal (oplax) transformations} are given by families of (horizontally globular) 2-cells 
\begin{equation} \eqlabel{m-hor}
\scalebox{0.86}{
\bfig
\putmorphism(-180,50)(1,0)[F(A)` G(A)`\alpha(A)]{550}1a
\putmorphism(-180,-270)(1,0)[F\s'(A)`G'(A) `\beta(A)]{550}1b
\putmorphism(-170,50)(0,-1)[\phantom{Y_2}``=]{320}1l
\putmorphism(350,50)(0,-1)[\phantom{Y_2}``=]{320}1r
\put(0,-140){\fbox{$\Theta_A$}}
\efig}
\end{equation} 
and axioms \axiom{m.ho.-1} and \axiom{m.ho.-2} obtained from \axiomref{m.ho-vl.-1} and \axiomref{m.ho-vl.-2} by ignoring the 2-cells 
$(\alpha_0)_f, (\beta_0)_f, \alpha_0^u$ and $\beta_0^u$, and 
\item {\em vertical modifications} or {\em modifications between vertical (lax) transformations} are given by families of (vertically globular) 2-cells 
\begin{equation} \eqlabel{m-vert}
\scalebox{0.86}{
\bfig
\putmorphism(-180,50)(1,0)[F(A)` G(A)`=]{550}1a
\putmorphism(-180,-270)(1,0)[F\s'(A)`G'(A) `=]{550}1b
\putmorphism(-170,50)(0,-1)[\phantom{Y_2}``\alpha_0(A)]{320}1l
\putmorphism(350,50)(0,-1)[\phantom{Y_2}``\beta_0(A)]{320}1r
\put(0,-140){\fbox{$\Theta_A$}}
\efig}
\end{equation} 
and axioms \axiom{m.vl.-1} and \axiom{m.vl.-2} obtained from \axiomref{m.ho-vl.-1} and \axiomref{m.ho-vl.-2} by ignoring the 2-cells 
$\delta_{\alpha,f}, \delta_{\beta,f}, \alpha^u$ and $\beta^u$.  
\end{itemize}

\begin{prop} \prlabel{essence}
Let %$\Dd$ be a monoidal double category and 
$\omega$ be the identity vertical modification between two vertical composites of vertical strict transformations 
$$\omega: \,\, \threefrac{\alpha_1}{...}{\alpha_k} \,\,\, \Rrightarrow \,\,\, \threefrac{\beta_1}{...}{\beta_l} \,\,$$ 
which act between lax double functors $F\Rightarrow G:\Bb\to\Dd$ between double categories,  
%globular 2-cell in a double category $\Dd$ with the following characteristics: 
%\begin{itemize}
%\item its component 1v-cells are 1v-cell components of composable vertical strict transformations and they all have companions, and
%%\item its component 1v-cells may be tensor products of other 1v-cells, and 
so that all 1v-cell components $\alpha_1(A),...\alpha_k(A)$ and $\beta_1(A), ..., \beta_l(A)$ for 0-cells $A$ in $\Bb$, have companions in $\Dd$. 
%\end{itemize}
Then: 
\begin{enumerate}
\item $\omega$ induces an invertible  %horizontally globular 2-cell $\hat\omega$, which defines a (horizontal) modification
horizontal modification $\hat\omega$ between the two (vertical compositions of the) induced horizontal natural transformations; 
\item the assignment between 2-cell components of $\omega$ and $\hat\omega$ is invertible; 
%\textcolor{rojo}{
%\item if the vertical strict transformations $\alpha_1,...\alpha_k$ and $\beta_1, ..., \beta_l$ are invertible, then 
%the modification $\hat\omega$ is invertible;}
\item if $\omega_1,...,\omega_m$ are %vertically globular 2-cells with the above characteristics, 
vertical modifications with the above characteristics, then any sensible equation formed by their 2-cell components  
$\hat\omega_1(A),...,\hat\omega_m(A)$, for any 0-cell $A$ in $\Bb$, %, which are components of the induced modifications, 
holds true. 
\end{enumerate}
\end{prop}

We will also need the following variation of the above claim. In it we %continue to 
treat $\omega$ %both as a modification and its component 2-cells. 
as component 2-cells of the modification in question.

\begin{prop} \prlabel{omega*}
Let $\omega$ be as in \prref{essence} and assume moreover that at least two of 1v-cells in one of its vertical edges are non-trivial. Then: 
%and at least one 1v-cell in the other edge is non-trivial. Then:
\begin{enumerate}
\item $\omega$ induces a 2-cell $\omega^*$ defining a modification in the sense of \deref{modif-hv};
\item there is a one-to-one correspondence between 2-cells $\omega^*$ and $\hat\omega$ from \prref{essence}, so that $\omega^*$ induces a modification if and only if so does $\hat\omega$;
\item if $\omega^*_1,...,\omega^*_m$ are 2-cells induced by vertically globular 2-cells $\omega_1,...,\omega_m$ with the characteristics as $\omega$ of this proposition, then any sensible equation formed by the 2-cells $\omega^*_1,...,\omega^*_m$ holds true. 
%$\hat\omega_1,...,\hat\omega_k$, which are components of the induced modifications, 
\end{enumerate}
\end{prop}

Similarly as in \prref{omega*}, a general vertically globular 2-cell $\omega_1$ in which $v,v'$ (as on the left below) are liftable 
semi-lifts to a square-formed 2-cell $\omega^*_1$. Also, a general square $\omega_2$ in which $u,u'$ (as in the middle below) 
are liftable lifts to a horizontally globular 2-cell $\hat\omega_2$.
\begin{equation} \eqlabel{omega-gen}
\scalebox{0.9}{
\bfig
 \putmorphism(-90,500)(1,0)[A` A`=]{440}1a
\putmorphism(350,500)(0,-1)[` B'`v']{400}1r
\putmorphism(-90,-300)(1,0)[C`C `=]{440}1a
\putmorphism(350,110)(0,-1)[\phantom{Y_2}``u']{400}1r
\putmorphism(-100,110)(0,-1)[\phantom{Y_2}``v]{400}1l
\putmorphism(-100,500)(0,-1)[` B`u]{400}1l
\put(60,70){\fbox{$\omega_1$}}
\efig}
\quad\mapsto\quad 
\scalebox{0.9}{
\bfig
\putmorphism(-150,170)(1,0)[A`B'`\hat v']{460}1a
\putmorphism(-150,-180)(1,0)[B`C`\hat v]{460}1a
\putmorphism(-150,170)(0,-1)[\phantom{Y_2}``u]{350}1l
\putmorphism(310,170)(0,-1)[\phantom{Y_2}``u']{350}1r
\put(0,0){\fbox{$\omega^*_1$}}
\efig}\hspace{0,2cm};
\qquad
\scalebox{0.9}{
\bfig
\putmorphism(-150,170)(1,0)[A`A'`f]{460}1a
\putmorphism(-150,-180)(1,0)[B`B'`g]{460}1a
\putmorphism(-150,170)(0,-1)[\phantom{Y_2}``u]{350}1l
\putmorphism(310,170)(0,-1)[\phantom{Y_2}``u']{350}1r
\put(0,0){\fbox{$\omega_2$}}
\efig}
\quad\mapsto\quad 
\scalebox{0.9}{
\bfig
 \putmorphism(-150,210)(1,0)[A`A'`f]{450}1a
\putmorphism(-160,200)(0,-1)[\phantom{Y_2}`A `=]{380}1l
 \putmorphism(300,210)(1,0)[\phantom{Y}`B'  `\hat u']{450}1a
\putmorphism(730,200)(0,-1)[\phantom{Y_2}`B'. `=]{380}1r
\putmorphism(-150,-200)(1,0)[`B`\hat u]{450}1a
\putmorphism(350,-200)(1,0)[``g]{330}1a
\put(220,0){\fbox{$\hat\omega_2$}}
\efig}
\end{equation}
Moreover, all the laws that such 2-cells $\omega_1, \omega_2$ obey are preserved in the analogous form by the lifting operation. 
This is so, as the lifting is done by inserting the structure 2-cells $\Epsilon$ (on the extreme right) and $\eta$ (on the extreme left), so that a lift of any kind of composition of 2-cells $a_1,..,a_k$ is the same kind of composition of the lifted 2-cells  $\hat a_1,..,\hat a_k$ (at the joins one inserts the identity 2-cells $\frac{\eta}{\Epsilon}$).

\section{Horizontally premonoidal double categories}

In \cite{Fem3} we introduced premonoidal double categories. As we discussed in Section 3.3 thereof, there are (at least) two ways 
to define them. In that reference we chose the vertical orientation of their structural transformations, in the style of \cite{Shul}, 
so that lifting vertical structures one obtains the corresponding horizontal ones and thus covers a bicategorical situation. However, 
that construction would not cover {\em all} premonoidal bicategories. Another way is to let the structural transformations be horizontal 
and thus enlarge the covered bicategorical premonoidal setting. We take this approach here. %by including the vertical cells. 

\subsection{Binoidal double categories and central cells}  \sslabel{binoidal}  % and center double category}

We will need the notion of central horizontal 1-cells (in \cite[Section 3.1]{Fem3} we introduced also the centrality structures on 1v- and 2-cells). 

\begin{defn}
We say that a double category $\Bb$ is {\em binoidal} if for all 0-cells $A,B\in\Bb$ there 
%are both horizontal and vertical pseudonatural transformations  
are pseudodouble functors $A\ltimes -$ and $-\rtimes B$ acting $\Bb\to\Bb$ and such that $A\ltimes B=
A\rtimes B=:A\bowtie B$. 
\end{defn}

\begin{defn}
Assume $\Bb$ is a binoidal double category. 
\begin{itemize}
\item A 1h-cell $f:A\to A'$ in $\Bb$ is said to be {\em left central} if there is a horizontal pseudonatural transformation 
$f\ltimes-:A\ltimes -\to A'\ltimes-$ such that $f\ltimes B=f\rtimes B$ for all $B\in\Bb$. % induces a pseudonatural transformation on the underlying horizontal 2-category $\HH(\Bb)$. 
Likewise, a 1h-cell $f:A\to A'$ in $\Bb$ is said to be {\em right central} if there is a horizontal pseudonatural transformation 
$-\rtimes f:-\ltimes A\to -\rtimes A'$ such that $B\rtimes f=B\ltimes f$ for all $B\in\Bb$. \\ % induces a pseudonatural transformation on $\HH(\Bb)$. \\
A 1h-cell is said to be {\em central} if it is both left and right central. 
\item A 1v-cell $v:A\to\tilde A$ in $\Bb$ is said to be {\em left central} if there is a vertical pseudonatural transformation 
$v\ltimes-:A\ltimes -\to \tilde A\ltimes-$ such that $v\ltimes B=v\rtimes B$ for all $B\in\Bb$. % induces a pseudonatural transformation on the underlying vertical 2-category $\V(\Bb)$. 
Likewise, a 1v-cell $v:A\to\tilde A$ in $\Bb$ is said to be {\em right central} if there is a vertical pseudonatural transformation 
$-\rtimes v:-\rtimes A\to -\rtimes\tilde A$ such that $B\ltimes v=B\rtimes v$ for all $B\in\Bb$. \\ % induces a pseudonatural transformation on $\V(\Bb)$. \\ 
A 1v-cell is said to be {\em central} if it is both left and right central.
\end{itemize}
\end{defn}

For reader's convenience we write down the 2-cell components of the pseudonatural transformations in play. 
For a horizontal pseudonatural transformation $f\ltimes-:A\ltimes -\to A'\ltimes-$, a 1h-cell $g:B\to B'$  and 
a 1v-cell $v:A\to\tilde A$ the 2-cell components $f\ltimes-\vert_g$ and $f\ltimes-\vert_v$ of the oplax transformation structure of 
$f\ltimes-$ have the form of the left diagrams below (for the lax transformation structure of $f\ltimes-$ the first 2-cell 
component is different: it is a 2-cell obtained by reading the same upper left diagram from bottom to top). 
Likewise, the 2-cell components $-\rtimes f\vert_g$ and $-\rtimes f\vert_v$ of the oplax transformation structure of 
$-\rtimes f:-\rtimes A\to -\rtimes A'$  has the form of the right diagrams below (and the differing 2-cell component for the lax structure 
is obtained by reading the upper right diagram from bottom to top). 
\begin{equation} \eqlabel{f lr}
\scalebox{0.86}{
\bfig
 \putmorphism(-170,500)(1,0)[A\ltimes B`A\ltimes B' `A\ltimes g]{620}1a
 \putmorphism(450,500)(1,0)[\phantom{A\ot B}`A'\ltimes B' `f\ltimes B']{680}1a
 \putmorphism(-150,50)(1,0)[A\ltimes B`A'\ltimes B `f\ltimes B]{600}1a
 \putmorphism(450,50)(1,0)[\phantom{A\ot B}`A'\ltimes B' `A'\ltimes g]{680}1a

\putmorphism(-180,500)(0,-1)[\phantom{Y_2}``=]{450}1r
\putmorphism(1100,500)(0,-1)[\phantom{Y_2}``=]{450}1r
\put(350,260){\fbox{$f\ltimes-\vert_g$}}
\efig}
\qquad\qquad
\scalebox{0.86}{
\bfig
 \putmorphism(-170,500)(1,0)[B\rtimes A`B'\rtimes A `g\rtimes A]{620}1a
 \putmorphism(450,500)(1,0)[\phantom{A\ot B}`B'\rtimes A' `B'\rtimes f]{680}1a
 \putmorphism(-150,50)(1,0)[B\rtimes A`B\rtimes A' `B\rtimes f]{600}1a
 \putmorphism(450,50)(1,0)[\phantom{A\ot B}`B'\rtimes A' `g\rtimes A']{680}1a

\putmorphism(-180,500)(0,-1)[\phantom{Y_2}``=]{450}1r
\putmorphism(1100,500)(0,-1)[\phantom{Y_2}``=]{450}1r
\put(350,260){\fbox{$-\rtimes f\vert_g$}}
\efig}
\end{equation}

$$
\scalebox{0.86}{
\bfig
 \putmorphism(-150,50)(1,0)[A\ltimes B`A'\ltimes B `f\ltimes B]{640}1a
\putmorphism(-150,-400)(1,0)[A\ltimes\tilde  B`A'\ltimes\tilde  B `f\ltimes \tilde B]{640}1a
\putmorphism(-160,50)(0,-1)[\phantom{Y_2}``A\ltimes v]{450}1l
\putmorphism(450,50)(0,-1)[\phantom{Y_2}``A'\ltimes v]{450}1r
\put(-40,-170){\fbox{$f\ltimes-\vert_v$}}
\efig}
\qquad\qquad\qquad
\scalebox{0.86}{
\bfig
 \putmorphism(-150,50)(1,0)[B\rtimes A`B\rtimes A' `B\rtimes f]{640}1a
\putmorphism(-150,-400)(1,0)[\tilde B\rtimes A ` \tilde B\rtimes A' ` \tilde B\rtimes f]{640}1a
\putmorphism(-160,50)(0,-1)[\phantom{Y_2}``v\rtimes A]{450}1l
\putmorphism(450,50)(0,-1)[\phantom{Y_2}``v\rtimes A']{450}1r
\put(-40,-170){\fbox{$-\rtimes f\vert_v$}}
\efig}
$$

\subsection{Horizontally premonoidal double categories}

The title of this subsection is intended to stress the difference between the notion of a premonoidal double category that will be introduced below, and that from \cite{Fem3}. In this paper all premonoidal double categories will be meant in the ``horizontal'' sense (except when we compare them to the vertical ones). When we drop the adjective, the horizontal sense is meant. 

\begin{defn} \delabel{premon} %\delabel{unital} 
Let $\Bb$ be a binoidal double category. 
We say that $\Bb$ is a (horizontally) {\em premonoidal} double category if there exist: 
\begin{enumerate}
\item a unit object $I$, 
\item pseudonatural equivalences $\lambda: I\ltimes -\Rightarrow\Id$ and $\rho: -\rtimes I\Rightarrow\Id$ and 
$$\alpha_{-,B,C}: (-\rtimes B)\rtimes C\Rightarrow -\rtimes(B\bowtie C)$$
$$\alpha_{A,-,C}: (A\ltimes -)\rtimes C\Rightarrow A\ltimes(-\rtimes C)$$
$$\alpha_{A,B,-}: (A\bowtie B)\ltimes -\Rightarrow A\ltimes(B\ltimes -)$$
for every $A,B,C\in\Bb$, so that  
 \begin{enumerate}[i)]
 \item all 1h-cell components of the above five equivalences are central, 
 \item 1h-cell components of the three equivalences $\alpha$ coincide, 
%$\alpha_{A,B,C}:=\alpha^1_A=\alpha^2_B=\alpha^3_C$ is their common 
\end{enumerate}
\item there are ten horizontal modifications 
$$
\scalebox{0.8}{
\bfig
 \putmorphism(-150,400)(1,0)[((-B)C)D`(-(BC))D `\alpha_{-,B,C}D]{1100}1a
 \putmorphism(1050,400)(1,0)[\phantom{A\ot B}`-((BC)D)`\alpha_{-,BC,D}]{1080}1a
 \putmorphism(2220,400)(1,0)[\phantom{A\ot B}` -(B(CD)) `-\rtimes\alpha_{B,C,D}]{1120}1a

 \putmorphism(-170,50)(1,0)[((-B)C)D`(-B)(CD) `\alpha_{-B,C,D}]{1620}1b
 \putmorphism(1560,50)(1,0)[\phantom{A\ot B}`-(B(CD)) `\alpha_{-,B,CD}]{1800}1b %(M(NP))E

\putmorphism(-180,400)(0,-1)[\phantom{Y_2}``=]{350}1r
\putmorphism(3290,400)(0,-1)[\phantom{Y_2}``=]{350}1l
\put(1380,230){\fbox{$p_{-,B,C,D}$}}
\efig}
$$

%2
$$
\scalebox{0.8}{
\bfig
 \putmorphism(-150,400)(1,0)[((A-)C)D`(A(-C))D `\alpha_{A,-,C}D]{1100}1a
 \putmorphism(1050,400)(1,0)[\phantom{A\ot B}`A((-C)D)`\alpha_{A,-C,D}]{1080}1a
 \putmorphism(2220,400)(1,0)[\phantom{A\ot B}` A(-(CD)) `A\alpha_{-,C,D}]{1120}1a

 \putmorphism(-170,50)(1,0)[((A-)C)D`(A-)(CD) `\alpha_{A-,C,D}]{1620}1b
 \putmorphism(1560,50)(1,0)[\phantom{A\ot B}`A(-(CD)) `\alpha_{A,-,CD}]{1800}1b %(M(NP))E

\putmorphism(-180,400)(0,-1)[\phantom{Y_2}``=]{350}1r
\putmorphism(3290,400)(0,-1)[\phantom{Y_2}``=]{350}1l
\put(1380,230){\fbox{$p_{A,-,C,D}$}}
\efig}
$$

%3
$$
\scalebox{0.8}{
\bfig
 \putmorphism(-150,400)(1,0)[((AB)-)D`(A(B-))D `\alpha_{A,B,-}D]{1100}1a
 \putmorphism(1050,400)(1,0)[\phantom{A\ot B}`A((B-)D)`\alpha_{A,B-,D}]{1080}1a
 \putmorphism(2220,400)(1,0)[\phantom{A\ot B}` A(B(-D)) `A\alpha_{B,-,D}]{1120}1a

 \putmorphism(-170,50)(1,0)[((AB)-)D`(AB)(-D) `\alpha_{AB,-,D}]{1620}1b
 \putmorphism(1560,50)(1,0)[\phantom{A\ot B}`A(B(-D)) `\alpha_{A,B,-D}]{1800}1b %(M(NP))E

\putmorphism(-180,400)(0,-1)[\phantom{Y_2}``=]{350}1r
\putmorphism(3290,400)(0,-1)[\phantom{Y_2}``=]{350}1l
\put(1380,230){\fbox{$p_{A,B,-,D}$}}
\efig}
$$

%4
$$
\scalebox{0.8}{
\bfig
 \putmorphism(-150,400)(1,0)[((AB)C)-`(A(BC))- `\alpha_{A,B,C}\ltimes-]{1100}1a
 \putmorphism(1050,400)(1,0)[\phantom{A\ot B}`A((BC)-)`\alpha_{A,BC,-}]{1080}1a
 \putmorphism(2220,400)(1,0)[\phantom{A\ot B}` A(B(C-)) `A\alpha_{B,C,-}]{1120}1a

 \putmorphism(-170,50)(1,0)[((AB)C)-`(AB)(C-) `\alpha_{AB,C,-}]{1620}1b
 \putmorphism(1560,50)(1,0)[\phantom{A\ot B}`A(B(C-)) `\alpha_{A,B,C-}]{1800}1b %(M(NP))E

\putmorphism(-180,400)(0,-1)[\phantom{Y_2}``=]{350}1r
\putmorphism(3290,400)(0,-1)[\phantom{Y_2}``=]{350}1l
\put(1380,230){\fbox{$p_{A,B,C,-}$}}
\efig}
$$

% left action
$$%\begin{equation} \eqlabel{effects}
\scalebox{0.8}{
\bfig
 \putmorphism(-150,400)(1,0)[(IA)-`A- `\lambda_A\ltimes-]{1430}1a

 \putmorphism(-180,50)(1,0)[(IA)-`I(A-) `\alpha_{I,A,-}]{650}1b %I(ME)
 \putmorphism(460,50)(1,0)[\phantom{A\ot B}`A- `\lambda_{A-}]{780}1b

\putmorphism(-180,400)(0,-1)[\phantom{Y_2}``=]{350}1r
\putmorphism(1280,400)(0,-1)[\phantom{Y_2}``=]{350}1l
\put(540,220){\fbox{$l_{A,-}$}}
\efig}
\qquad
\scalebox{0.8}{
\bfig
 \putmorphism(-150,400)(1,0)[(I-)B`-B `\lambda_{-}B]{1430}1a

 \putmorphism(-180,50)(1,0)[(I-)B`I(-B) `\alpha_{I,-,B}]{650}1b %I(ME)
 \putmorphism(460,50)(1,0)[\phantom{A\ot B}`-B `\lambda_{-B}]{780}1b
 %\putmorphism(1260,50)(1,0)[\phantom{A\ot B}`ME`]{800}1a

\putmorphism(-180,400)(0,-1)[\phantom{Y_2}``=]{350}1r
\putmorphism(1280,400)(0,-1)[\phantom{Y_2}``=]{350}1l
\put(540,220){\fbox{$l_{-,B}$}}
\efig}
$$

%right action
$$%\begin{equation} \eqlabel{effects}
\scalebox{0.8}{
\bfig
 \putmorphism(-150,400)(1,0)[(-B)I `-B `\rho_{-B}]{1430}1a

 \putmorphism(-180,50)(1,0)[(-B)I `E(MI) `\alpha_{-,B,I}]{650}1b %I(ME)
 \putmorphism(460,50)(1,0)[\phantom{A\ot B}`-B `-\rtimes\rho_B]{780}1b

\putmorphism(-180,400)(0,-1)[\phantom{Y_2}``=]{350}1r
\putmorphism(1280,400)(0,-1)[\phantom{Y_2}``=]{350}1l
\put(540,220){\fbox{$r_{-,B}$}}
\efig}
\qquad
\scalebox{0.8}{
\bfig
 \putmorphism(-150,400)(1,0)[(A-)I `A- `\rho_{A-}]{1430}1a

 \putmorphism(-180,50)(1,0)[(A-)I `A(-I) `\alpha_{A,-,I}]{650}1b %I(ME)
 \putmorphism(460,50)(1,0)[\phantom{A\ot B}`A- `A\rho_{-}]{780}1b

\putmorphism(-180,400)(0,-1)[\phantom{Y_2}``=]{350}1r
\putmorphism(1280,400)(0,-1)[\phantom{Y_2}``=]{350}1l
\put(540,220){\fbox{$r_{A,-}$}}
\efig}
$$

% m from left action
$$%\begin{equation} \eqlabel{effects}
\scalebox{0.8}{
\bfig
 \putmorphism(-210,400)(1,0)[(AI)- ` A(I-) `\alpha_{A,I,-}]{670}1a
 \putmorphism(460,400)(1,0)[\phantom{A\ot B}`A- `A\lambda_{-}]{800}1a
\putmorphism(-150,50)(1,0)[(AI)-`A-`\rho_A\ltimes-]{1430}1b

\putmorphism(-180,400)(0,-1)[\phantom{Y_2}``=]{350}1r
\putmorphism(1280,400)(0,-1)[\phantom{Y_2}``=]{350}1l
\put(540,220){\fbox{$m_{A,-}$}}
\efig}
\qquad
% m from right action
\scalebox{0.8}{
\bfig
 \putmorphism(-210,400)(1,0)[(-I)B`-(IB) `\alpha_{-,I,M}]{670}1a
 \putmorphism(460,400)(1,0)[\phantom{A\ot B}`-B `-\rtimes\lambda_B]{800}1a
\putmorphism(-150,50)(1,0)[(-I)B`-B`\rho_{-}B]{1430}1b

\putmorphism(-180,400)(0,-1)[\phantom{Y_2}``=]{350}1r
\putmorphism(1280,400)(0,-1)[\phantom{Y_2}``=]{350}1l
\put(540,220){\fbox{$m_{-,B}$}}
\efig}
$$
so that   
 \begin{enumerate}[i)]
 \item the following four groups of 1h-cell components of modifications coincide: a) of the four $p$'s, b) of $m_{A,-}$ and $m_{-,B}$, 
c) of $l_{A,-}$ and $l_{-,B}$, and d) of $r_{A,-}$ and $r_{-,B}$, for all $A,B\in\Dd$; 
\item 2-cell components of the above ten horizontal modifications obey the axioms as in a monoidal bicategory ((TA1)-(TA3) of \cite{GPS}). 
\end{enumerate}
\end{enumerate}
\end{defn}

The definition of a horizontally premonoidal double category is such that its underlying bicategory together with the underlying bicategorical transformations of the horizontal equivalence in point 2. and the underlying modifications of the horizontal modifications in point 3. is a premonoidal bicategory from \cite{HF1}. 

Centrality of the 1h-cell components of the equivalences in point 2. are required in order for the following to be pseudonatural 
transformations: $-\rtimes\alpha_{B,C,D}$ from the first and $\alpha_{A,B,C}\ltimes-$ from the fourth pentagon, 
$\lambda_A\ltimes-$ from $l_{A,-}$, \,\, $-\rtimes\rho_B$ from $r_{-,B}$, \,\, $\rho_A\ltimes-$ from $m_{A,-}$ and $-\rtimes\lambda_B$ from 
$m_{-,B}$.

\smallskip

As explained in the proof of \cite[Proposition 3.15]{Fem3}, and by \prref{lifting 1v to equiv} and \prref{essence} one has: 

\begin{prop}
A vertically premonoidal double category whose associativity and unitality constraints are liftable vertical 
transformations induces a horizontally premonoidal double category. 
\end{prop}

\section{Locally cubical bicategories of (monoidal) double categories}

A locally cubical bicategory was introduced in \cite[Definition 11]{GG} as a bicategory enriched in the monoidal
2-category $Dbl_2$ of pseudodouble categories, pseudodouble functors and vertical strict transformations. Concerning the cells 
it consists of, it contains objects and the hom-parts $[A,B]$ are double categories. 
As any biclosed monoidal bicategory can be viewed as a bicategory enriched over itself, the biclosed monoidal 2-category $Dbl_2$ is a bicategory enriched over 
$Dbl_2$, {\em i.e.} locally cubical bicategory \textfrak{Dbl}, that we will describe further below. 
It is interesting to note that a one-object locally cubical bicategory yields a vertically monoidal double category. 

%\bigskip

\begin{rem}
We may think of a locally cubical bicategory schematically as below. \vspace{0,2cm} %double category 
%
%\hspace{3,5cm} $\bullet$ 0-cells \vspace{0,3cm} 
%
%\hspace{0,5cm} $\bullet$  vertical 1-cells \hspace{2cm} $\bullet$  horizontal 1-cells \vspace{0,3cm} 
%
%\hspace{7cm} $\bullet$ squares (2-cells)  \vspace{0,2cm}
%
%\noindent where the diagonal part represents its underlying bicategory. Now, if we think of a double category as of a diamond: \vspace{0,2cm}
%
%\hspace{3,5cm} $\bullet$ 0 \vspace{0,3cm} 
%
%\hspace{2cm}  1v $\bullet$  \hspace{1,5cm} $\bullet$  1h \vspace{0,3cm} 
%
%\hspace{3,5cm} $\bullet$ 2  \vspace{0,2cm}
%
%\noindent then we can think of a locally cubical bicategory as  of a {\em hanging diamond} in the sense \vspace{0,2cm}

\hspace{3,5cm} $\bullet$ 0 \vspace{0,3cm} 

\hspace{3,5cm} $\bullet$ 1 \vspace{0,3cm} 

\hspace{2cm}  2v $\bullet$  \hspace{1,5cm} $\bullet$  2h \vspace{0,3cm} 

\hspace{3,5cm} $\bullet$ 3  \vspace{0,2cm}

%\hspace{7cm} $\bullet$ cubes (3-cells) %\pause \vspace{0,4cm}
We will use this in \sssref{perspective lcb}. 
\end{rem}

\smallskip

We will consider three locally cubical bicategories. Let us first understand their hom-parts. % {\em i.e.} ``diamond parts''. 

%\subsection{Three locally cubical bicategories} 

\subsection{Double categories of (monoidal) double functors} \sslabel{mon d cats}

Let $[\Cc, \Dd]$ be the double category of pseudodouble functors, vertical strict transformations, horizontal pseudonatural 
transformations and modifications. We will consider its vertically and horizontally monoidal versions. 

\medskip

Shulman introduced a definition of a monoidal double category using vertical transformations. 
We refer to this definition as to {\em vertically monoidal} double category. In \cite[Section 10.2]{Fem3} we introduced a 
horizontal version of the notion. We first recall these definitions.

\begin{defn} \cite[\text{Definition 2.9}]{Shul} \delabel{Shul} \\ %\vspace{-1cm} \vspace{-0,4cm} 
A {\em monoidal double category} is a double category $\Dd$ equipped with pseudodouble functors $\ot:\Dd\times\Dd\to\Dd$ and $I:*\to\Dd$ and 
invertible vertical strict transformations 
$$\alpha: \ot\comp(\Id\times\ot) \stackrel{\iso}{\to} \ot\comp(\ot\times\Id) $$
$$\lambda: \ot\comp(I\times\Id) \stackrel{\iso}{\to} \Id $$
$$\rho: \ot\comp(\Id\times I) \stackrel{\iso}{\to} \Id $$ 
satisfying the pentagonal and three triangular axioms (via identity vertical modifications). 
\end{defn}

\begin{defn} \cite[Definition 10.9]{Fem3} \delabel{hor mon D}
For a double category $\Dd$ we say that it is {\em horizontally monoidal} if there are pseudodouble functors 
$\ot:\Dd\times\Dd\to\Dd$ and $I:*\to\Dd$, horizontal equivalence transformations 
$$\alpha: \ot\comp(\Id\times\ot) \stackrel{\iso}{\to} \ot\comp(\ot\times\Id) $$
$$\lambda: \ot\comp(I\times\Id) \stackrel{\iso}{\to} \Id $$
$$\rho: \ot\comp(\Id\times I) \stackrel{\iso}{\to} \Id ,$$ 
and horizontal modifications $p,m,l,r$ whose 2-cell components satisfy axioms (TA1)-(TA3) as in \cite{GPS} (write the latter as equations of pasted horizontally globular 2-cells). 
\end{defn} 

The above definition is such that for a horizontally monoidal double category $(\Dd,\ot,\alpha,\lambda,\rho)$ 
the underlying horizontal bicategory $\HH(\Dd)$ is a monoidal bicategory. 

\medskip

We have:

\begin{thm} \thlabel{two monoidalities}
A monoidal double category $(\Dd, \ot, I, \alpha,\lambda,\rho)$ in which $\alpha,\lambda,\rho$ are liftable yields a horizontally monoidal double category $(\Dd, \ot, I, \hat\alpha, \hat\lambda,\hat\rho)$. 
\end{thm}

\subsubsection{Double category of vertically monoidal double functors}

We start with the double category of vertically monoidal double functors.

\begin{prop} \cite[Proposition 5.6]{GGV} \prlabel{dbl cat of mon} \\
Let $\Cc, \Dd$ be vertically monoidal double categories.
There is a double category ${}^v\w\MMonDbl[\Cc, \Dd]$ of lax monoidal (pseudo) double functors, pseudomonoidal
vertical strict transformations, monoidal horizontal pseudonatural transformations, and monoidal modifications. 
\end{prop}

The cells in the above double category are meant in the sense of Definitions 5.1, 5.2, 5.3 and 5.4 from \cite{GGV}. In particular: \vspace{-0,2cm}
%, this means that: \vspace{-0,2cm} %(in the sense of \cite[Definition 5.1]{GGV})
\begin{itemize} 
\item for lax monoidal pseudodouble functors $F:\Cc\to\Dd$ there exist vertical strict transformations 
$F^2: \ot\comp(F\times F)\Rightarrow F\comp(-\ot-)$ and $F^0:I_\Dd\Rightarrow F\comp I_\Cc$ (here $I_\Cc$ and $I_\Dd$ are the monoidal unit double functors), 
whereby for $F^2$ there are %cocycle {\em vertical} 
1v-cells $F^2_{A,B}$ and 2-cell components $F^2_{f,g}$ for 1h-cells $f,g$ that are square-formed 
$$
\scalebox{0.86}{
\bfig
\putmorphism(-100,50)(1,0)[F(A)\ot F(B)`F(A')\ot F(B')`F(f)\ot F(g)]{1650}1a
\putmorphism(-100,-400)(1,0)[F(A\ot B)`F(A'\ot B'), `F(f\ot g)]{1650}1a
\putmorphism(-80,50)(0,-1)[\phantom{Y_2}``F^2_{A,B}]{420}1l
\putmorphism(1560,50)(0,-1)[\phantom{Y_2}``F^2_{A',B'}]{420}1r
\put(500,-140){\fbox{$F^2_{f,g}$}}
\efig} \vspace{-0,2cm}
$$ 
and there is a 1v-cell $F^0:I_\Dd\to F(I_\Cc)$ and a 2-cell $F^*$ inevitably given by $\frac{\Id_{F^0}}{s_0}$ where $s_0$ is the unit structure isomorphism of $F$, 
so that the 1v-cells $F^2_{A,B}$ and $F^0$, on one hand, and their corresponding 2-cells, on the other hand, satisfy the respective 
associativity and two unitality laws making the 1-dimensional functors $F_0, F_1$ lax monoidal; 
\item pseudomonoidal vertical strict transformations $\sigma:F\Rightarrow \tilde F$ have 2-cell components $\sigma^2_{A,B}$ and 
$\sigma^0$ that are vertically globular and invertible
$$
\scalebox{0.9}{
\bfig
 \putmorphism(-90,500)(1,0)[F(A)\ot F(B)`F(A)\ot F(B) `=]{740}1a
\putmorphism(-120,500)(0,-1)[\phantom{Y_2}`F(A\ot B) `F^2_{A,B}]{400}1l
\putmorphism(650,500)(0,-1)[\phantom{Y_2}`\tilde F(A)\ot \tilde F(B) `\sigma(A)\ot\sigma(B)]{400}1r
\putmorphism(-120,100)(0,-1)[\phantom{Y_2}``\sigma_{A\ot B}]{400}1l
\putmorphism(650,100)(0,-1)[\phantom{Y_2}``\tilde F^2_{A,B}]{400}1r
\putmorphism(-90,-300)(1,0)[\tilde F(A\ot B)`\tilde F(A\ot B) `=]{740}1a
\put(130,70){\fbox{$\sigma^2_{A,B}$}}
\efig} \qquad
\scalebox{0.9}{
\bfig
 \putmorphism(-90,500)(1,0)[I` I`=]{440}1a
\putmorphism(350,500)(0,-1)[` `\tilde F^0]{800}1r
\putmorphism(-90,-300)(1,0)[\tilde F(I)`\tilde F(I), `=]{440}1a
\putmorphism(-100,500)(0,-1)[` F(I)`F^0]{400}1l
\putmorphism(-100,110)(0,-1)[\phantom{Y_2}``\sigma(I)]{400}1l
\put(60,70){\fbox{$\sigma^0$}}
\efig} \vspace{-0,2cm}
$$ 
%\vspace{-0,2cm}
\item monoidal horizontal pseudonatural transformations $\beta:F\Rightarrow G$ have 2-cell components $\beta^2_{A,B}$ and $\beta^0$ 
that are square-formed
$$
\scalebox{0.86}{
\bfig
\putmorphism(-100,50)(1,0)[F(A)\ot F(B)`G(A)\ot G(B)`\beta(A)\ot\beta(B)]{1650}1a
\putmorphism(-100,-400)(1,0)[F(A\ot B)`G(A\ot B) `\beta(A\ot B)]{1650}1a
\putmorphism(-80,50)(0,-1)[\phantom{Y_2}``F^2_{A,B}]{420}1l
\putmorphism(1560,50)(0,-1)[\phantom{Y_2}``G^2_{A,B}]{420}1r
\put(560,-140){\fbox{$\beta^2_{A,B}$}}
\efig}\qquad
\scalebox{0.9}{
\bfig
\putmorphism(-150,170)(1,0)[I` I`=]{460}1a
\putmorphism(-150,-180)(1,0)[F(I)`G(I).`\beta(I)]{460}1b
\putmorphism(-150,170)(0,-1)[\phantom{Y_2}``F^0]{350}1l
\putmorphism(310,170)(0,-1)[\phantom{Y_2}``G^0]{350}1r
\put(0,0){\fbox{$\beta^0$}}
\efig}
$$ 
\end{itemize}

\medskip

\begin{rem} \rmlabel{modif+coc}
Associativity and two unitality laws making the functor $F_0$ lax monoidal express that $F^2$ determines a 
1-dimensional left and right normalized 2-cocycle via 1v-cells $F^2_{A,B}$, whereas 
associativity and two unitality laws making the functor $F_1$ lax monoidal express that $F^2$ determines a 
2-dimensional left and right normalized 2-cocycle via 2-cell components $F^2_{f,g}$. 

From the definition of a monoidal horizontal pseudonatural transformation $\beta$ from \cite[Definition 5.3]{GGV} it is clear that 
it is made of a modification $\beta^2:\beta(-)\ot\beta(-)\Rrightarrow\beta(-\ot-)$ (the first and the last axiom from the definition 
are precisely the modification axioms \axiomref{m.ho-vl.-1} and \axiomref{m.ho-vl.-2}) that is simultaneously a left and right 
normalized 2-dimensional 2-cocycle. Similarly, $\sigma^2$ of a pseudomonoidal vertical strict transformation $\sigma$ is a 
vertical modification that is simultaneously a left and right normalized 2-dimensional vertical 2-cocycle. Modifications that are 
 simultaneously left and right normalized 2-dimensional 2-cocycles we will call {\em cocycle modifications}. 
% \textcolor{rojo}{
We defer the explicit definition and study of double categorical cocycles to an 
eventual future work. The meaning of a 1- and 2-dimensional 2-cocycles should be clear from the context. The adjective 2-dimensional 
has the meaning that the constraint is expressed by an equation of 2-cells. 
The axioms that the horizontal modifications $p, m, l$ of a horizontally monoidal double category from \deref{hor mon D} obey mean that 
$p$ is a 2-dimensional left and right normalized 4-cocycle. The axioms that the ten horizontal modifications of a horizontally premonoidal double category 
from \deref{premon} obey mean that they all together make a 2-dimensional left and right normalized 4-cocycle $p$. %}
\end{rem}

\begin{defn}
A {\em strongly lax monoidal pseudodouble functor} is a lax monoidal pseudodouble functor $F:\Cc\to\Dd$ for which the 
vertical strict transformations $F^2: \ot\comp(F\times F)\Rightarrow F\comp(-\ot-)$ and $F^0:I_\Dd\Rightarrow F\comp I_\Cc$ 
are invertible. 

A {\em liftable} strongly lax monoidal pseudodouble functor is a strongly lax monoidal pseudodouble functor $F$ whose 
underlying vertical transformations $F^2, F^0$ are liftable. 
\end{defn}

\begin{defn}
A {\em strongly monoidal horizontal pseudonatural transformation} is a 
monoidal horizontal pseudonatural transformation $\beta$ between strongly lax monoidal pseudodouble functors 
whose modification $\beta^2$ is invertible (the 2-cell components $\beta^2_{A,B}$ are vertically invertible). 
\end{defn}

\medskip

Observe that in the  double category ${}^v\w\MMonDbl[\Cc, \Dd]$ we have two bicategories sitting inside. One is the 2-category 
of lax monoidal pseudodouble functors, pseudomonoidal vertical strict transformations and vertical monoidal modifications, which 
we will denote by ${}^v\w\MMonDbl_v[\Cc, \Dd]_2$. The other is the bicategory of lax monoidal pseudodouble functors, monoidal horizontal pseudonatural transformations, and horizontal monoidal modifications that we denote by ${}^v\w\MMonDbl_h[\Cc, \Dd]_2$. 

\smallskip

We also have:

\begin{prop} \cite[Proposition 5.6]{GGV} \prlabel{lift in dbl cat of mon} \\
A pseudomonoidal vertical strict transformation $\sigma: F \Rightarrow \tilde F$ has a companion as a vertical
1-cell in ${}^v\w\MMonDbl[\Cc, \Dd]$ if and only if the underlying vertical transformation of $\sigma$ has a companion as
a vertical 1-cell of $[\Cc, \Dd]$, {\em i.e.} it is liftable.
\end{prop}

By \prref{lifting 1v to equiv} and \prref{essence} we have:

\begin{cor} \colabel{2-cat lift}
Let ${}^v\w\MMonDbl_{v^0}[\Cc, \Dd]_2$ denote the sub-2-category of ${}^v\w\MMonDbl_v[\Cc, \Dd]_2$ whose 1-cells are liftable. There is a faithful 
2-functor 
$${}^v\w\MMonDbl_{v^0}[\Cc, \Dd]_2\hookrightarrow {}^v\w\MMonDbl_h[\Cc, \Dd]_2.$$
\end{cor}

\subsubsection{Double category of horizontally monoidal double functors}

The above Corollary is about lifting vertical transformations and modifications inside of ${}^v\w\MMonDbl[\Cc, \Dd]$ into their horizontal counterparts. 
We next introduce a horizontal version of ${}^v\w\MMonDbl[\Cc, \Dd]$ that we will denote by ${}^h\w\MMonDbl[\Cc, \Dd]$. 
Let ${}^{v_0}\w\MMonDbl[\Cc, \Dd]$ denote the version of the double category ${}^v\w\MMonDbl[\Cc, \Dd]$ in which $\Cc$ and $\Dd$ have liftable 
vertical monoidal structures, {\em i.e.} their associativity and unity vertical transformations are liftable.
Then we will match 
a double subcategory of ${}^{v_0}\w\MMonDbl[\Cc, \Dd]$ and the double category ${}^h\w\MMonDbl[\Cc, \Dd]$ via a lifting kind of an assignment. 

\bigskip

Instead of giving a rigorous definition of ${}^h\w\MMonDbl[\Cc, \Dd]$, we will stress the differences with respect to 
${}^v\w\MMonDbl[\Cc, \Dd]$. To begin with, let $\Cc,\Dd$ be {\em horizontally} monoidal double categories, then ${}^h\w\MMonDbl[\Cc, \Dd]$ is the double category consisting of: \vspace{-0,2cm}
\begin{itemize} 
\item {\em horizontally} lax monoidal pseudodouble functors (this means that the cocycle 1-cells $F^2_{A,B}$ are {\em horizontal}, the %correspondingly the 
2-cell components $F^2_{f,g}$ are horizontally globular and invertible
$$
\scalebox{0.86}{
\bfig
 \putmorphism(-470,500)(1,0)[F(A)\ot F(B)`F(A')\ot F(B')`F(f)\ot F(g)]{1040}1a
 \putmorphism(760,500)(1,0)[\phantom{F(f)}`F(A'\ot B') `F^2_{A',B'}]{700}1a
 \putmorphism(-470,120)(1,0)[F(A)\ot F(B)`F(A\ot B)`F^2_{A,B}]{840}1a
 \putmorphism(460,120)(1,0)[\phantom{G(B)}`F(A'\ot B') `F(f\ot g)]{960}1a
\putmorphism(-480,500)(0,-1)[\phantom{Y_2}``=]{380}1r
\putmorphism(1440,500)(0,-1)[\phantom{Y_2}``=]{380}1r
\put(280,310){\fbox{$F^2_{f,g}$}}
\efig} \vspace{-0,2cm}
$$ 
and the associativity and two unitality laws for the 1-dimensional functors $F_0, F_1$ to be lax monoidal are governed by invertible horizontal modifications: % $\omega, \zeta, \xi$; 
$$ 
\scalebox{0.82}{
\bfig
\putmorphism(-390,450)(0,-1)[\phantom{Y_2}` `=]{450}1r 

\putmorphism(-390,450)(1,0)[(F(A)F(B))F(C)` F(AB)F(C)`F^2_{A,B}F(C)]{1170}1a 
\putmorphism(840,450)(1,0)[\phantom{A''\ot B'}`F((AB)C) ` F^2_{AB,C}]{1170}1a%%%%%
 \putmorphism(2040,450)(1,0)[\phantom{A''\ot B'}` F(A(BC)) `F(\alpha_{A,B,C})]{1160}1a %T(T_A(T_BT_C))

 \putmorphism(-390,0)(1,0)[(F(A)F(B))F(C)` F(A)(F(B)F(C)) `\alpha_{F(A),F(B),F(C)}]{1350}1b 
 \putmorphism(1120,0)(1,0)[\phantom{A''\ot B'}` F(A)F(BC) `F(A)F^2_{B,C}]{1040}1b 
 \putmorphism(2230,0)(1,0)[\phantom{A''\ot B'}` F(A(BC)) `F^2_{A,F(BC)}]{1020}1b 

\putmorphism(3200,450)(0,-1)[\phantom{Y_2}``=]{450}1r 

\put(1100,230){\fbox{$\omega_{A,B,C}$}}
\efig}
$$

$$ 
\scalebox{0.82}{
\bfig
\putmorphism(-390,450)(0,-1)[\phantom{Y_2}` `=]{450}1r 

\putmorphism(-390,450)(1,0)[F(A)I` F(A)F(I) `F(A)F^0]{740}1a 
\putmorphism(370,450)(1,0)[\phantom{A''\ot B'}` F(AI) ` F^2_{A,I} ]{700}1a %%%%%
 \putmorphism(1020,450)(1,0)[\phantom{A''\ot B'}` F(A) `F(\rho_A) ]{760}1a %T(T_A(T_BT_C))

 \putmorphism(-390,0)(1,0)[F(A)I` F(A) `\rho_{F(A)}]{2230}1b 

\putmorphism(1820,450)(0,-1)[\phantom{Y_2}``=]{450}1r 

\put(700,210){\fbox{$\delta_A$}}
\efig}
$$

$$ 
\scalebox{0.82}{
\bfig
\putmorphism(-390,450)(0,-1)[\phantom{Y_2}` `=]{450}1r 

\putmorphism(-390,450)(1,0)[IF(A)` F(I)F(A) `F^0F(A)]{740}1a 
\putmorphism(370,450)(1,0)[\phantom{A''\ot B'}` F(IA) ` F^2_{I,A} ]{700}1a %%%%%
 \putmorphism(1020,450)(1,0)[\phantom{A''\ot B'}` F(A) `F(\lambda_A) ]{760}1a %T(T_A(T_BT_C))

 \putmorphism(-390,0)(1,0)[IF(A)` F(A); `\lambda_{F(A)}]{2230}1b 

\putmorphism(1820,450)(0,-1)[\phantom{Y_2}``=]{450}1r 

\put(700,210){\fbox{$\gamma_A$}}
\efig}
$$

\item pseudomonoidal vertical strict transformations (the 2-cell components $\sigma^2_{A,B}$ and $\sigma^0$ are correspondingly square-formed)
$$
\scalebox{0.86}{
\bfig
\putmorphism(-100,50)(1,0)[F(A)\ot F(B)`F(A\ot B)`F^2_{A,B}]{1650}1a
\putmorphism(-100,-400)(1,0)[\tilde F(A)\ot\tilde F(B)`\tilde F(A\ot B) `\tilde F^2_{A,B}]{1650}1a
\putmorphism(-80,50)(0,-1)[\phantom{Y_2}``\sigma(A)\ot\sigma(B)]{420}1l
\putmorphism(1560,50)(0,-1)[\phantom{Y_2}``\sigma(A\ot B)]{420}1r
\put(550,-140){\fbox{$\sigma^2_{A,B}$}}
\efig}
\qquad\quad
\scalebox{0.86}{
\bfig
\putmorphism(-100,50)(1,0)[I` F(I) `F^0]{650}1a
\putmorphism(-100,-400)(1,0)[I` \tilde F(I) `\tilde F^0]{650}1a
\putmorphism(-80,50)(0,-1)[\phantom{Y_2}``=]{420}1l
\putmorphism(560,50)(0,-1)[\phantom{Y_2}``\sigma_I]{420}1r
\put(150,-160){\fbox{$\sigma^0$}}
\efig}
$$ 
\vspace{-0,2cm}
\item monoidal horizontal pseudonatural transformations (the 2-cell components $\beta^2_{A,B}$ and $\beta^0$ are correspondingly horizontally globular and invertible) 
$$
\scalebox{0.86}{
\bfig
 \putmorphism(-570,500)(1,0)[F(A)\ot F(B)`G(A)\ot G(B)`\beta(A)\ot\beta(B)]{1140}1a
 \putmorphism(740,500)(1,0)[\phantom{F(f)}`G(A\ot B) `G^2_{A,B}]{700}1a
 \putmorphism(-570,120)(1,0)[F(A)\ot F(B)`F(A\ot B)`F^2_{A,B}]{940}1a
 \putmorphism(460,120)(1,0)[\phantom{G(B)}`G(A\ot B) `\beta(A\ot B)]{960}1a
\putmorphism(-580,500)(0,-1)[\phantom{Y_2}``=]{380}1r
\putmorphism(1440,500)(0,-1)[\phantom{Y_2}``=]{380}1r
\put(280,310){\fbox{$\beta^2_{A,B}$}}
\efig}
\qquad\quad
\scalebox{0.86}{
\bfig
 \putmorphism(0,500)(1,0)[I`I`=]{540}1a
 \putmorphism(480,500)(1,0)[\phantom{F(f)}`I `G^0]{600}1a
 \putmorphism(0,120)(1,0)[I `F(I)`F^0]{500}1a
 \putmorphism(460,120)(1,0)[\phantom{G(B)}`G(I) `\beta_I]{620}1a
\putmorphism(0,500)(0,-1)[\phantom{Y_2}``=]{380}1r
\putmorphism(1080,500)(0,-1)[\phantom{Y_2}``=]{380}1r
\put(420,310){\fbox{$\beta^0$}}
\efig}
$$ 
and \vspace{-0,2cm}
\item monoidal modifications. \vspace{-0,2cm}
\end{itemize} 
The axioms that $\sigma^2_{A,B}$ and $\beta^2_{A,B}$ obey are correspondingly accommodated with respect to those from 
${}^v\w\MMonDbl[\Cc, \Dd]$. 

\begin{rem} \rmlabel{hor lax mon f}
Observe that {\em horizontally} lax monoidal pseudodouble functors $\Cc\to\Dd$ are precisely the lax monoidal pseudofunctors 
$\HH(\Cc)\to\HH(\Dd)$ between the underlying bicategories (in the sense of \cite[Section 2.2]{MV}, or {\em weak monoidal} homomorphisms 
of \cite[Definition 2]{DS}; $\omega, \delta, \gamma$ above are shuffled (double categorical) versions of $\omega, \zeta, \xi$ from 
\cite{MV} adjusted for the proof of \prref{lax monoidal-comm h}). 

Moreover, monoidal horizontal pseudonatural transformations $F\Rightarrow G:\Cc\to\Dd$ in ${}^h\w\MMonDbl[\Cc, \Dd]$ are monoidal pseudonatural transformations $\HH(F)\Rightarrow \HH(G):\HH(\C)\to\HH(\Dd)$ in the underlying bicategories satisfying additionally the following law
\begin{equation}\eqlabel{h mon tr law}
\scalebox{0.86}{
\bfig
\putmorphism(-550,400)(1,0)[F(A)F(B)`G(A)G(B)`\beta(A)\beta(B)]{900}1a
\putmorphism(430,400)(1,0)[\phantom{F(A)}`G(AB) `G^2_{A,B}]{560}1a
 \putmorphism(-570,0)(1,0)[F(A)F(B)`F(AB)`F^2_{A,B}]{920}1a
 \putmorphism(350,0)(1,0)[\phantom{F(A)}`G(AB) `\beta_{AB}]{660}1a

\putmorphism(-580,400)(0,-1)[\phantom{Y_2}``=]{400}1l
\putmorphism(1010,400)(0,-1)[\phantom{Y_2}``=]{400}1l
\put(240,210){\fbox{$\beta^2_{A,B}$}}

\putmorphism(380,0)(0,-1)[\phantom{Y_3}``F(uv)]{400}1l
\putmorphism(1010,0)(0,-1)[\phantom{Y_2}``G(uv)]{400}1r
 \putmorphism(350,-400)(1,0)[F(\tilde A\tilde B)` G(\tilde A\tilde B) `\beta(\tilde A\tilde B)]{600}1a
\put(600,-160){\fbox{$\beta^{u,v}$}}
\efig}
=
\scalebox{0.86}{
\bfig
\putmorphism(-550,400)(0,-1)[\phantom{Y_2}``F(uv)]{400}1l
\putmorphism(320,400)(0,-1)[\phantom{Y_3}``G(uv)]{400}1r
 \putmorphism(-550,400)(1,0)[F(A)F(B)`G(A)G(B)`\beta(A)\beta(B)]{900}1a
\put(-210,200){\fbox{$\beta^u\beta^v$}}

\putmorphism(-550,0)(1,0)[F(\tilde A)F(\tilde B)`G(\tilde A)G(\tilde B)`\beta(\tilde A)\beta(\tilde B)]{900}1b
\putmorphism(450,0)(1,0)[\phantom{F(A)}`G(\tilde A\tilde B) `G^2_{A,B}]{560}1a
 \putmorphism(-570,-400)(1,0)[F(\tilde A)F(\tilde B)`F(\tilde A\tilde B)`F^2_{A,B}]{820}1a
 \putmorphism(150,-400)(1,0)[\phantom{F(A)F(\tilde B)}`G(\tilde A\tilde B). `\beta(\tilde A\tilde B)]{820}1a

\putmorphism(-550,0)(0,-1)[\phantom{Y_2}``=]{400}1r
\putmorphism(1000,0)(0,-1)[\phantom{Y_2}``=]{400}1r
\put(160,-210){\fbox{$\beta^2_{\tilde A,\tilde B}$}}
\efig}
\end{equation}
At last, the underlying double functors for the monoidal horizontal pseudonatural transformations are 
horizontally lax monoidal pseudodouble functors. Then inspecting \cite[Definition 5.5]{GGV} it 
becomes clear that the modifications are precisely the bicategorical monoidal modifications with 2-cell components 
in the underlying bicategory $\HH(\Dd)$.
\end{rem}

Let $\K,\Ll$ be monoidal 2-categories and let $\MonBicat(\K,\Ll)$ denote the 2-category of lax monoidal pseudofunctors, 
monoidal pseudonatural transformations and monoidal modifications (as in \cite{ChG}). From the above is clear that there is a 
faithful 2-functor 
$$\HH({}^h\w\MMonDbl[\Cc, \Dd]) \stackrel{\F_{\equref{h mon tr law}}}{\to} \MonBicat(\HH(\Cc), \HH(\Dd))$$
from the underlying 2-category of the double category ${}^h\w\MMonDbl[\Cc, \Dd]$ forgetting the property \equref{h mon tr law} on 1-cells.

\subsubsection{Lifting of vertical to horizontal monoidal double structures}

Finally, we relate the vertically monoidal structures in ${}^v\w\MMonDbl[\Cc, \Dd]$ with the horizontally monoidal structures in ${}^h\w\MMonDbl[\Cc, \Dd]$ 
for double categories $\Cc$ and $\Dd$ which allow that.

\begin{prop} \prlabel{ver-hor-dbl}
Let $\Cc, \Dd$ be vertically monoidal double categories with liftable monoidal structures, {\em i.e.} their associativity and unity 
vertical transformations are liftable. Let ${}^{v_0}\w\MMonDbl[\Cc, \Dd]_0$ denote the double subcategory of ${}^v\w\MMonDbl[\Cc, \Dd]$ in which the objects are the liftable strongly lax monoidal pseudodouble functors, and the 1h-cells are strongly monoidal horizontal 
pseudonatural transformations. Then there is a double functor 
$$L_{\Cc,\Dd}:{}^{v_0}\w\MMonDbl[\Cc, \Dd]_0\hookrightarrow{}^h\w\MMonDbl[\Cc, \Dd].$$
\end{prop}

\begin{proof}
First of all, by the first assumption and \thref{two monoidalities} we will consider $\Cc, \Dd$ as horizontally monoidal double categories. %Supposing that we discussed $L_{\Cc,\Dd}$ on objects, we know by \prref{lifting 1v to equiv} that the underlying vertical strict transformations $\sigma$ from ${}^v\w\MMonDbl[\Cc, \Dd]_0$ lift to horizontal pseudonatural transformations $\hat\sigma$ in 
%${}^h\w\MMonDbl[\Cc, \Dd]$. 
The definition of $L_{\Cc,\Dd}$ on the cells is as follows: the underlying pseudodouble functors $F$, the underlying vertical transformations 
$\sigma$ and horizontal transformations $\beta$, as well as the modifications $a$ remain unchanged, we only change the monoidality structures of the cells in the following way: \vspace{-0,2cm} 
%1-cell components $\sigma(A) and $\beta(A)$ 
%, we should prove that \vspace{-0,2cm} 
%vertical cocycle 1-cells $F^2_{A,B}$ lift to horizontal cocycle 1-cells $\hat{F^2_{A,B}}$, and that 
\begin{itemize}
\item the square-formed lax monoidal structure 2-cell components $F^2_{f,g}$ of the objects lift to their horizontally globular associates 
$\widehat{F^2_{f,g}}$
$$
\scalebox{0.86}{
\bfig
 \putmorphism(-470,500)(1,0)[F(A)\ot F(B)`F(A')\ot F(B')`F(f)\ot F(g)]{1040}1a
 \putmorphism(760,500)(1,0)[\phantom{F(f)}`F(A'\ot B') `\widehat{F^2_{A',B'}}]{700}1a
 \putmorphism(-470,120)(1,0)[F(A)\ot F(B)`F(A\ot B)`\widehat{F^2_{A,B}}]{840}1a
 \putmorphism(460,120)(1,0)[\phantom{G(B)}`F(A'\ot B') `F(f\ot g)]{960}1a
\putmorphism(-480,500)(0,-1)[\phantom{Y_2}``=]{380}1r
\putmorphism(1440,500)(0,-1)[\phantom{Y_2}``=]{380}1r
\put(280,310){\fbox{$\widehat{F^2_{f,g}}$}}
\efig} \vspace{-0,2cm}
$$  \vspace{-0,2cm} 
\item the vertically globular 2-cell components $(\sigma^2_{A,B})^{-1}$ of 1v-cells lift to their square-formed associates 
$\widehat{\sigma^2_{A,B}}$
$$
\scalebox{0.86}{
\bfig
\putmorphism(-100,50)(1,0)[F(A)\ot F(B)`F(A\ot B)`\widehat{F^2_{A,B}}]{1650}1a
\putmorphism(-100,-400)(1,0)[F'(A)\ot F'(B)`F'(A\ot B) `\widehat{{F'}^2_{A,B}}]{1650}1a %\widehat{{\tilde F}^2_{A,B}}
\putmorphism(-80,50)(0,-1)[\phantom{Y_2}``\sigma(A)\ot\sigma(B)]{420}1l
\putmorphism(1560,50)(0,-1)[\phantom{Y_2}``\sigma(A\ot B)]{420}1r
\put(150,-160){\fbox{$\widehat{(\sigma^2_{A,B})^{-1}}$}}
\efig}
$$ 
 and \vspace{-0,2cm} 
\item the square-formed 2-cell components $\beta^2_{A,B}$ of 1h-cells lift to their horizontally globular associates 
$\widehat{\beta^2_{A,B}}$ 
$$
\scalebox{0.86}{
\bfig
 \putmorphism(-470,500)(1,0)[F(A)\ot F(B)`G(A)\ot G(B)`\beta(A)\ot\beta(B)]{1040}1a
 \putmorphism(760,500)(1,0)[\phantom{F(f)}`G(A\ot B) `\widehat{G^2_{A,B}}]{700}1a
 \putmorphism(-470,120)(1,0)[F(A)\ot F(B)`F(A\ot B)`\widehat{F^2_{A,B}}]{840}1a
 \putmorphism(460,120)(1,0)[\phantom{G(B)}`G(A\ot B) `\beta(A\ot B)]{960}1a
\putmorphism(-480,500)(0,-1)[\phantom{Y_2}``=]{380}1r
\putmorphism(1440,500)(0,-1)[\phantom{Y_2}``=]{380}1r
\put(280,300){\fbox{$\widehat{\beta^2_{A,B}}$}}
\efig}
$$ \vspace{-0,2cm} 
\end{itemize}
%and that the corresponding axioms are fulfilled. This all holds true by the argument below \equref{omega-gen}. 
and similarly the 2-cell components $F^*, (\sigma^0)^{-1}$ and $\beta^0$ at $I$ lift to $\hat F^*, 
\widehat{(\sigma^0)^{-1}}$ and $\hat{\beta^0}$, respectively. 
Thus defined images by $L_{\Cc,\Dd}$ of $F, \sigma, \beta$ and $a$ we will denote by $\check{F}, \check{\sigma}, \check{\beta}$ and 
$\check{a}$, respectively. 
By the argument below \equref{omega-gen} and \prref{lifting 1v to equiv} we have that the above liftings are well-defined and also that the corresponding axioms in ${}^h\w\MMonDbl[\Cc, \Dd]$ are fulfilled.

To verify that $L_{\Cc,\Dd}$ is a double functor we take two horizontal transformations 
\[
\xymatrix{ 
\Cc \ar@{:=}@/^1pc/[rr] |{\Downarrow \beta} \ar@/^{2pc}/[rr]^{F}  \ar[rr] 
  \ar@{:=}@/_1pc/[rr] |{\Downarrow \beta\s'}  \ar@/_{2pc}/[rr]_{H}  &&
	\Dd %\ar@{:=}@/^1pc/[rr] |{\Downarrow G} \ar@/^{2pc}/[rr]^{B}  \ar[rr] 
  %\ar@{:=}@/_1pc/[rr] |{\Downarrow G}  \ar@/_{2pc}/[rr]_{B''} && r
}\] 
By comparing the 2-cell components, it is directly seen that $\frac{L(\beta)}{L(\beta')}=L(\frac{\beta}{\beta'})$, as one inserts the identity 2-cell $\frac{\eta}{\Epsilon}=\Id$ holding among the companion-structure 2-cells. Exactly the same holds for two vertical transformations, and we have that $L_{\Cc,\Dd}$ is a (strict) double functor. 
\qed\end{proof}

\subsection{The announced locally cubical bicategories}

Now we turn to the announced locally cubical bicategories. 
At first, let \textfrak{Dbl} denote the locally cubical bicategory whose objects 
are double categories and the hom-part is $[\Cc, \Dd]$. 
%, pseudodouble functors, vertical strict transformations, horizontal pseudotransformations and modifications. 

\medskip

Next, let $\textfrak{VMonDbl}$ denote the locally cubical bicategory whose objects are %of 
vertically monoidal double categories and the hom-part is ${}^v\w\MMonDbl[\Cc, \Dd]$ from \prref{dbl cat of mon}. 
%lax monoidal pseudodouble functors, pseudomonoidal vertical strict transformations, monoidal horizontal pseudotransformations and monoidal modifications, where the ``diamond part'' is meant 

\medskip

Finally, let us denote by $\textfrak{HMonDbl}$ the locally cubical bicategory %of: \vspace{-0,2cm}
%\begin{itemize} 
%\item horizontally monoidal double categories, 
%\item horizontally lax monoidal pseudodouble functors (the cocycle 1-cells are horizontal and correspondingly the 2-cell components $F^2_{f,g}$ are horizontally globular), 
%\item pseudomonoidal vertical strict transformations (the 2-cell components $\sigma^2_{A,B}$ are correspondingly square-formed), 
%\item monoidal horizontal pseudotransformations (the 2-cell components $\beta^2_{A,B}$ are correspondingly horizontally globular), and 
%\item monoidal modifications, 
%\end{itemize}
%whereby the axioms that the cells of the ``diamond part'' obey are correspondingly accommodated with respect to those from $\textfrak{VMonDbl}$. 
whose objects are horizontally monoidal double categories and the hom-part is ${}^h\w\MMonDbl[\Cc, \Dd]$. 

\medskip

The 2-functor from \coref{2-cat lift} extends accordingly to an auto-homomorphism 
\begin{equation} \eqlabel{auto}
L_{tr}:\textfrak{VMonDbl}\to\textfrak{VMonDbl}
\end{equation}
on the locally cubical bicategory $\textfrak{VMonDbl}$ (see \cite[Section 3.5]{GSh} for the definition of a homomorphism of locally cubical bicategories). We are now interested in the following.

\begin{prop} \prlabel{ver-hor-cbl}
Let $\textfrak{V}^0\textfrak{MonDbl}_0$ be a locally cubical sub-bicategory of $\textfrak{VMonDbl}$ in which the vertically monoidal structures of objects
are liftable, and whose hom-objects are given by ${}^{v_0}\w\MMonDbl[\Cc, \Dd]_0$ (from \prref{ver-hor-dbl}). 
There is a homomorphism of locally cubical bicategories $L:\textfrak{V}^0\textfrak{MonDbl}_0\hookrightarrow\textfrak{HMonDbl}$ 
defined on objects by 
$$(\Dd, \ot, I, \alpha,\lambda,\rho)\mapsto(\Dd, \ot, I, \hat\alpha, \hat\lambda,\hat\rho)$$ 
(as in \thref{two monoidalities}) and on hom-double categories it is given by 
$$L_{\Cc,\Dd}:{}^{v_0}\w\MMonDbl[\Cc, \Dd]_0\hookrightarrow{}^h\w\MMonDbl[\Cc, \Dd].$$ 
\end{prop}

\begin{proof}
In order to define a homomorphism of locally cubical bicategories it remains to give 2-cells $\mu_{\Cc,\Dd,\Ee}$ and $\iota_\Cc$ in $Dbl_2$ for the compatibility of $L$ with the horizontal composition and unitality on the hom-double categories. These should be vertical transformations, which are given on objects, 1v- and 1h-cells of ${}^{v_0}\w\MMonDbl[\Dd, \Ee]_0\ot {}^{v_0}\w\MMonDbl[\Cc, \Dd]_0$. Respectively, we should 
give $\mu_{\Cc,\Dd,\Ee}$ and $\iota_\Cc$ when evaluated at pseudodouble functors $G\ot F$ (this should give a 1v-cell in 
${}^h\w\MMonDbl[\Cc, \Ee]$, {\em i.e.} a vertical transformation $\omega$), on vertical transformations $\sigma'\ot\sigma$ and on 
horizontal transformations $\beta'\ot\beta$ - the latter two should give 2-cells in ${}^h\w\MMonDbl[\Cc, \Ee]$, {\em i.e.} monoidal modifications. 
Since by giving these structure cells 
we are dealing with $\check G\ot\check F, \widecheck{G\ot F}, \,\, \check\sigma'\ot\check\sigma, \widecheck{\sigma'\ot\sigma}, \,\, 
\check\beta'\ot\check\beta, \widecheck{\beta'\ot\beta}$ evaluated at an object $C\in\Cc$, by the definition of $L$ on those cells the latter are 
simply equal to $(G\ot F)(C), (\sigma'\ot\sigma)(C)$ and $(\beta'\ot\beta)(C)$, respectively, so we may take for $\omega$ the identity vertical transformation and for the two desired modifications the identity modifications to complete the definition of $L$. Then all the desired axioms for $L$ to be a homomorphism of locally cubical bicategories are fulfilled and we have the claim. 
\qed\end{proof}

\section{Actions on double categories} 

In \ssref{mon d cats} we recalled the notions of vertically and horizontally monoidal double categories. 
In \cite[Section 10.3]{Fem3} we introduced vertical and horizontal actions of correspondingly monoidal double categories.  
We recall these definitions here.

\begin{defn} \delabel{ver-act} 
We say that a monoidal double category $\Mm$ acts (from the left) on a double category $\Ee$ if there is a pseudodouble functor 
$F:\Mm\times\Ee\to\Ee$, invertible vertical strict transformations with components 
$$\tilde\lambda_E: I\times E\to E\quad\text{and}\quad\tilde\alpha_{M,N,E}:(M\ot N)\times E\to M\times(N\times E)$$
with $M,N\in\Mm$ and $E\in\Ee$, and identity vertical modifications $\tilde p, \tilde l, \tilde m$, analogous to $p,l,m$ from \deref{hor mon D}.  \\
\end{defn}

This kind of action we will call a {\em vertical action}. To the contrast to it, we will differ what we will call a {\em horizontal action}. The difference is actually already anticipated by the fact that in the vertical action Shulman's kind of monoidality of the acting double category is meant, whereas in the horizontal action {\em horizontal monoidality} of the acting double category is assumed.
This ``horizontal'' action we define as follows. 

\begin{defn} \delabel{hor-act} 
By an action of a {\em horizontally monoidal} double category $\Mm$ on a double category $\Ee$ we mean the data comprised of: a 
pseudodouble functor $F:\Mm\times\Ee\to\Ee$, horizontal equivalences with components 
$$\tilde\lambda_E: I\times E\to E\quad\text{and}\quad\tilde\alpha_{M,N,E}:(M\ot N)\times E\to M\times(N\times E)$$
with $M,N\in\Mm$ and $E\in\Ee$, and horizontal modifications $\tilde p, \tilde m, \tilde l$ %, analogous to $p,m,l$ from \deref{hor mon D}.
with 2-cell components 
$$%\begin{equation} \eqlabel{effects}
\scalebox{0.8}{
\bfig
 \putmorphism(-150,400)(1,0)[((MN)P)E`(M(NP))E `\alpha^\Mm_{M,N,P}E]{1100}1a
 \putmorphism(1050,400)(1,0)[\phantom{A\ot B}`M((NP)E)`\tilde\alpha_{M,NP,E}]{1080}1a
 \putmorphism(2220,400)(1,0)[\phantom{A\ot B}` M(N(PE)) `M\tilde\alpha_{N,P,E}]{1120}1a

 \putmorphism(-170,50)(1,0)[((MN)P)E`(MN)(PE) `\tilde\alpha_{MN,P,E}]{1620}1b
 \putmorphism(1560,50)(1,0)[\phantom{A\ot B}`M(N(PE)) `\tilde\alpha_{M,N,PE}]{1800}1b %(M(NP))E

\putmorphism(-180,400)(0,-1)[\phantom{Y_2}``=]{350}1r
\putmorphism(3290,400)(0,-1)[\phantom{Y_2}``=]{350}1l
\put(1380,230){\fbox{$\tilde p^{M,N,P}_E$}}
\efig}
$$%\end{equation}

$$%\begin{equation} \eqlabel{effects}
\scalebox{0.8}{
\bfig
 \putmorphism(-150,400)(1,0)[(IM)E`ME `\lambda^\Mm_ME]{1430}1a

 \putmorphism(-210,50)(1,0)[(IM)E`I(ME) `\tilde\alpha_{I,M,E}]{670}1b %I(ME)
 \putmorphism(460,50)(1,0)[\phantom{A\ot B}`ME `\tilde\lambda_{ME}]{780}1b
 %\putmorphism(1260,50)(1,0)[\phantom{A\ot B}`ME`]{800}1a

\putmorphism(-180,400)(0,-1)[\phantom{Y_2}``=]{350}1r
\putmorphism(1280,400)(0,-1)[\phantom{Y_2}``=]{350}1l
\put(540,220){\fbox{$\tilde l^M_E$}}
\efig}
$$

$$%\begin{equation} \eqlabel{effects}
\scalebox{0.8}{
\bfig
 \putmorphism(-210,400)(1,0)[(MI)E`M(IE) `\tilde\alpha_{M,I,E}]{670}1a
 \putmorphism(460,400)(1,0)[\phantom{A\ot B}`ME `M\tilde\lambda_{E}]{800}1a
% \putmorphism(1260,400)(1,0)[\phantom{A\ot B}`ME`\rho^\Mm_ME]{800}1a
\putmorphism(-150,50)(1,0)[(MI)E`ME`\rho^\Mm_ME]{1430}1b

\putmorphism(-180,400)(0,-1)[\phantom{Y_2}``=]{350}1r
\putmorphism(1280,400)(0,-1)[\phantom{Y_2}``=]{350}1l
\put(540,220){\fbox{$\tilde m^M_E$}}
\efig}
$$
that satisfy axioms analogous to (TA1)-(TA2) of \cite{GPS} (see also \cite[Definition 4.4]{St}). 
\end{defn}

The axioms that the horizontal modifications $\tilde p, \tilde m, \tilde l$ in \deref{hor-act} obey mean that $\tilde p$ is a 2-dimensional left normalized 4-cocycle, according to \rmref{modif+coc}.

\begin{rem}\rmlabel{monad action}
The horizontal modifications $\tilde p, \tilde m, \tilde l$ in the above definition can be equivalently expressed as 
horizontal modifications
$$%\begin{equation} \eqlabel{effects}
\scalebox{0.8}{
\bfig
 \putmorphism(-150,400)(1,0)[M(N(PE))`M((NP)E) `M\tilde\alpha^{-1}_{N,P,E}]{1620}1a
 \putmorphism(1560,400)(1,0)[\phantom{A\ot B}`(M(NP))E `\tilde\alpha^{-1}_{M,NP,E}]{1800}1a

% \putmorphism(-170,400)(1,0)[M(N(PE))`M((NP)E) `M\tilde\alpha^{-1}_{N,P,E}]{1110}1a
% \putmorphism(1030,400)(1,0)[\phantom{A\ot B}`(M(NP))E `\tilde\alpha^{-1}_{M,NP,E}]{1080}1a
%  \putmorphism(2220,400)(1,0)[\phantom{A\ot B}`((MN)P)E `(\alpha^\Mm)^{-1}_{M,N,P}E]{1120}1a

%\putmorphism(-150,50)(1,0)[M(N(PE))`(MN)(PE) `\tilde\alpha^{-1}_{M,N,PE}]{1620}1b
% \putmorphism(1560,50)(1,0)[\phantom{A\ot B}`((MN)P)E`\tilde\alpha^{-1}_{MN,P,E}]{1800}1b

\putmorphism(-170,50)(1,0)[M(N(PE))`(MN)(PE) `\tilde\alpha^{-1}_{M,N,PE}]{1110}1b
 \putmorphism(1030,50)(1,0)[\phantom{A\ot B}`((MN)P)E`\tilde\alpha^{-1}_{MN,P,E}]{1080}1b
  \putmorphism(2220,50)(1,0)[\phantom{A\ot B}`(M(NP))E`\alpha^\Mm_{M,N,P}E]{1120}1b

\putmorphism(-180,400)(0,-1)[\phantom{Y_2}``=]{350}1r
\putmorphism(3290,400)(0,-1)[\phantom{Y_2}``=]{350}1l
\put(1380,180){\fbox{$\tilde p^{M,N,P}_E$}} % visina 220
\efig}
$$%\end{equation}

$$%\begin{equation} \eqlabel{effects}
\scalebox{0.8}{
\bfig
 \putmorphism(-210,400)(1,0)[ME`I(ME) `\tilde\lambda^{-1}_{ME}]{670}1a
 \putmorphism(460,400)(1,0)[\phantom{A\ot B}`(IM)E `\tilde\alpha^{-1}_{I,M,E}]{780}1a
 \putmorphism(1260,400)(1,0)[\phantom{A\ot B}`ME`\lambda^\Mm_ME]{800}1a
 \putmorphism(-150,50)(1,0)[ME`ME`=]{2230}1b

\putmorphism(-180,400)(0,-1)[\phantom{Y_2}``=]{350}1r
\putmorphism(2080,400)(0,-1)[\phantom{Y_2}``=]{350}1l
\put(780,210){\fbox{$\tilde l^M_E$}}
\efig}
$$

$$%\begin{equation} \eqlabel{effects}
\scalebox{0.8}{
\bfig
 \putmorphism(-210,400)(1,0)[ME`M(IE) `M\tilde\lambda^{-1}_{E}]{670}1a
 \putmorphism(460,400)(1,0)[\phantom{A\ot B}`(MI)E `\tilde\alpha^{-1}_{M,I,E}]{800}1a
 \putmorphism(1260,400)(1,0)[\phantom{A\ot B}`ME`\rho^\Mm_ME]{800}1a
\putmorphism(-150,50)(1,0)[ME`ME.`=]{2230}1b

\putmorphism(-180,400)(0,-1)[\phantom{Y_2}``=]{350}1r
\putmorphism(2080,400)(0,-1)[\phantom{Y_2}``=]{350}1l
\put(780,210){\fbox{$\tilde m^M_E$}}
\efig}
$$
The action can also be defined if $\tilde\alpha, \tilde\lambda$ are not invertible. We will use this in \seref{bistrong,comm} where instead we will consider 
horizontal transformations (resp. vertical strict transformations) $\beta: -\times(-\times-)\Rightarrow (-\ot-)\times-$ and $\nu:-\Rightarrow I\times-$ 
and modifications with 2-cell components of the form as in this remark. This kind of action we will call {\em lax action}. 
More about it we will comment in \ssref{equiv}. The proper actions, as in the two definitions above, we will use in \seref{Para}. 
\end{rem}

\begin{ex}
If $\Mm$ is the trivial double category, the above vertical and horizontal action of double categories are precisely the vertical and horizontal double 
monad from \cite[Definitions 7.3 and 7.2]{GGV}, respectively. 
\end{ex}

We will often write $\crta\ot$ for an action double functor $F:\Dd\times\Ee\to\Ee$, be it for vertically or horizontally 
monoidal double categories $\Dd$ and $\Ee$. Moreover, in the diagrams we will abbreviate the notation $F(M,E)=M\crta\ot E$ into 
simply $ME$, as we already did in the above definition.

\smallskip

\begin{rem} \rmlabel{gen action}
Observe that the action of the multiplication $\mu_A: T^2(A)\to T(A)$ and the unit $\eta_A:A\to T(A)$ of a double categorical monad $T$ 
from \cite{GGV} (both vertical and horizontal) in the setting of an action of a monoidal double category $\Mm$ become 
$$\tilde\alpha^{-1}_{M,N,A}: M\crta\ot(N\crta\ot A)\to(M\ot N)\crta\ot A \qquad \text{and} \qquad \tilde\lambda^{-1}_A:A\to J\crta\ot A$$
where $J\in\Mm$ is the monoidal unit. This is the reason why we will use prevalently $\beta:=\tilde\alpha^{-1}$ rather  than $\tilde\alpha$. 
\end{rem}

One has by \prref{lifting 1v to equiv}:

\begin{prop}  \cite[Proposition 10.16]{Fem3} \prlabel{vert->horiz act} \\
Given a vertical action $(F, \tilde\alpha, \tilde\lambda)$ of a monoidal double category 
$(\Dd, \ot, \alpha,\lambda,\rho)$ by a pseudodouble functor $F$ on a double category $\Ee$, so that $\alpha,\lambda,\rho$ and 
$\tilde\alpha, \tilde\lambda$ are liftable, then $F$ induces a horizontal action $(F, \hat{\tilde\alpha}, \hat{\tilde\lambda})$ of the horizontally monoidal double category $(\Dd, \ot, \hat\alpha, \hat\lambda,\hat\rho)$ on $\Ee$. 
\end{prop}

\subsection{Monoidal actions of double categories}

Monoidal actions of double categories require the existence of a braiding in the monoidal double category that is acting on another 
double category. Again, we differentiate vertical and double braidings, we start by introducing their definitions. 

\begin{defn}
A {\em braiding} in a  vertically monoidal double category $\Mm$ is an invertible vertical strict transformation 
$\Phi:-\ot-\Rightarrow (-\ot-)\tau: \Mm\times\Mm\to\Mm\times\Mm$, where $\tau$ is the flip double functor (flipping the order of the factors), 
with two identity vertical modifications $b_1, b_2$ with the following 2-cell components
$$
\scalebox{0.9}{
\bfig
 \putmorphism(-90,500)(1,0)[(AB)C` (AB)C`=]{440}1a

\putmorphism(-100,500)(0,-1)[` (BA)C`\Phi_{A,B}C]{400}1l
\putmorphism(-100,110)(0,-1)[\phantom{Y_2}`B(AC)`\alpha_{B,A,C}]{400}1l
\putmorphism(-100,-300)(0,-1)[\phantom{Y_2}``B\Phi_{A,C}]{400}1l

\putmorphism(350,500)(0,-1)[` A(BC)`\alpha_{A,B,C}]{400}1r
\putmorphism(350,110)(0,-1)[\phantom{Y_2}`(BC)A`\Phi_{A,BC}]{400}1r
\putmorphism(350,-300)(0,-1)[\phantom{Y_2}``\alpha_{B,C,A}]{400}1r

\putmorphism(-90,-710)(1,0)[B(CA)` B(CA)`=]{440}1a
\put(60,-130){\fbox{$b_1$}}
\efig}
\qquad
\scalebox{0.9}{
\bfig
 \putmorphism(-90,500)(1,0)[A(BC)` A(BC)`=]{440}1a

\putmorphism(-100,500)(0,-1)[` A(CB)`A\Phi_{B,C}]{400}1l
\putmorphism(-100,110)(0,-1)[\phantom{Y_2}`(AC)B`\alpha_{A,C,B}^{-1}]{400}1l
\putmorphism(-100,-300)(0,-1)[\phantom{Y_2}``\Phi_{A,C}B]{400}1l

\putmorphism(350,500)(0,-1)[` (AB)C `\alpha_{A,B,C}^{-1}]{400}1r
\putmorphism(350,110)(0,-1)[\phantom{Y_2}`C(AB)`\Phi_{AB,C}]{400}1r
\putmorphism(350,-300)(0,-1)[\phantom{Y_2}``\alpha_{C,A,B}^{-1}]{400}1r

\putmorphism(-90,-710)(1,0)[(CA)B` (CA)B.`=]{440}1a
\put(60,-130){\fbox{$b_2$}}
\efig}
$$
\end{defn}

This braiding we will refer to as a {\em vertical double braiding}, and for $\Mm$ we will say that it is (vertically) double braided. 
Similarly we will do in the horizontal setting:

\begin{defn}
A {\em braiding} in a horizontally monoidal double category $\Mm$ is a horizontal pseudonatural transformation 
$\Phi:-\ot-\Rightarrow (-\ot-)\tau: \Mm\times\Mm\to\Mm\times\Mm$ with invertible 1h-cell components and 
with two invertible horizontal modifications $b_1, b_2$ with the following 2-cell components
$$ 
\scalebox{0.82}{
\bfig
\putmorphism(-390,450)(0,-1)[\phantom{Y_2}` `=]{450}1r 

\putmorphism(-390,450)(1,0)[(AB)C` A(BC) `\alpha_{A,B,C}]{740}1a 
\putmorphism(330,450)(1,0)[\phantom{A''\ot B'}` (BC)A ` \Phi_{A,BC} ]{740}1a %%%%%
 \putmorphism(1020,450)(1,0)[\phantom{A''\ot B'}` B(CA) `\alpha_{B,C,A} ]{760}1a %T(T_A(T_BT_C))

 \putmorphism(-390,0)(1,0)[(AB)C` (BA)C `\Phi_{A,B}C]{740}1a 
\putmorphism(340,0)(1,0)[\phantom{A''\ot B'}` B(AC) ` \alpha_{B,A,C} ]{740}1a %%%%%
 \putmorphism(1040,0)(1,0)[\phantom{A''\ot B'}` B(CA) `B\Phi_{A,C} ]{740}1a %T(T_A(T_BT_C))

\putmorphism(1820,450)(0,-1)[\phantom{Y_2}``=]{450}1r 

\put(700,210){\fbox{$b_1$}}
\efig}
$$

$$ 
\scalebox{0.82}{
\bfig
\putmorphism(-390,450)(0,-1)[\phantom{Y_2}` `=]{450}1r 

\putmorphism(-390,450)(1,0)[A(BC)` (AB)C `\alpha_{A,B,C}^{-1}]{740}1a 
\putmorphism(370,450)(1,0)[\phantom{A''\ot B'}` C(AB) ` \Phi_{AB,C} ]{700}1a %%%%%
 \putmorphism(1020,450)(1,0)[\phantom{A''\ot B'}` (CA)B `\alpha_{C,A,B}^{-1} ]{760}1a %T(T_A(T_BT_C))

 \putmorphism(-390,0)(1,0)[A(BC)` A(CB) `A\Phi_{B,C} ]{740}1a 
\putmorphism(330,0)(1,0)[\phantom{A''\ot B'}` (AC)B ` \alpha_{A,C,B}^{-1} ]{740}1a %%%%%
 \putmorphism(1050,0)(1,0)[\phantom{A''\ot B'}` (CA)B. `\Phi_{A,C}B ]{730}1a %T(T_A(T_BT_C))

\putmorphism(1820,450)(0,-1)[\phantom{Y_2}``=]{450}1r 

\put(700,230){\fbox{$b_2$}}
\efig}
$$
\end{defn}

A braided double category yields a premonoidal double category in the following way.

\begin{prop} \prlabel{ex-braid-prem}
Let $(\Mm, \ot, \Phi)$ be a horizontally braided monoidal double category and let $M\in\Mm$ be a fixed object. We define a double category 
${}_M\Mm$ having the same objects as $\Mm$, its 1-cells in both directions on objects $A\to B$ and 2-cells on them are given by the corresponding 1-cells on objects $M\ot A\to M\ot B$ in $\Mm$ and 2-cells in $\Mm$ over those, respectively, while all compositions and identities are as in $\Mm$. Then ${}_M\Mm$ is a (horizontally) premonoidal double category and the canonical pseudodouble functor 
$M\ot-:\Mm\to {}_M\Mm$ preserves the premonoidal structure strictly. 
\end{prop}

\begin{proof}
We only give a sketch of the proof. We set $\gamma_{M,A,B}:=\alpha_{A,M,B}(\Phi_{M,A}\ot\Id_B)\alpha_{M,A,B}^{-1}$ and define a 
binoidal structure on ${}_M\Mm$ by 
$$
A\ltimes g:=\big(M(A B)\stackrel{\gamma_{M,A,B}}{\to}A(MB) \stackrel{ A g}{\to}A(MB') \stackrel{\gamma_{M,A,B}}{\to} M(AB')\big)
$$
$$f\rtimes B:=\big(M(AB)\stackrel{\alpha_{M,A,B}^{-1}}{\to}(MA) B \stackrel{fB}{\to}(MA')B \stackrel{\alpha_{A,M,B}}{\to} M(A'B)\big).$$
%By the way how ${}_M\Mm$ is defined, it is impossible to express naturality of the associativity constraint in more than one variable at a time. 
One defines the 2-cell components of the three $\alpha$'s in ${}_M\Mm$ as suitable horizontal compositions of the corresponding 2-cell components of the associativity $\alpha$ of $\Mm$ and the identity 2-cells on 1-cell components of $\alpha$ and $\Phi\ot 1, 1\ot\Phi$, 
necessary to accommodate the domain and codomain of 1h-cells in ${}_M\Mm$. This is a sort of ``horizontal conjugation''. 
Since $\Mm$ is monoidal, the canonical pseudodouble functor 
$M\ot-:\Mm\to {}_M\Mm$ factors through the center double category of ${}_M\Mm$. In particular, 1h-cell components of the structural 
equivalences are all central. For the last claim, observe that for example for $g:B\to B'$ in $\Mm$ the 2-cell component $\alpha_{A,g,C}$ 
in $\Mm$ is mapped into $M\ot\alpha_{A,g,C}$, which corresponds to the horizontal conjugation of $\alpha_{A,Mg,C}$ defining the corresponding 
structural equivalence of ${}_M\Mm$.  
\qed\end{proof}

For vertically braided $\Mm$ we will consider: 

\begin{prop} \prlabel{braid-mon.act.}
Let $\Mm$ be a vertically monoidal double category that is also double braided (with a double braiding $\Phi$). Then the monoidal structure 
$\ot:\Mm\times\Mm\to\Mm$ is a lax monoidal double functor. 
%\item Additionally, given a vertical action $F:\Mm\times\Dd\to\Dd$, with vertical structure transformations $\tilde\alpha, \tilde\lambda$ 
%and an interchanger $int: (M\rtr A)\ot(N\rtr B)\to (MN)\rtr(AB)$. Then $F$ is lax monoidal (with $F^2$ given by the interchanger and $F^0$ by 
%$\tilde\lambda_I$). 
%\end{enumerate} \vspace{-0,2cm}
\end{prop}

\begin{proof}
The 2-dimensional 2-cocycle determining $\ot^2$ is given by the 2-cell components $\Id_f\ot\Phi_{g,h}\ot\Id_l$ at 1h-cells $f,g,h,l$. 
The cocycle (associativity) axiom holds by the interchange law of double categories and left and right normalization (unitalities) 
hold trivially (identify the 1v-cell component $\ot^0:I\to I\ot I$ with $1^I\ot 1^I$). 
\qed\end{proof}

The next definition is inspired by 1-categorical \cite[Definition 5.1.1]{CG}. 
%\textcolor{rojo}{Check if ``action is a morphism in $\MonCat$ or similar - the category of monoidal double categories???}. 

\begin{defn} \delabel{v mon act}
Let $\Mm$ and $\Dd$ be vertically monoidal double categories endowed with a vertical action $F:\Mm\times\Dd\to\Dd$ 
with structural vertical transformations $\tilde\alpha, \tilde\lambda$, and assume that $\Mm$  is double braided. 
We say that $(F, \tilde\alpha, \tilde\lambda)$ is a {\em lax monoidal (vertical) action} if: \vspace{-0,2cm}
\begin{enumerate}
\item $F$ is a lax monoidal pseudodouble functor, and \vspace{-0,2cm}
\item the transformations $\tilde\alpha, \tilde\lambda$ 
are pseudomonoidal vertical transformations (see \cite[Definition 5.3]{GGV}). %and \vspace{-0,2cm}
%\item the vertical modifications $\tilde p, \tilde l, \tilde m$ are monoidal (\cite[Definition 5.5]{GGV}). 
\end{enumerate}
\end{defn}

\begin{rem}
For the  transformation $\tilde\alpha: (-\ot-)\crta\ot-\Rightarrow -\crta\ot(-\crta\ot-)$ to be pseudomonoidal it is necessary that 
$\Mm$ be braided. We will illustrate this in more detail below for the horizontal case. 
\end{rem}

%That  $F=\crta\ot:\Mm\times\Dd\to\Dd$  in the above definition is lax monoidal as a pseudodouble functor, according to 
%\cite[Definition 5.1]{GGV} it means that 
As a lax monoidal pseudodouble functor $F=\crta\ot:\Mm\times\Dd\to\Dd$ in the above definition 
is in particular equipped with a vertical transformation $F^2:F(-)\ot F(-)\Rightarrow F(-\ot-)$ with 1v-cell components - 
that we call {\em interchangers} - 
$$int_{M,D,N,E}: (M\crta\ot D)\ot(N\crta\ot E)\to (M\ot N)\crta\ot(D\ot E)$$ 
and 2-cells  
$$
\scalebox{0.86}{
\bfig
\putmorphism(-100,50)(1,0)[(M\crta\ot D)\ot(N\crta\ot E)`(M'\crta\ot D')\ot(N'\crta\ot E')`(m\crta\ot d)\ot(n\crta\ot e)]{1650}1a
\putmorphism(-100,-400)(1,0)[(M\ot N)\crta\ot(D\ot E)`(M'\ot N')\crta\ot(D'\ot E'). `(m\ot n)\crta\ot(d\ot e)]{1650}1a
\putmorphism(-80,50)(0,-1)[\phantom{Y_2}``int_{M,D,N,E}]{420}1l
\putmorphism(1560,50)(0,-1)[\phantom{Y_2}``int_{M',D',N',E'}]{420}1r
\put(500,-140){\fbox{$int_{m,d,n,e}$}}
\efig}
$$
%and an additional 1v-cell $F^0:I_\Dd\to I_\Mm\crta\ot I_\Dd$, so that \textcolor{rojo}{+ the other + axioms.}

For horizontally braided $\Mm$ we will consider the following definition (instead of structural horizontal transformations $\tilde\alpha, \tilde\lambda$ we will 
use $\beta, \nu$ for future use, as explained in \rmref{monad action}).

\begin{defn} \delabel{h mon act}
For horizontally monoidal double categories $\Mm$ and $\Dd$ endowed with a horizontal action $F:\Mm\times\Dd\to\Dd$ 
with structural horizontal transformations $\beta, \nu$, we say that $(F, \beta, \nu)$ is a {\em lax monoidal (horizontal) action} if \vspace{-0,2cm}
\begin{enumerate}
\item $\Mm$ is double braided (with braiding $\Phi$); 
\item $F$ is horizontally lax monoidal pseudodouble functor, \vspace{-0,2cm} 
\item the transformations $\beta, \nu$ are monoidal horizontal transformations %:=\tilde\alpha^{-1}   :=\tilde\lambda^{-1}
so that the structure modifications $\beta^2, \beta^0$ and $\nu^2, \nu^0$ are invertible horizontal modifications, and \vspace{-0,2cm} 
\item the horizontal modifications $\tilde p, \tilde l, \tilde m$ are monoidal (\cite[Definition 5.5]{GGV}).  
\end{enumerate}
\end{defn}

\smallskip

Horizontally lax monoidal structures of the pseudoduble functors $G=-\crta\ot(-\crta\ot-)$ and $\tilde G=(-\ot-)\crta\ot-$ at 
objects $(M,N,A),(M', N',B)\in\Mm\times\Mm\times\Dd$ are induced by \vspace{-0,2cm}
$$G^2=\big( (M(NA))(M'(N'B))\stackrel{F^2_{M,NA,M', N'B}}{\longrightarrow} (MM')((NA)(N'B))\stackrel{(MM')F^2_{N,A,N',B}}{\longrightarrow}
(MM')((NN')(AB)) \big)$$ %_{M,N,A;M', N',B}
$$\tilde G^2=\big(((MN)A)((M'N')B) \stackrel{F^2_{MN,A,M'N',B}}{\longrightarrow} ((MN)(M'N'))(AB)
\stackrel{(M\Phi_{N,M'}N')(AB)}{\longrightarrow}((MM')(NN'))(AB) \big).$$
Thus monoidality of the transformation $\beta$, concretely the 2-cocycle 2-cells $\beta^2$, require $\Mm$ to be braided.

\begin{ex}
If $\Mm$ is the trivial double category, the above vertically and horizontally monoidal action of double categories are precisely the pseudomonoidal 
vertical and monoidal horizontal double monad from \cite[Definitions 8.3 and 8.2]{GGV}, respectively. 
\end{ex}

\smallskip

\begin{rem} \rmlabel{hor modif of hor lax mon}
Both in the case of a vertical and a horizontal action, monoidal units of $\Mm$ and $\Dd$ we will denote as $J$ and $I$, respectively. 
The 1h-cell (and also 1v-cell) $F^0:I\to F(J\times I)$ we will denote by $\iota: I\to JI$. 

Now we make several observations regarding the above definition. 
Firstly, horizontal lax monoidality of $F$ according to \rmref{hor lax mon f} means that it is given by a lax monoidal 
bicategorical pseudofunctor $\HH(F):\HH(\Mm)\times\HH(\Dd)\to\HH(\Dd)$. In particular, 
there are horizontal interchangers, {\em i.e.} 1h-cells 
$F^2_{M,D,N,E}=int_{M,D,N,E}: (M\crta\ot D)\ot(N\crta\ot E)\to (M\ot N)\crta\ot(D\ot E)$, and there are horizontal modifications 
$\omega, \delta, \gamma$ as mentioned in the Remark.  

Secondly, the 2-cell components of the 
invertible horizontal modifications $\beta^2, \beta^0$ and $\nu^2, \nu^0$ have the following form
$$ 
\scalebox{0.82}{
\bfig
\putmorphism(-390,450)(0,-1)[\phantom{Y_2}` `=]{450}1r 

\putmorphism(-390,450)(1,0)[(M(NA))(P(QB))` (MP)((NA)(QB))`F^2_{M,NA,P,QB}]{1170}1a 
\putmorphism(1000,450)(1,0)[\phantom{A''\ot B'}`(MP)((NQ)(AB)) ` (MP)F^2_{N,A,Q,B}]{1170}1a%%%%%
 \putmorphism(2370,450)(1,0)[\phantom{A''\ot B'}` ((MP)(NQ))(AB) `\beta_{MP,NQ,AB}]{1160}1a %T(T_A(T_BT_C))

 \putmorphism(-390,0)(1,0)[(M(NA))(P(QB))` ((MN)A)((PQ)B) `\beta_{M,N,A}\beta_{P,Q,B}]{1350}1b 
 \putmorphism(1180,0)(1,0)[\phantom{A''\ot B'}` ((MN)(PQ))(AB) `F^2_{MN,A,PQ,B}]{1110}1b 
 \putmorphism(2500,0)(1,0)[\phantom{A''\ot B'}` ((MP)(NQ))(AB) `F^2_{M,N,P,Q}(AB)]{1180}1b 

\putmorphism(3600,450)(0,-1)[\phantom{Y_2}``=]{450}1r 

\put(1100,230){\fbox{$(\beta^2)^{M,N,P,Q}_{A,B}$}}
\efig}
$$

$$ 
\scalebox{0.82}{
\bfig
\putmorphism(-390,450)(0,-1)[\phantom{Y_2}` `=]{450}1r 

\putmorphism(-390,450)(1,0)[I` JI`\iota]{370}1a 
\putmorphism(-150,450)(1,0)[\phantom{A''\ot B'}` J(JI) ` J\iota ]{650}1a %%%%%
 \putmorphism(400,450)(1,0)[\phantom{A''\ot B'}` (JJ)I `\beta_{MP,NQ,AB}]{900}1a %T(T_A(T_BT_C))
 \putmorphism(1230,450)(1,0)[\phantom{A''\ot B'}` JI `\lambda^\Mm_{J}]{650}1a %T(T_A(T_BT_C))

 \putmorphism(-390,0)(1,0)[I` JI `\iota]{2280}1b 

\putmorphism(1900,450)(0,-1)[\phantom{Y_2}``=]{450}1r 

\put(700,230){\fbox{$\beta^0$}}
\efig}
$$

$$ 
\scalebox{0.82}{
\bfig
\putmorphism(-390,450)(0,-1)[\phantom{Y_2}` `=]{450}1r 

\putmorphism(-390,450)(1,0)[AB` (JA)(JB)`\nu_A\nu_B]{570}1a 
\putmorphism(210,450)(1,0)[\phantom{A''\ot B'}` (JJ)(AB) ` F^2_{J,A,J,B} ]{770}1a %%%%%
 \putmorphism(1000,450)(1,0)[\phantom{A''\ot B'}` J(AB) `\lambda^\Mm(AB)]{900}1a %T(T_A(T_BT_C))

 \putmorphism(-390,0)(1,0)[AB` J(AB) `\nu_{AB}]{2270}1b 

\putmorphism(1900,450)(0,-1)[\phantom{Y_2}``=]{450}1r 

\put(620,230){\fbox{$(\nu^2)_{A,B}$}}
\efig}
$$

$$ 
\scalebox{0.82}{
\bfig
\putmorphism(-390,450)(0,-1)[\phantom{Y_2}` `=]{450}1r 
\putmorphism(-390,450)(1,0)[I` JI`\nu_I]{800}1a 
\putmorphism(-390,0)(1,0)[I` JI. `\iota]{800}1b 
\putmorphism(400,450)(0,-1)[\phantom{Y_2}``=]{450}1r 
\put(-50,230){\fbox{$\nu_0$}}
\efig}
$$
According to \rmref{hor lax mon f}, the monoidal horizontal transformations $\beta, \nu$ are monoidal pseudonatural transformations 
in the underlying bicategories satisfying additionally the law \equref{h mon tr law}. 

Due to the same remark, $\tilde p, \tilde l, \tilde m$ are precisely the bicategorical modifications with 2-cell components 
in the underlying bicategory $\HH(\Dd)$.
\end{rem}

\bigskip

Given a setting of \deref{v mon act}, if we assume that the following are liftable: 
\begin{itemize}
\item the 1v-cell components of the vertical braiding; 
\item the monoidal structures $\alpha^\Mm,\lambda^\Mm,\rho^\Mm$ of $\Mm$ and $\alpha^\Dd,\lambda^\Dd,\rho^\Dd$ of $\Dd$; 
 \vspace{-0,2cm}
\item the vertical transformations ({\em i.e.} their 1v-cell components) $\tilde\alpha, \tilde\lambda$ and 
$F^2, F^0$ of the vertical action; and that these four are also invertible vertical transformations,
\end{itemize}
then we obtain:
\begin{itemize}
\item $\Mm$ and $\Dd$ are horizontally monoidal (according to \thref{two monoidalities}) 
(and their underlying horizontal bicategories $\HH(\Mm)$ and $\HH(\Dd)$ are monoidal);  \vspace{-0,2cm}
\item $\Mm$ is horizontally braided (by \prref{lifting 1v to equiv} and \prref{essence});  \vspace{-0,2cm}
\item a horizontal lax monoidal action $F:\Mm\times\Dd\to\Dd$ (due to \prref{vert->horiz act} and \prref{lifting 1v to equiv})
(and hence an action of the underlying horizontal bicategories that is a lax monoidal pseudofunctor);  \vspace{-0,2cm}
\item the associator $\hat{\tilde\alpha}$ and the unitor $\hat{\tilde\lambda}$ of the horizontal action $F:\Mm\times\Dd\to\Dd$ are monoidal 
horizontal transformations (due to \cite[Proposition 5.7]{GGV}), and hence monoidal transformations on the underlying horizontal bicategories;  \vspace{-0,2cm}
\item the monoidal horizontal modifications $\widehat{\tilde p}, \widehat{\tilde l}, \widehat{\tilde m}$ (due to \prref{essence}). 
%and hence monoidal transformations on the underlying horizontal bicategories. 
\end{itemize}

\begin{cor} \colabel{v.act-h.act}
Let $\Mm$ and $\Dd$ be vertically monoidal double categories endowed with a lax monoidal vertical action $F:\Mm\times\Dd\to\Dd$ 
with structural vertical transformations $\tilde\alpha, \tilde\lambda$ and $F^2, F^0$, and assume that $\Mm$ is braided by 
a vertical braiding $\Phi$.    
If the monoidal structures $\alpha^\Mm,\lambda^\Mm,\rho^\Mm$ of $\Mm$ and $\alpha^\Dd,\lambda^\Dd,\rho^\Dd$ of $\Dd$, and  
the vertical transformations $\tilde\alpha, \tilde\lambda, \,\, F^2, F^0, \,\, \Phi$ are all invertible and liftable, then 
$(F, \hat{\tilde\alpha}, \hat{\tilde\lambda}, \hat{F^2}, \hat{F^0})$ determines a lax monoidal horizontal action of $\Mm$ on $\Dd$. 
\end{cor}

\subsubsection{The perspective of locally cubical bicategories} \ssslabel{perspective lcb}

Let $\V\w\MonDbl_v$ denote the tricategory whose objects are vertically monoidal double categories and whose hom-objects are 
the 2-categories ${}^v\w\MMonDbl_v[\Cc, \Dd]_2$ (recall from before \coref{2-cat lift}). 
We refer to it as the underlying vertical tricategory of $\textfrak{VMonDbl}$. 

Similarly, let us denote by $\HH\w\MonDbl_h$ the underlying horizontal tricategory of $\textfrak{HMonDbl}$. 

Then we can say that a vertical action (from \deref{v mon act}) is an action in the tricategory $\V\w\MonDbl_v$. Likewise, 
we can say that a horizontal action (from \deref{h mon act}) is an action in the tricategory $\HH\w\MonDbl_h$. 

\medskip

%In \coref{2-cat lift} $\MonDbl_{v^0}^v[\Cc, \Dd]_2$ denoted the sub-2-category of $\MonDbl_v^v[\Cc, \Dd]_2$ whose 1-cells are liftable.
Let $\V^0\w\MonDbl_{v^0,0}$ denote the sub-tricategory of $\V\w\MonDbl_v$ in which the vertically monoidal structures of objects, the 
lax monoidality 1-cocycles, as well as the 2-cells (the vertical transformations) are invertible and liftable. Then the claim of \coref{v.act-h.act} 
comes as a result of the mapping an action in $\V^0\w\MonDbl_{v^0,0}$ by the composition 
$$\V^0\w\MonDbl_{v^0,0} \hspace{0,4cm}\stackrel{L_{tr,0}^0}{\longrightarrow}\hspace{0,4cm} \V^0\w\MonDbl_{v^0,0} 
\hspace{0,4cm}\stackrel{L_{v^0}}{\longrightarrow}\hspace{0,4cm} \HH\w\MonDbl_h$$
whereby $L_{tr,0}^0$ is the corresponding restriction of \equref{auto}, and $L_{v^0}$ is the corresponding restriction of $L$ 
from \prref{ver-hor-cbl}. From the point of view of the locally cubical bicategories we can draw this graphically as follows 
\begin{figure}[H] 
\includegraphics[scale=0.8]{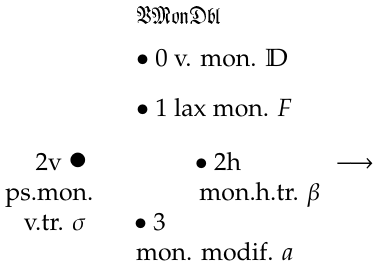} \hspace{0,14cm}
\includegraphics[scale=0.8]{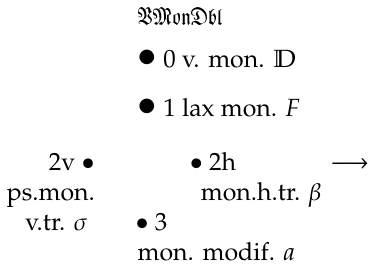} \hspace{0,14cm}
\includegraphics[scale=0.8]{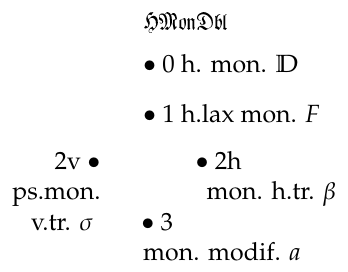} 
\end{figure}
where enlarged bullets represent the cells which are assumed to be liftable.

\section{Bistrong and commutative actions} \selabel{bistrong,comm}

In \cite[Section 10.2]{Fem3} we introduced the notion of a strength on a double monad (it appeared implicitly in \cite{GGV}). 
In this section we generalize that strength to a notion of a strength for an action of a monoidal double category, following the line of 
``strength-like morphisms'' from \cite[Section 5]{CG} for actegories. We also generalize the notion of a 
commutative monad from bicategories from \cite{HP, HF2} to double categories.

\subsection{Bistrong and commutative actions} %\sslabel{bistrong,comm}

We build the necessary environment to introduce the notion of a bistrength for an action of a monoidal double category: we will need strengths from two sides. 

\begin{defn}
Let $\Mm$ and $\Dd$ be vertically monoidal double categories with a vertical action $\crta\ot:\Mm\times\Dd\to\Dd$. 
A {\em left strength} of $\crta\ot$ is an invertible vertical strict transformation with 1v-cell components 
$$t_{A,B}^M: A\ot(M\crta\ot B)\to M\crta\ot(A\ot B)$$ %M:(M\crta\ot A)
for $M\in\Mm, A,B\in\Dd$,  and 
identity vertical modifications with the following identity vertically globular 2-cell components, expressing compatibility of the strength 
$t$ with:
 \begin{enumerate}[a)]
\item monoidal structure of $\Dd$ \vspace{-0,2cm}
$$x^M_A: \lambda_{M\crta\ot A}\Rightarrow\frac{t_{I,A}^M}{ M\crta\ot\lambda_A} \qquad \qquad
y_{A,B,C}^M: \frac{t_{AB,C}^M}{M\crta\ot\alpha_{A,B,C}}\Rightarrow\threefrac{\alpha_{A,B,M\crta\ot C}}{At_{B,C}^M}{ t_{A,BC}^M};$$
\item action structure of $\crta\ot$ \vspace{-0,2cm} %identity vertical modifications with identity vertically globular 2-cell components 
$$w_{A,B}^{M, N}: \frac{A\ot\tilde\alpha_{M,N,B}^{-1}}{t^{M\ot N}_{A,B}}\Rightarrow
\threefrac{t_{A,N\crta\ot B}^M}{M\crta\ot t^N_{A,B}}{\tilde\alpha_{M,N,A\ot B}^{-1}} \qquad \qquad
z_{A,B}:\tilde\lambda_{A\ot B}^{-1}\Rightarrow\frac{A\ot\tilde\lambda_B^{-1}}{t_{A,B}^I}.$$
\end{enumerate}

A {\em right strength} of $\crta\ot$ is defined similarly, with 1v-cell components 
$$s_{A,B}^M:(M\crta\ot A)\ot B\Rightarrow M\crta\ot(A\ot B)$$ %A\ot (M\crta\ot B)
and modifications that we denote by $x', y', z', w'$. 

For strengths on a vertical action we will also say that they are {\em vertical strengths}. 

We say that the action $\crta\ot$ with left and right vertical strengths $t$ and $s$ is {\em bistrong} if there is an identity vertical modification $q$ with components 
$$q_{A,B,C}^M:\threefrac{\alpha_{A,M\crta\ot B,C}}{As_{B,C}^M}{ t_{A,BC}^M} \Rightarrow\threefrac{t_{A,B}^M\ot C}{s_{AB,C}^M}{M\crta\ot\alpha_{A,B,C}}$$
where $\alpha$ is the associativity of $\Dd$. 
\end{defn}

%\begin{defn}
%Given a vertical action $\crta\ot:\Mm\times\Dd\to\Dd$. By an {\em interchanger} we mean a strict vertical invertible transformation with 1v-cell 
%components $int: (M\crta\ot A)\ot(N\crta\ot B)\to (MN)\crta\ot(AB)$ (we omit the tensor product of $\Mm$) \textcolor{rojo}{that satisfies the following 5 axioms: ...}
%\end{defn}

%\begin{prop}
%Let $\Mm$ and $\Dd$ be vertically monoidal double categories with a pseudodouble functor $F:\Mm\times\Dd\to\Dd$. 
%Assuming that $\Mm$ is braided (with double braiding $\Phi$) and there exists an interchanger $int$, then 
%$F$ has is lax monoidal structure. 
%\end{prop}

\bigskip

The above would be a {\em vertical bistrength}. As several times so far, 
we will differentiate two more instances of strengths and accordingly of bistrengths. 
These are so-to-say ``square'' or ``mixed'' ones, and horizontal ones, where the vertical notions induce the other two, assuming that the corresponding strengths are given by invertible and liftable vertical transformations. We are not going to define explicitly the mixed versions of the notions, as we will only need them in passing to obtain other results. We only highlight that mixed strengths are 
given for a vertically monoidal $\Mm$ by horizontal transformations with respect to the vertically monoidal structure of $\Dd$ so that 
their structure modifications are mixed ones, {\em i.e.} of the form as in \deref{modif-hv}. We will denote these modifications by 
$x^*, y^*, w^*, z^*$ (for left strengths) and 
$x'^*, y'^*, w'^*, z'^*$ (for right strengths) alluding to the modifications obtained in \prref{omega*}, see also  
\equref{omega-gen}. They are required to obey analogous axioms as in the case of horizontal strengths, whose axioms 
coincide with those of the bicategorical strengths (as only globular cells appear in them). A mixed bistrength is given by a mixed 
modification $q^*$ that satisfies two axioms analogous to those of \cite[Figure 6, p.9]{HF2} 
(they stand for compatibility of $q^*$ with $z^*, z'^*$ and with $w^*, w'^*$).

\small

%Assume that $t$ and $s$ are a left and a right vertical strength on a vertical action $\crta\ot:\Mm\times\Dd\to\Dd$ for vertically monoidal double categories $\Mm$ and $\Dd$.  
%Given a monoidal double category $(\Dd,\ot,\alpha,\lambda,\rho)$ and a vertical double monad $T$ on $\Dd$. 
Vertical strengths $t,s$ induce mixed strengths $\hat t, \hat s$ provided $t,s$ and the monoidal structure of $\Dd$ 
%, the vertical action $(\crta\ot,\tilde\alpha, \tilde\lambda)$ and vertical strengths $t,s$ as in the definition above 
are liftable. %According to  \prref{vert->horiz act}  then $\crta\ot$ 
%induces a horizontal action $(\crta\ot, \hat{\tilde\alpha}, \hat{\tilde\lambda})$ of the horizontally monoidal double category 
%$(\Mm, \ot, \hat\alpha^\Mm, \hat\lambda^\Mm,\hat\rho^\Mm)$ on $\Dd$ with horizontal strengths $\hat t, \hat s$.  
%Given a horizontally monoidal double category $(\Dd,\ot,\alpha,\lambda,\rho)$ and a horizontal double monad $S$ on $\Dd$. Assume that $t$ and $s$ are a left and a right horizontal strength on $S$. We 
Given a vertical bistrength with an identity vertical modification $q$, it induces a mixed modification with components 
$$
\scalebox{0.86}{
\bfig
 \putmorphism(-620,500)(1,0)[(A(M\crta\ot B))C`(M\crta\ot(AB))C`\hat t_{A,B}^M\ot C]{950}1a
 \putmorphism(500,500)(1,0)[\phantom{F(f)}`M\crta\ot((AB)C) `\hat s_{AB,C}^M]{700}1a
 \putmorphism(-620,120)(1,0)[A((M\crta\ot B)C) ` AM\crta\ot (BC) `A\ot\hat s_{B,C}^M]{950}1a
 \putmorphism(460,120)(1,0)[\phantom{G(B)}`M\crta\ot(A(BC)) `\hat t_{A,BC}^M]{740}1a
\putmorphism(-580,500)(0,-1)[\phantom{Y_2}``\alpha]{380}1l
\putmorphism(1110,500)(0,-1)[\phantom{Y_2}``M(\alpha)]{380}1r
\put(210,310){\fbox{$q_{A,B,C}^{*M}$}}
\efig}
$$
whereby the $\alpha$'s on the vertical arrows are from $\Mm$. 

\medskip

\begin{defn} \delabel{hor-strength}
Let $\Mm$ and $\Dd$ be horizontally monoidal double categories with a horizontal action $\crta\ot:\Mm\times\Dd\to\Dd$. 
A {\em left (horizontal) strength} on $\crta\ot$ consists of:
\begin{enumerate} [a)]
\item a horizontal pseudonatural transformation with 1h-components \\ $t^M_{A,B}:A\ot M(B)\to M(A\ot B)$ for $M\in\Mm, A,B\in\Dd$; 
\item invertible horizontal modifications with horizontally globular 2-cell components 
$$x_A^M: [t^M_{I,A}\vert M\lambda_A]\Rightarrow\lambda_{M\crta\ot A} \qquad 
y^M_{A,B,C}:[t^M_{AB,C}\vert M(\alpha_{A,B,C})] \Rightarrow[\alpha_{A,B,M\crta\ot C}\vert At^M_{B,C}\vert t^M_{A,BC}]$$
which fulfill three axioms as in \cite[Figure 3, p. 6-7]{HF2} (connecting them to the horizontal modifications $p,m,l$ %$\pi,\mu,\lambda$ 
of the horizontally monoidal double category $\Dd$);
\item invertible horizontal modifications with horizontally globular 2-cell components 
$$w^{M,N}_{A,B}: [t^M_{A,N\crta\ot B}\vert M\crta\ot t^N_{A,B}\vert \tilde\alpha^{-1}_{M,N,A\ot B}]\Rightarrow[A\tilde\alpha^{-1}_{M,N,B}\vert 
t^{M\ot N}_{A,B}] \qquad 
z_{A,B}:[A\tilde\lambda^{-1}_B\vert t^J_{A,B}]\Rightarrow\tilde\lambda^{-1}_{A\ot B}$$
which fulfill seven axioms analogous to those of \cite[Figures 4 and 5, p. 7-8]{HF2} 
(three axioms for connecting $w$ and $z$ to the horizontal modifications $\tilde p, \tilde l, \tilde m$ of the horizontal action $\crta\ot$; two for connecting both $z$ and $w$ with $x$, and two for connecting both $z$ and $w$ with $y$). 
\end{enumerate}
\end{defn}

\medskip

Finally, given horizontally monoidal double categories $\Mm$ on $\Dd$ a bistrong horizontal action of a $\Mm$ on $\Dd$ is a horizontal action 
$\crta\ot$ endowed with a left and a right horizontal strengths $t$ and $s$ and an invertible horizontal modification $q$ satisfying the same two axioms as in \cite[Figure 6, p.9]{HF2} (compatibility of $q$ with $z, z'$ and with $w, w'$).

\medskip

Observe that the 2-cell components of a left vertical strength $t$ at 1h-cells have the form as on the left below, whereas 
the 2-cell components of the induced horizontal transformation $\hat t$ at 1v-cells have the form as on the right below 
$$
\scalebox{0.82}{
\bfig
\putmorphism(-110,170)(1,0)[A(MB)`A'(M'B')`f(mg) ]{900}1a
\putmorphism(-160,200)(0,-1)[``t^M_{A,B}]{410}1l %((AB)C)D`(A(BC))D
\putmorphism(830,200)(0,-1)[``t^{M'}_{A',B'}]{410}1r
\putmorphism(-110,-200)(1,0)[M(AB)`M'(A'B')`m(fg)]{900}1a
\put(190,20){\fbox{$t^m_{f,g}$}}
\efig}
\qquad\qquad
\scalebox{0.82}{
\bfig
\putmorphism(-110,170)(1,0)[A(MB)`M(AB)`\hat t^M_{A,B} ]{900}1a
\putmorphism(-160,200)(0,-1)[``u(Uv)]{410}1l %((AB)C)D`(A(BC))D
\putmorphism(830,200)(0,-1)[``U(uv)]{410}1r
\putmorphism(-110,-200)(1,0)[A'(M'B')`M'(A'B')`\hat t^{M'}_{A',B'}]{900}1a
\put(190,20){\fbox{$\hat t^{U,u,v}$}}
\efig}
$$
whereby we write for short $MA=M\crta\ot A$ for $M\in\Mm, A\in\Dd$, and $f,g$ are 1h-cells in $\Dd$, $m$ is a 1h-cell in $\Mm$ and $u,v$ 
are 1v-cells in $\Dd$ and $U$ is a 1v-cell in $\Mm$.

\bigskip

We next define a {\em commutative} action of double categories and we will show in the next subsection 
that a commutative bistrong action is equivalently %to the action being 
a lax monoidal action. %, in the sense of \deref{v mon act}. 

\begin{defn} 
Let $\Mm$ and $\Dd$ be vertically monoidal double categories whereby $\Mm$ is braided by $\Phi$. 
A vertically bistrong (vertical) action $(\crta\ot,\tilde\alpha,\tilde\lambda; 
t,s,q)$ is called {\em commutative} if there is the identity vertical modification $c$ with 2-cell components
$$c^{MA,MB}: \fourfrac{t^N_{MA,B}}{Ns^M_{A,B}}{\tilde\alpha^{-1}_{N,M,AB}}{\Phi_{N,M}(AB)}
\Rightarrow \threefrac{s^M_{A,NB}}{Mt^N_{A,B}}{\tilde\alpha^{-1}_{M,N,AB}}.$$
\end{defn}

\begin{defn} \delabel{h comm}
Let $\Mm$ and $\Dd$ be horizontally monoidal double categories whereby $\Mm$ is braided by $\Phi$. 
A horizontally bistrong (horizontal) action $(\crta\ot,\tilde\alpha,\tilde\lambda; t,s,q)$ is called {\em commutative} if there is 
an invertible horizontal modification $c$ with 2-cell components
%$$c^{MA,NB}: [s^M_{A,NB} \vert Mt^N_{A,B} \vert \tilde\alpha^{-1}_{M,N,AB}] \Rightarrow
%[t^N_{MA,B} \vert Ns^M_{A,B} \vert \tilde\alpha^{-1}_{N,M,AB} \vert \Phi_{N,M}(AB)]$$
$$ 
\scalebox{0.82}{
\bfig
\putmorphism(-390,450)(1,0)[(MA)(NB)` M(A(NB))`s^M_{A,NB}]{1220}1a 
\putmorphism(-390,450)(0,-1)[\phantom{Y_2}` `=]{450}1r 
\putmorphism(890,450)(1,0)[\phantom{A''\ot B'}`M(N(AB)) ` Mt^N_{A,B}]{1320}1a%%%%%
 \putmorphism(2270,450)(1,0)[\phantom{A''\ot B'}` (MN)(AB) `\tilde\alpha^{-1}_{M,N,AB}]{1200}1a %T(T_A(T_BT_C))

 \putmorphism(-390,0)(1,0)[(MA)(NB)` N((MA)B)`t^N_{MA,B}]{850}1b %%%
 \putmorphism(540,0)(1,0)[\phantom{A''\ot B'}` N(M(AB)) `Ns^M_{A,B}]{900}1b 
 \putmorphism(1540,0)(1,0)[\phantom{A''\ot B'}` (NM)(AB) `\tilde\alpha^{-1}_{N,M,AB}]{900}1b 
 \putmorphism(2520,0)(1,0)[\phantom{A''\ot B'}` (MN)(AB) `\Phi_{N,M}(AB)]{1000}1b 
\putmorphism(3450,450)(0,-1)[\phantom{Y_2}``=]{450}1r 

\put(1380,230){\fbox{$c^{MA,NB}$}}
\efig}
$$
satisfying seven axioms analogous to those from \cite[Appendix E, Figure 8]{HF2}.
\end{defn}

\begin{rem} \rmlabel{diff}
We comment the axioms that $c$ in the above definition needs to fulfill. They %The commutativity axioms 
live in the underlying bicategory of $\Dd$ and are comprised of the 2-cell components of: the strength structure $y,z,w$ and 
$y',z',w'$ of the left and the right strength, the bistrength $q$, commutativity $c$, modifications $b_1, b_2$ from the braiding $\Phi$ 
($b_1$ appears in the axiom (6), and $b_2$ in (7)), and modifications $\tilde p, \tilde l, \tilde m$ 
of the action $(\crta\ot,\tilde\alpha,\tilde\lambda)$ - in the form as in \rmref{monad action}. 
Given the size and complexity of the diagrams expressing the axioms we restrain ourselves to only stress the differences between component diagrams in the cases of a bistrong bicategorical monad $T$ from \cite{HF2} and our horizontally bistrong action of double categories from \deref{hor-act}. First of all, 
as we highlighted in \rmref{gen action}, the action of the multiplication $\mu_A: T^2(A)\to T(A)$ and the unit $\eta_A:A\to T(A)$ of the monad $T$ become 
$$\tilde\alpha^{-1}_{M,N,A}: M\crta\ot(N\crta\ot A)\to(M\ot N)\crta\ot A \qquad \text{and} \qquad \tilde\lambda^{-1}_A:A\to J\crta\ot A$$
where $J\in\Mm$ is the monoidal unit; strength $t_{A,B}: A\ot T(B)\to T(A\ot B)$ of a double monad $T$ is substituted by 
strength $t^M_{A,B}:A\ot (M\crta\ot B)\to M\crta\ot(A\ot B)$ of the action $\crta\ot:\Mm\times\Dd\to\Dd$ and similarly happens 
with $y,z,w,q,\,y',z',w'$; 
% the rest of the structures are changed accordingly and the 2-cell components of the 
commutativity $c$ takes the form as in \deref{h comm}, and the strength component $w$ in our setting (from \deref{hor-strength}) has the form 
$$ 
\scalebox{0.82}{
\bfig
\putmorphism(-390,450)(0,-1)[\phantom{Y_2}` `=]{450}1r 

\putmorphism(-390,450)(1,0)[A(M(NB))` M(A(NB))`t^M_{A,NB}]{870}1a 
\putmorphism(540,450)(1,0)[\phantom{A''\ot B'}`M(N(AB)) ` Mt^N_{A,B}]{900}1a%%%%%
 \putmorphism(1500,450)(1,0)[\phantom{A''\ot B'}` (MN)(AB) `\tilde\alpha^{-1}_{M,N,AB}]{850}1a %T(T_A(T_BT_C))

 \putmorphism(-390,0)(1,0)[A(M(NB))` A((MN)B)`A\tilde\alpha^{-1}_{M,N,B}]{1350}1b 
 \putmorphism(1040,0)(1,0)[\phantom{A''\ot B'}` (MN)(AB). `t^{MN}_{A,B}]{1300}1b 

\putmorphism(2300,450)(0,-1)[\phantom{Y_2}``=]{450}1r 

\put(860,230){\fbox{$w^{M,N}_{A,B}$}}
\efig}
$$
\end{rem}

\subsection{Equivalence of lax monoidal and commutative action} \sslabel{equiv}

In the above definitions of this section and in \deref{h mon act} we used the inverses $\beta, \nu$ of the action transformations 
$\tilde\alpha: (-\ot-)\crta\ot-\Rightarrow-\crta\ot(-\crta\ot-)$ 
and $\tilde\lambda: J\crta\ot\Id_\Dd\Rightarrow\Id_\Dd$ to treat the actions of a monoidal double category $\Mm$ on $\Dd$. In the 
result that we are proving in this subsection it is this direction of the transformations that will be sufficient for the proof. To stress this fact 
we will say for an action given by (non-invertible) $\beta: -\times(-\times-)\Rightarrow (-\ot-)\times-$ and $\nu:-\Rightarrow I\times-$ that it is a 
{\em lax action}, as announced in \rmref{monad action}. (A lax action resembles skew-monoidality in a category. The three modifications $\tilde p, \tilde m, \tilde l$ correspond to the following three axioms in a skew-monoidal category: the pentagon, the one relating $\alpha_{A,I,B}, \lambda_B, \rho_A$, and the one relating 
$\alpha_{I,A,B}$ and $\lambda_A, \lambda_{AB}$. The remaining two axioms in a skew-monoidal category are: one relating $\alpha_{A,B,I}$ and $\rho_B,\rho_{AB}$, and another one relating $\lambda_I$ and $\rho_I$, and in them the $\rho$'s do 
not make sense in a left action.) %Namely, in the former axiom $A,B$ live in the monoidal acting category and $I$ lives in the category acted upon, so the action should be from the right. For the latter axiom $\rho_I$ does not make sense in a left action.)

\smallskip

To shed the doubts if a {\em lax} action is sufficient to formulate a right strength $s$, we type the diagrams for the structure 2-cells of $w', z'$ of $s$ for reader's convenience: 
$$ 
\scalebox{0.82}{
\bfig
\putmorphism(-390,450)(0,-1)[\phantom{Y_2}` `=]{450}1r 

\putmorphism(-390,450)(1,0)[(M(NA))B` M((NA)B)`s^M_{NA,B}]{870}1a 
\putmorphism(540,450)(1,0)[\phantom{A''\ot B'}`M(N(AB)) ` Ms^N_{A,B}]{900}1a%%%%%
 \putmorphism(1500,450)(1,0)[\phantom{A''\ot B'}` (MN)(AB) `\tilde\alpha^{-1}_{M,N,AB}]{850}1a %T(T_A(T_BT_C))

 \putmorphism(-390,0)(1,0)[(M(NA))B` ((MN)A)B`\tilde\alpha^{-1}_{M,N,A}B]{1350}1b 
 \putmorphism(1040,0)(1,0)[\phantom{A''\ot B'}` (MN)(AB) `s^{MN}_{A,B}]{1300}1b 

\putmorphism(2300,450)(0,-1)[\phantom{Y_2}``=]{450}1r 

\put(860,230){\fbox{${w'}^{M,N}_{A,B}$}}
\efig}
$$
$$
\scalebox{0.82}{
\bfig
\putmorphism(50,200)(0,-1)[AB`AB `]{380}1r
\putmorphism(130,210)(1,0)[``\nu_A\w\ot B]{410}1a
\putmorphism(680,210)(1,0)[(JA)B`J(AB)`s^J_{A,B}]{680}1a
\putmorphism(1370,200)(0,-1)[\phantom{Y_2}` `=]{380}1r
\putmorphism(80,-200)(1,0)[\phantom{Y}` J(AB). `\nu_{AB}]{1260}1a
\put(580,20){\fbox{$z'_{A,B}$}}
\efig}
$$

\smallskip

We will first prove the horizontal case. 

\begin{prop} \prlabel{lax monoidal-comm h}
Let $\Mm$ and $\Dd$ be horizontally monoidal double categories with a horizontal lax action  $(F, \beta,\nu)$ %$F:\Mm\times\Dd\to\Dd$ 
and assume that $\Mm$ is braided with double braiding $\Phi$. 
% with vertical structure transformations $\tilde\alpha, \tilde\lambda$, should be a {\em
The following are equivalent: 
\begin{itemize}
\item there is a structure $(F, F^2,F^0)$ of a lax monoidal horizontal lax action, 
\item %} (this means we have 1v-cell components 
%$$(M\ltr C)\ot(M'\ltr C')\longrightarrow (M\ot M')\ltr (C\ot C')$$ 
%where $M\ltr C:=F(M,C)$ - \textcolor{rojo}{this should be equivalent to having 
there is a structure $(F, s,t, q,c)$ of a bistrong and commutative horizontal lax action.
\end{itemize}
\end{prop}

\begin{proof}
Assume  $(F, \beta,\nu)$ is lax monoidal. Define 
$$t^M_{A,B}:=\big(A(MB)\stackrel{\nu_A(MB)}{\to} (JA)(MB) \stackrel{F^2_{JA,MB}}{\to} (JM)(AB)\stackrel{\lambda^\Mm_M(AB)}{\to}M(AB)\big)$$
$$s^M_{A,B}:=\big((MA)B\stackrel{(MA)\nu_B}{\to} (MA)(JB) \stackrel{F^2_{MA,JB}}{\to} (MJ)(AB)\stackrel{\rho^\Mm_M(AB)}{\to}M(AB)\big)$$
and their 2-cell components at 1v-cells $m:M\to \tilde M, a: A\to \tilde A, b:B\to\tilde B$ by the horizontal juxtaposition of the the corresponding 
2-cell components of $\nu\ot 1, F^2$ and $\lambda^\Mm\ot 1$ for $t$, and similarly for $s$. Their 2-cell components at 1h-cells are defined 
similarly by composing the same type of components of $\nu\ot 1, F^2, \lambda^\Mm\ot 1$. Then $s,t$ inherit the property of being pseudonatural 
transformations from their constituting pseudonatural transformations. 

The strength structures $x,y,z,w, \, x',y',z',w'$, the bistrength $q$ and commutativity $c$ are given analogously as in 
\cite[Appendix E, Figure 8]{HF2}, we only list the horizontal modifications of whose constituting 2-cell components the 
latter structures are defined (suppressing the naturality isomorphisms necessarily appearing to make the compositions possible): 
\begin{table}[H]
\begin{center}
\begin{tabular}{ c c } %  m{7em} | m{3.4cm}
%commutative bistrong structure 
define & $\vert\vert$ known data \\ \hline
$x$  & $\nu^0, \gamma$ \\ \hline %[1ex]   
$y$  & $\nu^2, \omega$ \\ \hline %[1ex]   
$z$ & $\nu^2$ \\ \hline %[1ex]   
$w$ & $\tilde m, \beta^2$ \\ \hline %[1ex]   % p of [HF]
$q$ & $\omega$ \\ \hline %[1ex]    
$c$ & $\quad\tilde l, \tilde m, \beta^2, (\beta^2)^{-1}$ \\ \hline %[1ex]   
\end{tabular}
\caption{Defining commutative bistrong action out of lax monoidal action}
\label{table:1}
\end{center}
\end{table}
\noindent Here the horizontal modifications $\tilde l, \tilde m$ are meant in the form as in \rmref{monad action}. 
We mark that bistrength $q$ is denoted by $b$ in \cite{HF2}, and that the invertibility of the modification $\beta^2$ is only 
used in the definition of $c$, as displayed in Table \ref{table:1}. 
The proof that the above define a structure of a commutative bistrong horizontal action is 
completely analogous to the proof in the bicategorical case (\cite[Theorem 7.5]{HF2}), given that in the structures of bistrength and commutativity 
vertical cells do not appear (neither in their laws). The differences only lie in the outlook of the component 2-cells 
of the given structures, as explained in \rmref{diff}, but the steps in the diagrammatic computation use precisely the same laws 
as in the bicategorical proof (including the modification laws for $b_1, b_2$ of the braiding $\Phi$). 

The converse is proved applying the same strategy: %For the other way around, assume that  is 
given a bistrong and commutative horizontal action $F$, to prove that $(F, \beta,\nu)$ 
is a lax monoidal horizontal action, we define 
$$F^2_{M,A,N,B}:=\big( (MA)(NB)\stackrel{s^M_{A,NB}}{\to} M(A(NB)) \stackrel{Mt^N_{A,B}}{\to} M(N(AB))\stackrel{\beta_{M,N,AB}}{\to}(MN)(AB)\big)$$ 
and its 2-cell components by the composition of the 2-cell components of $s,1\ot t$ and $\beta$ (then $F^2$ satisfies the cocycle modification 
property since so does $\beta$), we set 
$\iota:=\nu_I: I\to JI$ and we give the invertible horizontal modifications $\beta^2,\beta^0, \nu^2, \nu^0$ 
(for monoidality of $\beta$) and $\omega, \delta, \gamma$ (for horizontal lax monoidality of $F$) analogously as in 
\cite[Appendix E, Figure 8]{HF2}, which we depict in the Table below, being  $\tilde l, \tilde m, \tilde p$ as in \rmref{monad action}. 
\begin{table}[H]
\begin{center}
\begin{tabular}{ c c } %  m{7em} | m{3.4cm}
%commutative bistrong structure 
define & $\vert\vert$ known data \\ \hline
$\beta^2$  & $w', w, c, \tilde p$ \\ \hline %[1ex]   m
$\beta^0$  & $\tilde m$ \\ \hline %[1ex]   % p of [HF]
$\nu^2$ & $z, z', \tilde l$ \\ \hline %[1ex] njihovo n?  
$\nu^0$ & $\Id_{\nu_I}$ \\ \hline %[1ex]   
$\omega$ & $y',w',q, y, \tilde p, w$ \\ \hline %[1ex]    
$\delta$ & $\quad z, x', \tilde m$ \\ \hline %[1ex]   
$\gamma$ & $\quad z',x, \tilde l$ \\ \hline %[1ex]   
\end{tabular}
\caption{Defining lax monoidal action from commutative bistrong action}
\label{table:2}
\end{center}
\end{table}
%\rmref{hor modif of hor lax mon}
\qed\end{proof}

\bigskip

The vertical version of the above result is proved in a similar way, with the difference that there are less rules to be verified. 
To begin with, the vertical modifications $\tilde p, \tilde l, \tilde m$ of a vertical action $(F, \beta,\nu)$ are identities. 
Furthermore, the strength modifications $x,y,z,w$, the bistrength $q$ and commutativity $c$ are identities, and so are the vertical analogues of 
$\omega, \delta, \gamma$. A priori, the only non-trivial vertical modifications are the cocycle modifications 
$\beta^2, \beta^0, \nu^2, \nu^0$. However, a vertical bistrong commutative action by construction determines trivial cocycle modifications 
defining a lax monoidal action (see Table \ref{table:2}). This will be reflected in the claim of \prref{F is lax monoidal}. 
So, one only needs to verify that the vertical strict transformations $s,t,F^2$ are well-defined (this is proved similarly as above), and 
the checking of horizontal modification axioms in the above proof is substituted by checking that the relevant identity vertical modification 
is given by the corresponding identity of the composite 1v-cells.

\begin{prop} \prlabel{F is lax monoidal}
Let $\Mm$ and $\Dd$ be vertically monoidal double categories with a vertical lax action $(F, \beta,\nu)$ %$F:\Mm\times\Dd\to\Dd$ 
and assume that $\Mm$ is braided with double braiding $\Phi$. 
% with vertical structure transformations $\tilde\alpha, \tilde\lambda$, should be a {\em
The following are equivalent: 
\begin{itemize}
\item there is a structure $(F, F^2, F^0)$ of a lax monoidal vertical lax action with trivial cocycle modifications $\beta^2, \beta^0, \nu^2, \nu^0$, 
\item %} (this means we have 1v-cell components 
%$$(M\ltr C)\ot(M'\ltr C')\longrightarrow (M\ot M')\ltr (C\ot C')$$ 
%where $M\ltr C:=F(M,C)$ - \textcolor{rojo}{this should be equivalent to having 
there is a structure $(F, s,t, q,c)$ of a bistrong and commutative vertical lax action.
\end{itemize}
\end{prop}

\medskip

For the record we only list the 2-cell components of the cocycle modifications:
%\begin{proof}
%Assume $(F, \beta,\nu)$ is a lax monoidal vertical action. For $M\in\Mm, A\in\Dd$ the action of $M$ on $A$ we will denote simply by $MA$. 
%The monoidal unit in $\Mm$ we will denote by $J$. The monoidal structure %of $F$ 
%transformation $F^0$ we will denote by $F^0=\iota:I\to JI$. %$F^2$ and 
%That the transformations $\beta: -\crta\ot(-\crta\ot-)\Rightarrow(-\ot-)\crta\ot-:\Mm\times\Mm\times\Dd\to\Dd$ and $\nu:\Id\Rightarrow J-$
%are pseudomonoidal (according to \deref{v mon act} and \rmref{modif+coc}) it means there are cocycle modifications with 2-cell components 
$$\beta^2: \threefrac{\beta_{M,N,A}\beta_{P,Q,B}}{F^2_{MN,A,PQ,B}}{F^2_{M,N,P,Q}(AB)}\Rightarrow\threefrac{F^2_{M,NA,P,QB}}{(MP)F^2_{N,A,Q,B}}
{\beta_{MP,NQ,AB}}: \quad (M(NA))(P(QB))\to((MP)(NQ))(AB)$$
$$\beta^0:\iota\Rightarrow\fourfrac{\iota}{J\iota}{\beta}{\lambda^\Mm_{J}}: \quad I\to JI$$
$$\nu^2:\nu_{AB}\Rightarrow\threefrac{\nu_A\nu_B}{F^2_{J,A,J,B}}{\lambda^\Mm(AB)}:\quad AB\to J(AB)$$
$$\nu^0:\iota\Rightarrow\nu_I:\quad I\to JI$$
with %for pairs of triples of objects 
$(M,N,A), (P,Q,B)\in\Mm\times\Mm\times\Dd$. %and where $F^2$ is a 2-dimensional 2-cocycle. This means that both $(b_2,b_0)$ and $(n_2,n_0)$ satisfy five axioms each: associativity, two unitalities, naturality with respect to 1h-cells and coherence with interchangers for $ -\crta\ot(-\crta\ot-)$ and $(-\ot-)\crta\ot-$. 
%\qed\end{proof}

% \cite[Proposition 3.5]{LS}: ''Skew monoidales, skew warpings and quantum categories'' i [Booker, Street, Prop.8.3] ne rade u opstijem 
% slucaju akcija monoidalne duple kategorije $\Mm$

\medskip

The results of this subsection concern {\em lax} actions. We remark that invertibility of vertical strict transformations, in particular 
of the action transformations $\tilde\alpha$ and $\tilde\lambda$ and the cocycle $F^2$ (from lax monoidality of pseudodouble functors $F$), 
is required only for the sake of lifting vertical to horizontal structures, as in \prref{ver-hor-dbl}, due to \prref{lifting 1v to equiv}. 
%results that we need: a) a proper action rather than merely a lax action, 
%and b) {\em strongly} , {\em i.e.} to be an invertible vert. strict transformation.}

\section{Para construction on double categories} \selabel{Para}

In this final section we introduce the Para and coPara double category as a generalization of the 
Para construction for 1-categories. 
%The latter is a categorical framework with applications in computer science. %Our interest in it comes from the fact that it 
%Given that the Para construction presents a generalization of the Kleisli category, our motivation is to extend our double categorical results on the Kleisli double category from \cite{Fem3} to the setting of coPara for double categories, facilitating in this way a wider range of applications of this construction 
%extensively used in computer science. 
%
We will focus on the horizontal stuctures and will discuss vertical versions of the results at the end. 

\smallskip

Let $\Mm$ be a horizontally monoidal double category and $\Dd$ a double category, and suppose that $\Mm$ acts on $\Dd$ by a horizontal action 
$\crta\ot:\Mm\times\Dd\to\Dd$. We introduce a pseudodouble category $\Para_\Mm(\Dd)$ as follows: \\
\ul{0-cells:} objects of $\Dd$; \\
\ul{1v-cells:} 1v-cells of $\Dd$; \\
\ul{1h-cells:} pairs $(M,f):D\to E$ where $M$ is an object of $\Mm$ and $f: M\crta\ot D\to E$ is a 1h-cell in $\Dd$; \\
\ul{2-cells:} as on the left below %pairs $(u,a)$ where $u: M\to M'$ is a 1v-cell in $\Mm$ and $a$ is a 
are given by 2-cells in $\Dd$ as on the right below
$$
\bfig
\putmorphism(-160,120)(1,0)[ D`E`(M,f)]{550}1a
\putmorphism(-160,130)(0,-1)[`` d]{380}1l
\putmorphism(370,130)(0,-1)[``e]{380}1r
\putmorphism(-160,-250)(1,0)[ D'`E'`(M',f')]{550}1a
%\put(60,-30){\fbox{$a$}}
\efig
\qquad\qquad
\bfig
\putmorphism(-160,120)(1,0)[M\crta\ot D`E`f]{550}1a
\putmorphism(-160,130)(0,-1)[``u\crta\ot d]{380}1l
\putmorphism(370,130)(0,-1)[``e]{380}1r
\putmorphism(-160,-250)(1,0)[M'\crta\ot D'`E'`f']{550}1a
\put(60,-30){\fbox{$a$}}
\efig
$$
where $u: M\to M'$ is a 1v-cell in $\Mm$. % and $d:D\to D'$ and $e:E\to E'$ are 1v-cells in $\Dd$. 

Identity 1h-cell on $A\in\Dd$ is given by $(I,\tilde\lambda_A): A\to A$, namely $\tilde\lambda_A:I\crta\ot A\to A$. 
Composition of 1h-cells $(M,f):D\to E$ %(given by $f: M\crta\ot D\to E$) 
and $(N,g):E\to F$ %(given by $g: N\crta\ot E\to F$) 
is given by the following composition of 1h-cells in $\Dd$ 
$$(N\ot M)\crta\ot D\stackrel{\tilde\alpha}{\longrightarrow} N\crta\ot(M\crta\ot D) \stackrel{N\crta\ot f}{\longrightarrow} N\crta\ot E
\stackrel{g}{\longrightarrow} F.$$
Accordingly, horizontal composition of 2-cells $a:(f,f';d,e)$ and $b:(g,g';e,f)$ is given by the horizontal composition of 2-cells in $\Dd$ 
%$$[\tilde\alpha^{v,u,d}\,\, \vert\,\, v\crta\ot a \,\,\vert\,\, b].$$
$$\scalebox{0.82}{
\bfig
 \putmorphism(-950,210)(1,0)[(NM)D` N(MD)`\tilde\alpha_{N,M,D}]{960}1a %((AB)C)T(D') 
\putmorphism(-860,200)(0,-1)[\phantom{Y_2}` `(uv)d]{380}1l
\putmorphism(-110,200)(0,-1)[\phantom{Y_2}` `u(vd)]{380}1r
\putmorphism(120,210)(1,0)[\phantom{Y}` NE  `N\crta\ot f]{580}1a
\putmorphism(740,210)(1,0)[\phantom{Y}` F  `g]{460}1a
\putmorphism(680,200)(0,-1)[\phantom{Y_2}` `ue]{380}1r

\putmorphism(1200,200)(0,-1)[\phantom{Y_2}` `u(vd)]{380}1r
\putmorphism(-970,-200)(1,0)[(N'M')D' `N'(M'D'))`\tilde\alpha_{N',M',D'}]{980}1a % (\tilde A\tilde B)\tilde C
\putmorphism(220,-200)(1,0)[``N'\crta\ot f']{380}1a
\putmorphism(680,-200)(1,0)[N'E'` F'. `g']{530}1a
\put(-600,30){\fbox{$\alpha^{u,v,d}$}}
\put(240,30){\fbox{$ u\crta\ot a$}}
\put(900,30){\fbox{$b$}}
\efig}
$$
Vertical composition of 2-cells is given by the one in $\Dd$. 

\begin{rem}
Observe that it is for the composition of 1h-cells that we need a {\em horizontal} action, and hence {\em horizontally} monoidal 
double category $\Mm$. 
\end{rem}

Given that we are generalizing our results from \cite{Fem3}, we will indeed work with the dual coPara construction in which 
1h-cells are of the form $A\to M\crta\ot B$. For the double category $\coPara_\Mm(\Dd)$ to be well-defined, concretely the 
composition of 1h-cells $f:A\to M\crta\ot B$ and $g:B\to N\crta\ot C$ (and similarly for 2-cells), we will need the inverse direction of the associativity of the action $\crta\ot:\Mm\times\Dd\to\Dd$
$$A\stackrel{f}{\longrightarrow} M\crta\ot B \stackrel{M\crta\ot g}{\longrightarrow} M\crta\ot(N\crta\ot D)
\stackrel{\beta}{\longrightarrow} (M\ot N)\crta\ot D.$$
It is straightforwardly seen that in order to get unitality of the composition in $\coPara_\Mm(\Dd)$ the right direction for the unitality of the action is $\tilde\lambda^{-1}_A=\nu_A:A\to I\crta\ot A$. It clearly also determines the identity 1h-cell on $A$ in $\coPara_\Mm(\Dd)$. Accordingly, 
we will be indeed interested in lax actions $\crta\ot:\Mm\times\Dd\to\Dd$ %and we will write the colax associativity and unitality 
%(vertical) transformations as $\beta$ and $\nu$, respectively, throughout (in place of $\tilde\alpha$ and $\tilde\lambda$ from 
%\seref{bistrong,comm}).} 
as in the previous section.

\begin{defn}  \delabel{Para embed} 
Let $\Dd$ be a double category and let $(\crta\ot, \Mm, \Dd,\beta, \nu)$ be a horizontal lax action of a horizontally monoidal double 
category $\Mm$ on $\Dd$. The canonical embedding $P:\Dd\to \coPara_\Mm(\Dd)$ is the pseudodouble functor that is the identity on objects and vertical 1-cells, it sends a horizontal 1-cell $f: A\to B$ into a 1h-cell $P(f): A\to B$ in $\coPara_\Mm(\Dd)$ determined by 
$\nu_B\comp f: A\to J\crta\ot B$ in $\Dd$, and correspondingly a 2-cell $\phi$ to the horizontal composition of 2-cells %$[\omega\vert\Id_{\eta_B}]$. 
$$
\scalebox{0.86}{
\bfig
 \putmorphism(-50,150)(1,0)[A`\phantom{F(A)} `f]{520}1a
\putmorphism(-50,-200)(1,0)[\tilde A`\tilde B `g]{460}1b
\putmorphism(-50,150)(0,-1)[``u]{350}1l
\put(90,-20){\fbox{$\phi$}}
\putmorphism(410,150)(0,-1)[B` `v]{350}1l
 \putmorphism(350,150)(1,0)[\phantom{F(A)}`JB `\nu(B)]{600}1a
 \putmorphism(450,-200)(1,0)[`J\tilde B.`\nu(\tilde B)]{500}1b
\putmorphism(930,150)(0,-1)[\phantom{Y_2}``Jv]{350}1r
\put(590,-30){\fbox{$\nu^v$}}
\efig}
$$
\end{defn}

For $\Mm$ the trivial double category $P:\Dd\to \coPara_\Mm(\Dd)$ recovers the  canonical embedding $K:\Dd\to\Kk l(T)$ 
into the Kleisli double category from \cite[Definition 9.2]{GGV}.

\subsection{Actions on coPara} \sslabel{act Para}

From now on by a horizontal lax action $(\crta\ot, \Mm, \Dd,\beta, \nu)$ we will mean a horizontal lax action 
$\crta\ot:\Mm\times\Dd\to\Dd$ of $\Mm$ on $\Dd$ whereby both double categories are horizontally monoidal. 

Analogously as in \cite[Proposition 10.28]{Fem3} for the Kleisli double category, we have:
%\bigskip

\begin{prop} \prlabel{horiz str-icon}
A horizontal left strength $t$ on a horizontal lax action $(\crta\ot, \Mm, \Dd,\beta, \nu)$ 
induces a horizontal icon $\theta:-\crta\ot P(-)\Rightarrow P(-\ot-): \Dd\times\Dd\to\coPara_\Mm(\Dd)$ whose 
invertible 2-cell components $\theta_{f,g}$ for 1h-cells $f:A\to A', g:B\to B'$ are given by 
$$\theta_{f,g}=
\scalebox{0.82}{
\bfig
 \putmorphism(-650,210)(1,0)[AB` A'B'`f\ot g]{640}1a %((AB)C)T(D') 
\putmorphism(-660,200)(0,-1)[\phantom{Y_2}` `=]{380}1l
\putmorphism(-50,200)(0,-1)[\phantom{Y_2}` `]{380}1r
\putmorphism(-50,200)(0,-1)[\phantom{Y_2}` `=]{380}0r
\putmorphism(100,210)(1,0)[``A'\w\ot\nu_{B'}]{410}1a
\putmorphism(680,210)(1,0)[A'(JB')`J(A'B')`t^J_{A',B'}]{680}1a

\putmorphism(1600,200)(0,-1)[\phantom{Y_2}``]{380}0r
\putmorphism(1370,200)(0,-1)[\phantom{Y_2}` `=]{380}1r
\putmorphism(-670,-200)(1,0)[AB `A'B'`f\ot g]{620}1a %%%
\putmorphism(20,-200)(1,0)[\phantom{Y}` J(A'B')  `\nu_{A'B'}]{1360}1a

\put(-500,30){\fbox{$\Id_{f\ot g}$}}
\put(580,30){\fbox{$z_{A',B'}$}}
\efig}
$$
and the square-formed 2-cell components $\theta^{u,v}=\Id^{u\ot v}$ for 1v-cells $u,v$. 
\end{prop}

The proof of the next theorem is analogous to that of \cite[Theorem 10.29]{Fem3}. We will only outline the involved structures. 

\begin{thm} \thlabel{hor-str-act}
A left horizontal strength $t$ on a horizontal lax action $(\crta\ot, \Mm, \Dd,\beta, \nu)$ 
%on a horizontally monoidal double category $(\Dd, \ot, \alpha, \lambda, \rho)$ 
induces an action $\rtr:\Dd\times\coPara_\Mm(\Dd)\to\coPara_\Mm(\Dd)$.  
\end{thm}
 
A pseudodouble functor $\crta\ot:\Dd\times\coPara_\Mm(\Dd)\to\coPara_\Mm(\Dd)$ is defined on objects and 1v-cells by the action of the monoidal product of $\Dd$, for 1h-cells $f:A\to A'\in\Dd$ and $g:B\to M\crta\ot B'\in\coPara_\Mm(\Dd)$ we define 
$$f\crta\ot g:=(A\ot B\stackrel{f\ot g}{\to}A'\ot (M\crta\ot B') \stackrel{t_{A',B'}^M}{\to}M\crta\ot(A'\ot B'))$$
and for 2-cells $\sigma\in\Dd, \delta\in\coPara_\Mm(\Dd)$ we set 
$$
\scalebox{0.86}{
\bfig
 \putmorphism(-50,250)(1,0)[AB`A'(MB') `f\ot g]{720}1a
\putmorphism(-50,-200)(1,0)[\tilde A\tilde B`C(ND) `\tilde f\ot\tilde g]{720}1b
\putmorphism(-50,250)(0,-1)[``uu']{450}1l
\put(30,30){\fbox{$\sigma\ot\delta$}}
\putmorphism(730,250)(0,-1)[` `v(Uv')]{450}1l
 \putmorphism(750,250)(1,0)[\phantom{F(A)}`M(A'B') `t^M_{A',B'}]{600}1a
 \putmorphism(850,-200)(1,0)[`N(CD)`t^N_{C,D}]{500}1b
\putmorphism(1280,250)(0,-1)[\phantom{Y_2}``U(vv')]{450}1r
\put(900,30){\fbox{$t^{U,v,v'}$}}
\efig}
$$
whereby in objects and 1v-cells we omit writing the symbols $\ot$ and $\crta\ot$. 
The compositor 2-cell for $\rtr$ is given by the globular 2-cell 
$$ 
\scalebox{0.82}{
\bfig
 \putmorphism(300,450)(1,0)[A'(MB')` M(A'B') `t^M_{A',B'}]{860}1a
 \putmorphism(1310,450)(1,0)[` M(A''(NB''))  `M(f'g')]{780}1a

\putmorphism(340,450)(0,-1)[\phantom{Y_2}``=]{450}1l
\putmorphism(1150,0)(0,-1)[\phantom{Y_2}``=]{450}1l

 \putmorphism(-470,-450)(1,0)[\phantom{A''\ot B'}` A''(M(NB'')) `(f'f) ((Mg')g)]{1460}1b %%%
 \putmorphism(1110,-450)(1,0)[\phantom{A''\ot B'}` A''((MN)B'') `A''\beta_{M,N,B''}]{1340}1b 
 \putmorphism(2540,-450)(1,0)[\phantom{A''\ot B'}`  `t^{MN}_{A'',B''}]{1100}1b 

\put(980,220){\fbox{$(t^M)^{-1}_{f',g'}$}}
\put(280,-250){\fbox{$lax \ot$}}

\putmorphism(-390,0)(1,0)[AB` A'(MB')`f g]{700}1a 
\putmorphism(-390,0)(0,-1)[\phantom{Y_2}`AB `=]{450}1r 
\putmorphism(340,0)(1,0)[\phantom{A''\ot B'}`A''(M(NB'')) ` f'(Mg')]{850}1a%%%%%
 \putmorphism(1280,0)(1,0)[\phantom{A''\ot B'}`  `t^M_{A'',NB''}]{520}1a %T(T_A(T_BT_C))

 \putmorphism(2160,0)(1,0)[\phantom{A''\ot B'}` M(N(A''B'')) `M(t^N_{A'',B''})]{810}1a %T(T_A(T_BT_C))
 \putmorphism(3070,0)(1,0)[\phantom{A''\ot B'}`(MN)(A''B'')  `\beta_{A'',B''}]{800}1a %T(T_A(T_BT_C))

\putmorphism(2060,450)(0,-1)[\phantom{Y_2}`M(A''(NB'')) ` =]{450}1r 
\putmorphism(3900,0)(0,-1)[\phantom{Y_2}`(MN)(A''B'').`=]{450}1r 
\put(2360,-250){\fbox{$w^{M,N}_{A'',B''}$}}
\efig}
$$
The unitor is defined via $z$.

The action constraints $\tilde\alpha^L, \tilde\lambda$ for $\rtr$ 
are defined on objects and 1v-cells as $\alpha$ and $\lambda$ (of $\Dd$), and on 1h-cells as  
$$\tilde\alpha^L_{f,g,h}=\frac{[\Id_{(fg)h}\,\vert\,y^M_{A',B',C'}]}{[\alpha_{f,g,h}\,\vert\,\Id_{[A't^M_{B',C'}\vert t^M_{A',B'C'}]}]}
\quad\text{ and }\quad
\tilde\lambda_h=\frac{[\Id_{Ih}\vert  x^M_{A'}]}{\lambda_h}$$
for a 1h-cell $(M,h)$ in $\coPara_\Mm(\Dd)$, {\em i.e} $h:C\to MC$ in $\Dd$ 
(observe that all the 2-cells involved are globular, including $\alpha_{f,g,h}$ and $\lambda_h$). 

The horizontal modifications $\tilde l, \tilde m^L, \tilde p^L$ are defined by 
$$%\begin{equation} \eqlabel{hor-p,l,m-teta}
\frac{\theta_{\lambda_A, B}}{P(l_{AB})}=\tilde l_{AB}, \qquad 
\threefrac{[\Id_{P(\alpha_{A,I,B})} \vert\theta_{A,\lambda_B}]}{P(m_{AB})}{(\theta_{\rho_A, B})^{-1}}=\tilde m^L_{AB}, \qquad 
\frac{[\theta_{P(\alpha_{A,B,C}),D} \vert P(\alpha_{A,BC,D})\vert\theta_{A,P(\alpha_{B,C,D})}]}{P(p_{ABCD})}=
\tilde p^L_{ABCD}
$$%\end{equation}
where $\theta$ is from \prref{horiz str-icon} and $l,m,p$ are the horizontal modifications of the horizontally monoidal structure on $\Dd$.  
Here both $\tilde l, \tilde m^L, \tilde p^L$ and $l,m,p$ are meant in the form as in \deref{hor-act}.  In particular, one has 
$\tilde m^L_{AB}=P(m_{AB})$. 

\bigskip
 
Similarly to left actions $\crta\ot:\Mm\times\Dd\to\Dd$ from \deref{hor-act} there are right actions $\crta\ot:\Dd\times\Mm\to\Dd$ with 
horizontal equivalences with components 
$$\tilde\alpha^R_{E,M,N}:(E\crta\ot M)\crta\ot N\to E\crta\ot(M\crta\ot N)\quad\text{and}\quad\tilde\rho_E: E\crta\ot I\to E$$
with $M,N\in\Mm$ and $E\in\Ee$, and horizontal modifications $\tilde p^R, \tilde m^R, \tilde l^R$ with components 
$$%\begin{equation} \eqlabel{effects}
\scalebox{0.8}{
\bfig
 \putmorphism(-150,400)(1,0)[((EM)N)P`(E(MN))P `\tilde\alpha_{E,M,N}P]{1100}1a
 \putmorphism(1050,400)(1,0)[\phantom{A\ot B}`E((MN)P)`\tilde\alpha_{E,MN,P}]{1080}1a
 \putmorphism(2220,400)(1,0)[\phantom{A\ot B}` E(M(NP)) `E\alpha^\Mm_{M,N,P}]{1120}1a

 \putmorphism(-170,50)(1,0)[((EM)N)P`(EM)(NP) `\alpha_{EM,N,P}]{1620}1b
 \putmorphism(1560,50)(1,0)[\phantom{A\ot B}`E(M(NP)) `\tilde\alpha_{E,M,NP}]{1800}1b %(M(NP))E

\putmorphism(-180,400)(0,-1)[\phantom{Y_2}``=]{350}1r
\putmorphism(3290,400)(0,-1)[\phantom{Y_2}``=]{350}1l
\put(1380,230){\fbox{$\tilde p^{M,N,P}_E$}}
\efig}
$$%\end{equation}

$$%\begin{equation} \eqlabel{effects}
\scalebox{0.8}{
\bfig
 \putmorphism(-150,400)(1,0)[(EM)I `EM `\tilde\rho_{EM}]{1430}1a

 \putmorphism(-210,50)(1,0)[(EM)I `E(MI) `\tilde\alpha_{E,M,I}]{670}1b %I(ME)
 \putmorphism(460,50)(1,0)[\phantom{A\ot B}`EM `E\rho^\Mm_M]{780}1b

\putmorphism(-180,400)(0,-1)[\phantom{Y_2}``=]{350}1r
\putmorphism(1280,400)(0,-1)[\phantom{Y_2}``=]{350}1l
\put(540,220){\fbox{$\tilde r^M_E$}}
\efig}
$$

$$%\begin{equation} \eqlabel{effects}
\scalebox{0.8}{
\bfig
 \putmorphism(-210,400)(1,0)[(EI)M`E(IM) `\tilde\alpha_{E,I,M}]{670}1a
 \putmorphism(460,400)(1,0)[\phantom{A\ot B}`EM `E\lambda^\Mm_M]{800}1a
% \putmorphism(1260,400)(1,0)[\phantom{A\ot B}`ME`\rho^\Mm_ME]{800}1a
\putmorphism(-150,50)(1,0)[(EI)M`EM.`\tilde\rho_EM]{1430}1b

\putmorphism(-180,400)(0,-1)[\phantom{Y_2}``=]{350}1r
\putmorphism(1280,400)(0,-1)[\phantom{Y_2}``=]{350}1l
\put(540,220){\fbox{$\tilde m^M_E$}}
\efig}
$$

Now, similarly to \thref{hor-str-act}, a right horizontal strength $s$ on a horizontal lax action $(\crta\ot, \Mm, \Dd,\beta, \nu)$ 
induces a right action $\ltr:\coPara_\Mm(\Dd)\times\Dd\to\coPara_\Mm(\Dd)$ with the action constraints $\tilde\alpha^R, \tilde\rho$ for $\ltr$ 
defined on objects and 1v-cells as $\alpha$ and $\rho$ (of $\Dd$), and on 1h-cells by
$$\tilde\alpha^R_{f,g,h}=\frac{[\Id_{(fg)h}\,\vert\,y'^M_{A',B',C'}]}{[\alpha_{f,g,h}\,\vert\,\Id_{s^M_{A',B'C'}}]}
\quad\text{ and }\quad
\tilde\rho_f=\frac{[\Id_{fI}\vert  {x'}^M_{A'}]}{\rho_f}$$
with $f: A\to MA'$, whereby all the appearing 2-cells are globular, %2-cells $\alpha_{f,g,h}$ and $\rho_h$, 
and horizontal modifications $\tilde r, \tilde m^R, \tilde p^R$ defined by 
\begin{equation} \eqlabel{hor-p,l,m-teta-s}
\frac{\theta_{A,\rho_B}}{P(r_{AB})}=\tilde r_{AB}, \qquad 
%\threefrac{[\Id_{P(\alpha_{A,I,B})} \vert\theta_{A,\lambda_B}]}{P(m_{AB})}{(\theta_{\rho_A, B})^{-1}}
P(m_{AB})
=\tilde m^R_{AB}, \qquad 
\frac{[\theta_{P(\alpha_{A,B,C}),D} \vert P(\alpha_{A,BC,D})\vert\theta_{A,P(\alpha_{B,C,D})}]}{P(p_{ABCD})}=
\tilde p^R_{ABCD}. 
\end{equation}

\begin{thm} \thlabel{Kl-prem}
If a horizontal lax action $(\crta\ot, \Mm, \Dd,\beta, \nu)$ is bistrong, then the pseudodouble category 
$\coPara_\Mm(\Dd)$ is premonoidal. 
\end{thm}

\begin{proof}
Let $(t,s,q)$ be a bistrength with usual notations. By the above theorem we have a left and a right action of $\Dd$ on 
$\coPara_\Mm(\Dd)$ with structures $\tilde\alpha^L, \tilde\lambda, \tilde l^L, \tilde m^L, \tilde p^L$ for $\rtr$ and 
$\tilde\alpha^R, \tilde\rho, \tilde r^R, \tilde m^R, \tilde p^R$ for $\ltr$, where by construction it is 
$\tilde\alpha^L_{A,B,C}=\tilde\alpha^R_{A,B,C}, \,\, \tilde m^L_{A,B}=\tilde m^R_{A,B}, \,\, \tilde p^L_{ABCD}=\tilde p^R_{ABCD}$. 

We define a binoidal structure on $\coPara_\Mm(\Dd)$ by setting 
\begin{equation} \eqlabel{bino}
A\ltimes g:=(A\ot B\stackrel{A\ot g}{\to}A\ot MB' \stackrel{ t^M_{A,B'}}{\to}M(AB'))
\end{equation}
$$f\rtimes B:=(A\ot B\stackrel{f\ot B}{\to}NA'\ot B \stackrel{s^N_{A',B}}{\to}N(A'B))$$
on 1h-cells $(N,f):A\to A'$ and $(M,g):B\to B'$ in $\coPara_\Mm(\Dd)$. On objects and 1v-cells we define it by the action of the monoidal product of $\Dd$ 
and on 2-cells similarly as in \thref{hor-str-act}. That thus defined assignments $A\ltimes-, -\rtimes B: 
\coPara_\Mm(\Dd)\to\coPara_\Mm(\Dd)$ determine pseudodouble functors is proved analogously as for $\rtr$ of \thref{hor-str-act}. 
We now define three $\overline\alpha$'s as in \deref{premon} which should live in $\coPara_\Mm(\Dd)$. We set: 
$$\overline\alpha_{-,B,C}:=\tilde\alpha^R_{-,B,C}: (-\rtimes B)\rtimes C\Rightarrow -\rtimes(B\bowtie C)$$
$$\overline\alpha_{A,B,-}:=\tilde\alpha^L_{A,B,-}: (A\bowtie B)\ltimes -\Rightarrow A\ltimes(B\ltimes -)$$
$$\hspace{-1cm}\overline\alpha_{A,-,C}:%=[\alpha_{A,-,C}\,\vert\, q^*_{A,-,C}]: 
(A\ltimes -)\rtimes C\Rightarrow A\ltimes(-\rtimes C)$$
%and $\overline\alpha_{A,-,C}: (A\ltimes -)\rtimes C\Rightarrow A\ltimes(-\rtimes C)$ 
where the latter we define on objects and 1v-cells as $\alpha$ and by 
$$\overline\alpha_{A,g,C}:=
\frac{[\Id_{(Ag)C}\,\vert\,q^M_{A,B',C}]}{[\alpha_{A,g,C}\,\vert\,\Id_{[As^M_{B',C'}\vert t^M_{A,B'C}]}]}$$ 
%and $\overline\alpha_{A,v,C}:=\alpha_{A,v,C}$ at a 1h-cell $g:B\to MB'$ and a 1v-cell $v$. 
at a 1h-cell $g:B\to MB'$. 
We already know that $\tilde\alpha^R_{-,B,C}$ and $\tilde\alpha^L_{A,B,-}$ are horizontal pseudonatural transformations, and that 
$\tilde\alpha^R_{-,B,C}, \overline\alpha_{A,-,C}, \tilde\alpha^L_{A,B,-}$ on objects $A,B,C$, respectively, coincide with $\alpha_{A,B,C}$. 
Also $\overline\alpha_{A,-,C}$ is a horizontal pseudonatural transformation, because $\alpha_{A,-,C}$ is, and $q^{-}_{A,-,C}$ is a 
horizontal modification and $s,t$ are horizontal transformations. For 
$\overline\lambda: I\ltimes -\Rightarrow\Id$ and $\overline\rho: -\rtimes I\to \Id$ 
we set to be $\tilde\lambda: I\rtr -\Rightarrow\Id$ and $\tilde\rho: -\ltr I\Rightarrow \Id$ from the left and the right action.  

To prove that the common 1h-cell component $\alpha_{A,B,C}$ of $\overline\alpha_{-,B,C}, \overline\alpha_{A,B,-}$ and $\overline\alpha_{A,-,C}$ 
is central, we first introduce some notation. From (horizontal) monoidality of $\Dd$ we have that $-\ot-:\Dd\times\Dd\to\Dd$ is a pseudodouble functor. In particular for 1h-cells $f:A\to A'$ and $g:B\to B'$ let us denote %for the sake of this proof 
by $\ot^{lax}_{(1g)(fB')}: (f\ot B')(A\ot g)\Rightarrow f\ot g$ and $\ot^{colax}_{(fB)(A'g)}: f\ot g \Rightarrow (A'\ot g)(f\ot B)$ 
the corresponding 2-cell components of the pseudodouble functor structure $(-\ot-)^2$ of $\ot$. We then define 
$$
\scalebox{0.86}{
\bfig
\putmorphism(-1000,230)(0,-1)[\phantom{Y_2}` `\alpha_{A,B,C}\ltimes-\vert_{k^P}:=]{380}0l

\putmorphism(100,210)(1,0)[((AB)C)(MD')` M(((AB)C)D')  `t^M_{(AB)C,D'}]{1100}1a
 \putmorphism(1280,210)(1,0)[\phantom{((AB)C)D}` M((A(BC))D')`M(\alpha_{A,B,C}D')]{1150}1a %((AB)C)T(D') 

\putmorphism(-100,210)(0,-1)[\phantom{Y_2}` `=]{380}1l
\putmorphism(2300,210)(0,-1)[\phantom{Y_2}` `=]{380}1l

\putmorphism(-960,-200)(0,-1)[\phantom{Y_2}` `=]{380}1l
\putmorphism(1230,-200)(0,-1)[\phantom{Y_2}` `=]{380}1r

\putmorphism(-1070,-200)(1,0)[((AB)C)D `((AB)C)(MD')`A(BC)\ot k]{1130}1a %%%
\putmorphism(110,-200)(1,0)[\phantom{((AB)C)MD}`(A(BC))(MD')`\alpha_{A,B,C}(MD')]{1170}1a 
\putmorphism(1280,-200)(1,0)[\phantom{((AB)C)(MD)}`M((A(BC))D) ` t^M_{A(BC),D'}]{1100}1a

\putmorphism(-1070,-580)(1,0)[((AB)C)D ``\alpha_{A,B,C}D]{840}1a %%%
\putmorphism(0,-580)(1,0)[(A(BC))D `(A(BC))(MD') ` A(BC)\ot k]{1230}1a %%%

\put(-200,-380){\fbox{$\ot^{lax}_{(1k)(\alpha_{A,B,C}1)}$}}
\put(1000,0){\fbox{$t^M_{\alpha_{A,B,C},1_{D'}}$}}
\efig}
$$%\quad\text{and}\quad 
and 
$$
\scalebox{0.86}{
\bfig
\putmorphism(-1000,230)(0,-1)[\phantom{Y_2}` `-\rtimes\alpha_{B,C,D}\vert_{f^P}:=]{380}0l

\putmorphism(100,210)(1,0)[(NA')((BC)D)` N(A'((BC)D))  `s^N_{A',(BC)D}]{1100}1a
 \putmorphism(1280,210)(1,0)[\phantom{((AB)C)D}` N(A'(B(CD))) ` N(A'\alpha_{B,C,D})]{1150}1a %((AB)C)T(D') 

\putmorphism(-100,210)(0,-1)[\phantom{Y_2}` `=]{380}1l
\putmorphism(2300,210)(0,-1)[\phantom{Y_2}` `=]{380}1l

\putmorphism(-960,-200)(0,-1)[\phantom{Y_2}` `=]{380}1l
\putmorphism(1230,-200)(0,-1)[\phantom{Y_2}` `=]{380}1r

\putmorphism(-1070,-200)(1,0)[A((BC)D) `(NA')((BC)D)`f\ot (BC)D]{1130}1a %%%
\putmorphism(110,-200)(1,0)[\phantom{((AB)C)MD}`(NA')(B(CD))`(NA')\alpha_{B,C,D}]{1170}1a 
\putmorphism(1280,-200)(1,0)[\phantom{((AB)C)(MD)}`N(A'(B(CD))) ` s^N_{A', B(CD)}]{1100}1a

\putmorphism(-1070,-580)(1,0)[A((BC)D) ``A\alpha_{B,C,D}]{840}1a %%%
\putmorphism(0,-580)(1,0)[A(B(CD)) `(NA')(B(CD)) ` f\ot B(CD)]{1230}1a %%%

\put(-200,-380){\fbox{$\ot^{colax}_{(f1)(1\alpha_{B,C,D})}$}}
\put(1000,0){\fbox{$s^N_{1_{A'},\alpha_{B,C,D}}$}}
\efig}
$$
2-cell components of to be horizontal pseudonatural transformations $\alpha_{A,B,C}\ltimes-$ and $-\rtimes\alpha_{B,C,D}$ 
at 1h-cells $k^P=(M,k):D\to D'$ and $f^P=(N,f):A\to A'$ in  $\coPara_\Mm(\Dd)$. At 1v-cells $u:A\to\tilde A$ and $w: D\to \tilde D$ we set 
$\alpha_{A,B,C}\ltimes-\vert_w:=\Id_{\alpha_{A,B,C}}\ot 1^w$ and $-\rtimes\alpha_{B,C,D}\vert_u:=1^u\ot\Id_{\alpha_{B,C,D}}$. 
Since $(-\ot-)^2, t,s, \alpha$ are horizontal pseudonatural transformations, so are $\alpha_{A,B,C}\ltimes-$ and $-\rtimes\alpha_{B,C,D}$, 
and hence the 1h-cells $\alpha_{A,B,C}$ are central. 

\medskip

For the ten horizontal modifications from \deref{premon} we define four $\crta p$'s and two of each of $\crta m, \crta r, \crta l$ as the corresponding restrictions of $\tilde p_{ABCD}^L=\tilde p_{ABCD}^R, \,\, \tilde m_{AB}^L=\tilde m_{AB}^R, 
 \,\, \crta l_{AB}:=\tilde l^A_B,  \,\, \crta r_{AB}:=\tilde r^B_A$, respectively.  Then clearly they all are horizontal modifications 
and obey the necessary axioms. 
\qed\end{proof}

\subsection{Premonoidality versus monoidality}

For a double category $\Dd$, in \cite[Theorem 7.10]{Fem3} we proved that there is a 1-1 correspondence between its vertically premonoidal structures 
(in the sense of \cite[Definition 3.11]{Fem3}) such that the binoidal structure comes from a pseudodouble quasi-functor 
(equivalently, ``purely central'' binoidal structure/double category), 
and its vertically monoidal structures. The definition of a pseudodouble quasi-functor we defer to the Appendix. Recall that 
the difference between vertical and horizontal versions of the notions, of both premonoidal and monoidal structures, is that the structural transformations 
$\alpha, \lambda, \rho$ and modifications $p,m,l,r$ are vertical in the one and horizontal in the other version. 
The analogous 1-1 correspondence holds true with the horizontal versions of the notions. We discuss this below. 

\smallskip

In \cite[Theorem 5.7]{Fem2} we proved a double category equivalence involving lax double functors, whose pseudodouble functor version is 
a double category equivalence 
$$\F''\colon q\x\Ps_{hop}^{st}(\Aa\times\Bb,\Cc) \to \Ps_{hop}(\Aa\times\Bb,\Cc).$$
On the left-hand side the objects are binary pseudodouble quasi-functors, the 1h- and 1v-cells are the horizontal oplax transformations and vertical 
strict transformations of those, respectively, and the 2-cells are their modifications, while on the right-hand side is the double category of 
pseudodouble binary functors and their corresponding transformations and modifications. Based on this result and considering the vertical category part 
of the above double equivalence, in \cite[Theorem 7.5]{Fem3} we proved category equivalencies 
\begin{equation} \eqlabel{ps-eq}
\F_3\colon q_3\x\operatorname{Ps}_{vst}^{st}(\Bb\times\Bb\times\Bb,\Bb) \to \operatorname{Ps}_{vst}(\Bb\times\Bb\times\Bb,\Bb)
\end{equation}
and 
$$\F_4\colon q_4\x\operatorname{Ps}_{vst}^{st}(\Bb\times\Bb\times\Bb\times\Bb,\Bb) \to 
\operatorname{Ps}_{vst}(\Bb\times\Bb\times\Bb\times\Bb,\Bb).$$
On them the following theorem relies, meant with vertical (pre)monoidal structures.

\begin{thm} \thlabel{premon-mon} \cite[Theorem 7.10]{Fem3} \\
There is a purely central premonoidal double category structure $(H, \theta, I)$ on $\Dd$ with trivial 2-cells $(u,U)$ if and only if
there is a monoidal double category structure $(\ot, \Sigma, I)$ on $\Dd$. 
\end{thm}

Here $\theta$ is a vertical strict transformation from the left in \equref{ps-eq}, and it is comprised of three vertical strict transformations 
$\crta\alpha$ of vertical premonoidality, and $\Sigma$ is the vertical strict transformation of the associativity of the vertical monoidal product $\ot$. 

\smallskip

The analogous equivalences to $\F_3$ and $\F_4$ hold if one substitutes vertical strict for horizontal pseudonatural transformations: consider 
the 1h-cells in the double equivalence $\F''$ above (for more details, consult \cite[Section 5.4]{Fem2}), and extend to three and four variables. 
Thus we may claim:

\begin{thm} \thlabel{premon-mon-hor} 
There is a purely central horizontally premonoidal double category structure $(H, \theta, I)$ on $\Dd$ with trivial 2-cells $(u,U)$ if and only if
there is a horizontally monoidal double category structure $(\ot, \Sigma, I)$ on $\Dd$. 
\end{thm}

Here purely central means that the binoidal structure $H$ comes from a pseudodouble quasi-functor (see \prref{char df}).

\subsection{Extensions of the canonical actions}

In \cite[Definition 11]{HF} the notion of an extension of the canonical action of a monoidal bicategory on itself 
to the Kleisli bicategory of a monad was introduced. It corresponds to the notion of a {\em 0-strict morphism of actions} of 
\cite[Definition 20]{HF1}. In \cite[Definition 10.25]{Fem3} we extended this concept to the setting 
of horizontally monoidal double categories and horizontal double monads. We extend now the latter in that we substitute 
double monads by a horizontal lax action  of a horizontally monoidal double category. 

\smallskip

First of all, by a {\em horizontal icon} we mean a horizontal transformation such that its all 1h-cell components are identities. 

\begin{defn} \delabel{left ext}
Let $\Mm$ and $\Dd$ be horizontally monoidal double categories with a horizontal lax action $(F, \beta,\nu)$, and  
let $(\Dd, \ot, \alpha, \lambda, \rho; p,l,m)$ be the notation for the monoidal structure on $\Dd$.  
A {\em left horizontal extension of the canonical action of $\Dd$ on itself along $P$} is a pair $(\rtr, \theta)$, where \vspace{-0,2cm}
\begin{enumerate}
\item $\rtr:\Dd\times\coPara_\Mm(\Dd)\to\coPara_\Mm(\Dd)$ 
is a horizontal action with structure pseudonatural transformations $\tilde{\mathfrak a}, \tilde{\mathfrak l}$ and structure horizontal modifications 
$\mathfrak p, \mathfrak m, \mathfrak l$, and \vspace{-0,2cm}
\item $\theta:-\rtr P(-)\Rightarrow P(-\ot-)$ is an invertible horizontal icon, 
\end{enumerate} 
so that $\theta_{1_A,1_B}=\rtr^0_{A,B}$ (the unitor of the action on $A\rtr B$) and the identities 
\begin{equation} \eqlabel{hor-p,l,m-teta}
\frac{\theta_{\lambda_A, B}}{P(l_{AB})}=\mathfrak l_{AB}, \qquad 
\threefrac{[\Id_{P(\alpha_{A,I,B})} \vert\theta_{A,\lambda_B}]}{P(m_{AB})}{(\theta_{\rho_A, B})^{-1}}=\mathfrak m_{AB}, \qquad 
\frac{[\theta_{P(\alpha_{A,B,C}),D} \vert P(\alpha_{A,BC,D})\vert\theta_{A,P(\alpha_{B,C,D})}]}{P(p_{ABCD})}=
\mathfrak p_{ABCD}
\end{equation}
\begin{equation} \eqlabel{hor-lambda-teta}
\frac{[\theta_{1_I,f}\vert\Id_{P(\lambda_{A'})}]}{P(\lambda_f)}=\tilde{\mathfrak l}_{P(f)}: \quad
\tilde{\mathfrak l}_{A'}(I\rtr P(f))\Rightarrow P(f)\tilde{\mathfrak l}_A
\end{equation} 
\begin{equation} \eqlabel{hor-alfa-teta}
\frac{[\theta_{fg,h}\vert\Id_{P(\alpha_{A',B',C'})}]}{P(\alpha_{f,g,h})}=
\threefrac{\tilde{\mathfrak a}_{f,g,P(h)}}{[\Id_{\alpha_{A,B,C}}\vert f\rtr\theta_{g,h}]}{[\Id_{\alpha_{A,B,C}}\vert 
\theta_{f,gh}]}: \quad \alpha_{A',B',C'}\big( (fg)P(h)\big) \Rightarrow P(f(gh))\alpha_{A,B,C}
\end{equation} 
hold, whereby $P(\lambda_A)=\tilde{\mathfrak l}_A$ and $P(\alpha_{A,B,C})=\tilde{\mathfrak a}_{A,B,C}$. 
\end{defn}

A {\em right} horizontal extension of the canonical action of $\Dd$ on itself along $P$ is analogously defined as a pair $(\ltr, \xi)$, 
where  $\ltr:\coPara_\Mm(\Dd)\times\Dd\to\coPara_\Mm(\Dd)$ is a right action and $\xi: P(-)\ltr-\Rightarrow P(-\ot-)$ is an invertible 
horizontal icon. 

\smallskip

The proof of \cite[Theorem 10.31]{Fem3} for horizontal double monads (including the preparatory Section 10.4.3 from {\em loc.cit.}) 
passes {\em mutatis mutandi} to the setting of the coPara double category. Namely, we have:

\begin{thm} \thlabel{iff}
%Let $(\Dd, \ot, \alpha, \lambda, \rho)$ be a horizontally monoidal double category and let $S$ be a horizontal double monad on $\Dd$. 
Let $\Mm$ and $\Dd$ be horizontally monoidal double categories with a horizontal lax action  $(F, \beta,\nu)$. 
There is a one-to-one correspondence between left horizontal strengths on $F$ and extensions of the canonical action of $\Dd$ on itself (given by horizontal actions $\rtr:\Dd\times\coPara_\Mm(\Dd)\to\coPara_\Mm(\Dd)$ and horizontal icons $\theta:-\rtr P(-)\Rightarrow P(-\ot-)$). 
\end{thm}

From \cite[Section 10.4.3]{Fem3} we highlight the following notational conventions, originating from \cite[Proposition 4.3]{UV}, passing them to the setting of the coPara double category. We denote by $1_{MA}^\bullet:MA\to A$ a 1h-cell in $\coPara_\Mm(\Dd)$ given by $1_{MA}$ in $\Dd$. Moreover, 
by $1_A\rtr 1^\bullet_{MB}: A\rtr MB \to M(A\rtr B)$, and similarly by $1^\bullet_{MA}\ltr 1_B: MA\ltr B \to M(A\ltr B)$, we will denote 1h-cells in $\Dd$ 
determining the corresponding 1h-cells $A\rtr MB \to A\rtr B$ and $MA\ltr B \to A\ltr B$ in $\coPara_\Mm(\Dd)$, respectively. 

\smallskip

%Then for 1h-cells $f:A\to A'$ in $\Dd$ and $g^K:B\to B'$ in $\Kk l(T)$ we have that $(f\rtr g^K)^K: A\rtr B\to A'\rtr B'$ in $\Kk l(T)$ 
%is given by a 1h-cell $A\rtr B\to T(A'\rtr B')$ in $\Dd$ that can be understood as the composition 
%\begin{equation} \eqlabel{rtr explain}
%f\rtr g^K=\big(A\rtr B \stackrel{f\rtr g}{\longrightarrow} A'\rtr T(B') \stackrel{1_{A'}\rtr 1^\bullet_{T(B')}}{\xrightarrow{\hspace{1cm}}} T(A'\rtr B')\big).
%\end{equation}

The above result (and its right-hand side version) is about left and right actions of $\Dd$ on $\coPara_\Mm(\Dd)$. We will next investigate extensions of the monoidal structures from $\Dd$ to $\coPara_\Mm(\Dd)$ and in \ssref{mon Para} their correlation with left and right strengths.

\begin{defn} 
We will say that a monoidal structure on $\coPara_\Mm(\Dd)$ {\em extends} the monoidal structure of $\Dd$ if there are
\begin{enumerate}
\item a horizontally monoidal product $\ot^P:\coPara_\Mm(\Dd)\times\coPara_\Mm(\Dd)\to\coPara_\Mm(\Dd)$ 
with structure pseudonatural transformations $\tilde{\mathfrak a}, \tilde{\mathfrak l}, \tilde{\mathfrak r}$ and structure horizontal modifications $\mathfrak p, \mathfrak m, \mathfrak l, \mathfrak r$, and \vspace{-0,2cm}
\item an invertible horizontal icon $\theta^P:P(-)\ot^P P(-)\Rightarrow P(-\ot-)$, 
\end{enumerate} 
obeying axioms similar as in \deref{left ext}. 
\end{defn}

\begin{prop} \prlabel{ext-Freyd}
%\textcolor{rojo}{
If a monoidal structure on $\coPara_\Mm(\Dd)$ extends the monoidal structure of $\Dd$, then it determines a left and a right extension of the canonical action 
of $\Dd$ on itself along $P$, so that it is a biaction (there exists a horizontal equivalence $\kappa: (-\rtr P(-))\ltr-\Rightarrow -\rtr(P(-)\ltr-)$ obeying 
five axioms as in a Freyd action from \cite[Definition 22]{HF1} for 1h-cells, and the analogous five axioms for 1v-cells). %}
\end{prop}

\begin{proof}
Given a monoidal structure $\ot^P$ on $\coPara_\Mm(\Dd)$, the left and right actions $\rtr$ and $\ltr$ are induced by precomposing with $P$ at a suitable variable. 
Similarly as in \prref{horiz str-icon}, the two icons are defined via the unitors of $\rtr$ and $\ltr$ by 
$$\theta_{f,g}^L=
\scalebox{0.82}{
\bfig
 \putmorphism(-650,210)(1,0)[AB` A'B'`f\ot g]{640}1a %((AB)C)T(D') 
\putmorphism(-660,200)(0,-1)[\phantom{Y_2}` `=]{380}1r
\putmorphism(-50,200)(0,-1)[\phantom{Y_2}` `]{380}1r
\putmorphism(-50,200)(0,-1)[\phantom{Y_2}` `=]{380}0r
\putmorphism(80,210)(1,0)[``A'\w\ot\nu_{B'}]{410}1a
\putmorphism(650,210)(1,0)[A'(JB')`J(A' B')`1_{A'}\rtr 1_{JB'}^\bullet]{730}1a

\putmorphism(1600,200)(0,-1)[\phantom{Y_2}``]{380}0r
\putmorphism(1370,200)(0,-1)[\phantom{Y_2}` `=]{380}1r
\putmorphism(-670,-200)(1,0)[AB `A'B'`f\ot g]{620}1a %%%
\putmorphism(20,-200)(1,0)[\phantom{Y}` J(A' B')  `\nu_{A' B'}]{1360}1a

\put(-500,30){\fbox{$\Id_{f\ot g}$}}
\put(530,20){\fbox{$\rtr^0_{A',B'}$}}
\efig}
$$
and
$$
\theta_{f,g}^R=
\scalebox{0.82}{
\bfig
 \putmorphism(-650,210)(1,0)[AB` A'B'`f\ot g]{640}1a %((AB)C)T(D') 
\putmorphism(-660,200)(0,-1)[\phantom{Y_2}` `=]{380}1r
\putmorphism(-50,200)(0,-1)[\phantom{Y_2}` `]{380}1r
\putmorphism(-50,200)(0,-1)[\phantom{Y_2}` `=]{380}0r
\putmorphism(80,210)(1,0)[``\nu_{A'}\w\ot B']{410}1a
\putmorphism(650,210)(1,0)[(JA')B'`J(A' B')`1_{JA'}^\bullet\ltr 1_{B'}]{730}1a

\putmorphism(1600,200)(0,-1)[\phantom{Y_2}``]{380}0r
\putmorphism(1370,200)(0,-1)[\phantom{Y_2}` `=]{380}1r
\putmorphism(-670,-200)(1,0)[AB `A'B'`f\ot g]{620}1a %%%
\putmorphism(20,-200)(1,0)[\phantom{Y}` J(A' B').  `\nu_{A' B'}]{1360}1a

\put(-500,30){\fbox{$\Id_{f\ot g}$}}
\put(580,20){\fbox{$\ltr^0_{A',B'}$}}
\efig}
$$
The associativity equivalence $\tilde{\mathfrak a}$ of $\coPara_\Mm(\Dd)$ induces by restriction to $\Dd\times\Dd\times\Dd$ an equivalence 
with 2-cell components 
$$
\scalebox{0.82}{
\bfig
\putmorphism(-650,200)(0,-1)[\phantom{Y_2}` `]{380}1r
\putmorphism(-650,200)(0,-1)[\phantom{Y_2}` `=]{380}0r
\putmorphism(-600,210)(1,0)[(AB)C`(A'(MB'))C'`(f\rtr P(g))\ltr h]{1160}1a
\putmorphism(810,210)(1,0)[`A'((MB')C')`\alpha_{A',MB',C'}]{730}1a

\putmorphism(1900,200)(0,-1)[\phantom{Y_2}``]{380}0r
\putmorphism(1540,200)(0,-1)[\phantom{Y_2}` `=]{380}1r
\putmorphism(-620,-200)(1,0)[(AB)C` A(BC) `\alpha_{A,B,C}]{900}1a
\putmorphism(410,-200)(1,0)[` A'((MB')C') `f\rtr (P(g)\ltr h)]{1100}1a

\put(280,20){\fbox{$\kappa_{f, P(g),h}$}}
\efig}
$$
in $\coPara_\Mm(\Dd)$. 
\qed\end{proof}

\subsection{Monoidality of coPara} \sslabel{mon Para}

Relying on equivalence results for bicategorical monads Theorems 7.3-7.5 from \cite{HF2}, we are going to prove that the assertion of \prref{lax monoidal-comm h} extends to include two more equivalent properties. Namely: 

\begin{thm} \thlabel{main eq}
Let $\Mm$ and $\Dd$ be horizontally monoidal double categories with a horizontal lax action  $(F, \beta,\nu)$ %$F:\Mm\times\Dd\to\Dd$ 
and assume that $\Mm$ is braided with a double braiding $\Phi$. The following are equivalent: 
\begin{enumerate}
\item there is a structure $(F, F^2,F^0)$ of a lax monoidal horizontal lax action, \vspace{-0,2cm}
\item there is a structure $(F, s,t, q,c)$ of a bistrong and commutative horizontal lax action, \vspace{-0,2cm}
\item there is a purely central horizontally premonoidal structure on $\coPara_\Mm(\Dd)$ with trivial 2-cells $(u,U)$, \vspace{-0,2cm}
\item there is a horizontally monoidal structure on $\coPara_\Mm(\Dd)$ extending that on $\Dd$.
\end{enumerate}
\end{thm}

\begin{proof}
%We will prove this theorem by showing the implications $(2\Rightarrow 3), (3\Rightarrow 4)$ and $(4\Rightarrow 2)$. %of the second and third point. 
We first show the implication $(2\Rightarrow 3)$. 
In \thref{Kl-prem} we proved that a horizontal lax bistrong action $(\crta\ot, \Mm, \Dd,\beta, \nu; s, t,q)$ yields a horizontally premonoidal pseudodouble category 
structure on $\coPara_\Mm(\Dd)$. We resume the origin of the pieces of a premonoidal structure induced from such action in the following Table: 
\begin{table}[H]
\begin{center}
\begin{tabular}{ c c } %  m{7em} | m{3.4cm}
define premonoidal &  from lax bistrong action \\ \hline\hline
%binoidal  & $t,s$ \\ \hline %[1ex]   
$A\ltimes-$  & $t^-_{A,-}(A\ot-)$ \\ \hline %[1ex]   
$-\rtimes B$  & $s^-_{-,B}(-\ot B)$ \\ \hline %[1ex]   
%three $\crta\alpha$'s  & $\alpha, y, y',q$ \\ \hline %[1ex]   
%$\crta\lambda, \crta\rho$ & $\lambda, \rho, x$ \\ \hline %[1ex]   
$\overline\alpha_{-,B,C}$  & $\alpha_{-,B,C}, y'$ \\ \hline %[1ex]   
$\overline\alpha_{A,B,-}$  & $\alpha_{A,B,-}, y$ \\ \hline %[1ex]   
$\overline\alpha_{A,-,C}$  & $\alpha_{A,-,C}, q$ \\ \hline %[1ex]   
$\overline\lambda$ & $\lambda, x$ \\ \hline
$\overline\rho$ & $\rho, x'$ \\ \hline
centrality transf. & $s,t$ \\ \hline %[1ex]   % p of [HF]
\end{tabular}
\caption{Defining premonoidal structure from lax bistrong action}
\label{table:3}
\end{center}
\end{table} \vspace{-0,5cm}
\noindent %We add that $w,z$ from the left and right strength are defined in terms of $\beta, \nu$ from the lax action. 
%
%According to \thref{premon-mon-hor}, $\coPara_\Mm(\Dd)$ is horizontally monoidal 
%(with structural transformations $\alpha, \lambda, \rho$ and modifications $p,m,l,r$) 
%if and only if 
%the binoidal structure defined in \equref{bino} and below defines a pseudodouble quasi-functor. 
%it is  horizontally premonoidal whereby the binoidal structure is given by a pseudodouble quasi-functor. 
%
We will prove (2$\Rightarrow$ 3) by showing that commutativity $c$ makes the binoidal structure defined in \equref{bino} and below purely central. By 
\cite[Theorem 6.2]{Fem3} (recalled in \thref{lr central} in the Appendix) for pure centrality it is enough to prove points ii)-v) from 3. of \thref{lr central}. The desired coincidences of 1- and 2-cells from point 3. ii)-iv) we set by definition via the known images of the pseudodouble functors $A\ltimes-$ and $-\rtimes B$. We now define the four types of 2-cells as in item 2. of \prref{char df}, and then it will remain to show that we have the desired transformations 
(determined by those 2-cells) and modifications from point 3. of \thref{lr central}. For $U\ltimes-\vert_u=(-\rtimes u\vert_U)^{-1}$ we set to be identities. A desired 2-cell 
$$
\scalebox{0.86}{
\bfig
 \putmorphism(-150,50)(1,0)[AB`AB'`A\ltimes g]{600}1a
 \putmorphism(450,50)(1,0)[\phantom{A\ot B}`A'B' `f\rtimes B']{680}1a
\putmorphism(-180,50)(0,-1)[\phantom{Y_2}``=]{450}1r
\putmorphism(1100,50)(0,-1)[\phantom{Y_2}``=]{450}1r
\put(350,-190){\fbox{$(f,g)$}}
 \putmorphism(-150,-400)(1,0)[AB`A'B`f\rtimes B]{600}1a
 \putmorphism(450,-400)(1,0)[\phantom{A\ot B}`A'B' `A'\ltimes g]{680}1a
\efig}
$$ 
in $\coPara_\Mm(\Dd)$ %has the form of a horizontally globular 2-cell 
%$$AB\stackrel{Ag}{\to} A(MB') \stackrel{t^M_{A,B'}}{\to} M(AB')\stackrel{M(fB')}{\to} M((NA')B') \stackrel{Ms^N_{A',B'}}{\to} M(N(A'B'))
%\stackrel{\beta_{M,N,A'B'}}{\to}(MN)(A'B') \vspace{-0,4cm}$$ 
%$$\hspace{9cm} \stackrel{\Phi_{M,N}(A'B')}{\to} (NM)(A'B') $$
%$$\Downarrow$$
%$$AB\stackrel{fB}{\to} (NA')B \stackrel{s^N_{A',B}}{\to} N(A'B)\stackrel{N(A'g)}{\to} N(A'(MB')) \stackrel{Nt^M_{A',B'}}{\to} N(M(A'B'))
%\stackrel{\beta_{N,M,A'B'}}{\to}(NM)(A'B')$$
%in $\Dd$. We set it to be
we set to be the 2-cell
$$ 
\scalebox{0.78}{ 
\bfig
 \putmorphism(370,450)(1,0)[A(MB')` M(AB') `t^M_{A,B'}]{830}1a
 \putmorphism(1360,450)(1,0)[` M((NA')B')  `M(fB')]{730}1a

\putmorphism(410,450)(0,-1)[\phantom{Y_2}``=]{450}1l

\put(980,220){\fbox{$(t^M)^{-1}_{f,1_{B'}}$}}
\put(240,-240){\fbox{$ps \ot$}}

\putmorphism(-260,0)(1,0)[AB` A(MB')`A g]{630}1a 
\putmorphism(-260,0)(0,-1)[\phantom{Y_2}` `=]{430}1r 
\putmorphism(1220,0)(0,-1)[\phantom{Y_2}``=]{450}1l
\putmorphism(360,0)(1,0)[\phantom{A''\ot B'}`(NA')(MB') ` f(MB')]{840}1a%%%%%
 \putmorphism(1280,0)(1,0)[\phantom{A''\ot B'}`  `t^M_{NA',B'}]{520}1a %T(T_A(T_BT_C))
 \putmorphism(2150,0)(1,0)[\phantom{A''\ot B'}` M(N(A'B')) `M(s^N_{A',B'})]{820}1a %T(T_A(T_BT_C))
 \putmorphism(3070,0)(1,0)[\phantom{A''\ot B'}`(MN)(A'B')  `\beta_{M,N,A'B'}]{780}1a %T(T_A(T_BT_C))
 \putmorphism(3960,0)(1,0)[\phantom{A''\ot B'}`(NM)(A'B')  `\Phi_{M,N}(A'B')]{800}1a %T(T_A(T_BT_C))

\putmorphism(2060,450)(0,-1)[\phantom{Y_2}`M((NA')B') ` =]{450}1r 
\putmorphism(4680,0)(0,-1)[\phantom{Y_2}`(MN)(A''B'')`=]{450}1r 
\put(2460,-250){\fbox{$(c^{NA',MB'})^{-1}$}}

\putmorphism(-260,-450)(1,0)[AB` (NA')B`f B]{630}1a 
\putmorphism(360,-450)(1,0)[\phantom{A''\ot B'}`(NA')(MB') ` (NA')g]{800}1a%%%%%
 \putmorphism(1210,-450)(1,0)[\phantom{A''\ot B'}` N(A'(MB')) `s^N_{A',MB'}]{940}1b 
 \putmorphism(2230,-450)(1,0)[\phantom{A''\ot B'}` N(M(A'B')) `Nt^M_{A',B'}]{940}1b 
 \putmorphism(3260,-450)(1,0)[\phantom{A''\ot B'}` (NM)(A'B') `\beta_{N,M,A'B'}]{890}1b 

\putmorphism(410,-450)(0,-1)[\phantom{Y_2}`(NA')B`=]{450}1l
\putmorphism(2060,-450)(0,-1)[\phantom{Y_2}``=]{450}1l
\put(1020,-670){\fbox{$s^N_{1_{A'},g}$}}
\putmorphism(380,-900)(1,0)[\phantom{A''\ot B'}`N(A'B) ` s^N_{A',B'}]{760}1a%%%%%
 \putmorphism(1150,-900)(1,0)[\phantom{A''\ot B'}` N(A'(MB')) `N(A'g)]{860}1a 
\efig}
$$
in $\Dd$, and for the remaining desired 2-cells we set 
$$
\scalebox{0.86}{
\bfig
\putmorphism(-150,250)(1,0)[A\rtimes B`A'\rtimes B`f\rtimes B]{600}1a
\putmorphism(-150,-200)(1,0)[A\rtimes\tilde B`A'\rtimes\tilde B `f\rtimes\tilde B]{640}1a
\putmorphism(-180,250)(0,-1)[\phantom{Y_2}``A\rtimes u]{450}1l
\putmorphism(450,250)(0,-1)[\phantom{Y_2}``A'\rtimes u]{450}1r
\put(0,30){\fbox{$f\ltimes u$}}
\efig}
:=
\scalebox{0.82}{
\bfig
\putmorphism(-710,170)(1,0)[AB``fB ]{460}1a
\putmorphism(-710,200)(0,-1)[``Au]{410}1l 
\putmorphism(-110,170)(1,0)[(NA')B`N(A'B)`s^N_{A',B} ]{900}1a
\putmorphism(-160,200)(0,-1)[``]{410}1r 
\putmorphism(-180,200)(0,-1)[``(NA')u]{410}0r 
\putmorphism(830,200)(0,-1)[``N(A'u)]{410}1r
\putmorphism(-710,-200)(1,0)[A\tilde B``f\tilde B ]{460}1b %\tilde A(MB')
\putmorphism(-110,-200)(1,0)[(NA')\tilde B`N(A'\tilde B)`s^N_{A',\tilde B}]{900}1b
\put(-530,-20){\fbox{$fu$}}
\put(250,-20){\fbox{$s^{N,1_{A'},u}$}}
\efig}
$$
and 
$$\scalebox{0.86}{
\bfig
\putmorphism(-150,250)(1,0)[A\ltimes B`A\ltimes B'`A\ltimes g]{600}1a
\putmorphism(-150,-200)(1,0)[ \tilde A\ltimes B`\tilde A\ltimes B' `\tilde A\ltimes g]{640}1a
\putmorphism(-180,250)(0,-1)[\phantom{Y_2}``U\rtimes B]{450}1l
\putmorphism(450,250)(0,-1)[\phantom{Y_2}``U\rtimes B']{450}1r
\put(0,30){\fbox{$U\ltimes g$}}
\efig}
:=
\scalebox{0.82}{
\bfig
\putmorphism(-710,170)(1,0)[AB``Ag ]{460}1a
\putmorphism(-710,200)(0,-1)[``UB]{410}1l 
\putmorphism(-110,170)(1,0)[A(MB')`M(AB')`t^M_{A,B'} ]{900}1a 
\putmorphism(-160,200)(0,-1)[``U(MB')]{410}1r 
\putmorphism(830,200)(0,-1)[``M(UB')]{410}1r
\putmorphism(-710,-200)(1,0)[\tilde AB``\tilde Ag]{460}1b
\putmorphism(-110,-200)(1,0)[\tilde A(MB')`M(\tilde AB').`t^M_{\tilde A,B'}]{900}1b
\put(-540,-20){\fbox{$Ug$}}
\put(250,-20){\fbox{$t^{M,U,1_{B'}}$}}
\efig}
$$
The axioms \axiomref{h.o.t.-3} and \axiomref{h.o.t.-4} for $f\ltimes-$ and $-\rtimes g$ are easily seen to hold. The proof of 
the axioms \axiomref{v.l.t.\x 1} and \axiomref{v.l.t.\x 2} for $U\ltimes-$ and $-\rtimes u$ are somewhat more laborous, but direct 
(they involve the structural 2-cells $w$ (resp. $w'$), those of $t$ (resp. $s$), the lax and colax compositors of $-\ot-$, and $\beta$). 
The axioms \axiomref{v.l.t.\x 3} and \axiomref{v.l.t.\x 4} trivially hold because here the 2-cells $U\ltimes-\vert_u$ are trivial, 
and \axiomref{v.l.t.\x 5} holds by naturality of $\ot$ in $\Dd$. The latter axiom is simultaneously \axiomref{m.ho-vl.-2} for 
$a\ltimes-$ and $-\rtimes b$, which we need to prove to be modifications. Similarly, the axiom \axiomref{m.ho-vl.-1} for them to be modifications 
is simultaneously the axiom \axiomref{h.o.t.-5} for $f\ltimes-$ and $-\rtimes g$. Its proof is straightforward: one applies naturality of $\ot$, 
the axiom \axiomref{h.o.t.-5} of $s^N$, and the modification axiom \axiomref{m.ho-vl.-1} of $c^{-1}$. 
The only remaining axioms to be proved are \axiomref{h.o.t.-1} and \axiomref{h.o.t.-2} for $f\ltimes-$ and $-\rtimes g$. They are entirely contained 
in the underlying bicategory and they hold similarly as in \cite[Theorem 7.3]{HF2}.

Assume 3. Then by \thref{premon-mon-hor} we know that $\coPara_\Mm(\Dd)$ has a monoidal product $\ot^P$. %\textcolor{rojo}{
The lax structure of the pseudodouble functor $\ot^P$ determines an invertible horizontal icon $\theta^P: P(-)\ot^P P(-)\Rightarrow P(-\ot-)$ 
as follows
$$
\theta_{f,g}=
\scalebox{0.82}{
\bfig
 \putmorphism(-650,210)(1,0)[AB` A'B'`f\ot g]{640}1a %((AB)C)T(D') 
\putmorphism(-660,200)(0,-1)[\phantom{Y_2}` `=]{380}1r
\putmorphism(-50,200)(0,-1)[\phantom{Y_2}` `]{380}1r
\putmorphism(-50,200)(0,-1)[\phantom{Y_2}` `=]{380}0r
\putmorphism(60,210)(1,0)[``\nu_{A'}\w\ot \nu_{B'}]{480}1a
\putmorphism(740,210)(1,0)[(JA')(JB') `J(A' B')`\nabla_{A',B'}]{840}1a

\putmorphism(1540,200)(0,-1)[\phantom{Y_2}` `=]{380}1r
\putmorphism(-670,-200)(1,0)[AB `A'B'`f\ot g]{620}1a %%%
\putmorphism(20,-200)(1,0)[\phantom{Y}` J(A' B')  `\nu_{A' B'}]{1610}1a

\put(-500,30){\fbox{$\Id_{f\ot g}$}}
\put(450,0){\fbox{$(-\ot^P-)^0_{A',B'}$}}
\efig}
$$
%Namely, we need $\theta^P_{f,g}: \nu_{A'}f\ot^P \nu_{B'}g \Rightarrow\nu_{A\ot B}(f\ot g)$... Problem: $3\Rightarrow 4$. 
so that $\ot^P$ extends $\ot$ of $\Dd$. Here $\nabla_{A',B'}$ is a 1h-cell in $\Dd$ such that $1_{A'}\ot^P 1_{B'}$ in $\coPara_\Mm(\Dd)$ is determined by 
$\nabla_{A',B'}(\nu_{A'}\w\ot \nu_{B'})$ (we skip technical details concerning how pseudodouble functor structure of $\ot^P$ is defined out of the binoidal structure, in \cite[Section 5.1]{Fem3}). 

\medskip

To prove $(4\Rightarrow 2)$, we start from a horizontally monoidal structure $\ot^P$ on $\coPara_\Mm(\Dd)$ extending $\ot$ of $\Dd$. 
By \prref{ext-Freyd} we have a left and a right extension of the canonical action of $\Dd$ on itself, and by \thref{iff} and its right-hand 
side version there are a left and a right strength $t$ and $s$ of the action of $\Mm$ on $\Dd$. We define a bistrength $q$ via the equivalence 
$\kappa$ from \prref{ext-Freyd}, specifically as its 2-cells components 
$$q:=
\scalebox{0.82}{
\bfig
\putmorphism(-650,200)(0,-1)[\phantom{Y_2}` `]{380}1r
\putmorphism(-650,200)(0,-1)[\phantom{Y_2}` `=]{380}0r
\putmorphism(-600,210)(1,0)[(A(MB))C`(M(AB))C`(1_A\rtr 1^\bullet_{MB})\ot 1_C]{1260}1a
\putmorphism(900,210)(1,0)[`M((AB)C)`1_{AB}^\bullet\ltr 1_C]{730}1a %\alpha_{A,MB,C}
\putmorphism(1830,210)(1,0)[`M(A(BC))`M\alpha_{A,B,C}]{780}1a 

\putmorphism(2800,200)(0,-1)[\phantom{Y_2}``]{380}0r
\putmorphism(2600,200)(0,-1)[\phantom{Y_2}` `=]{380}1r
\putmorphism(-620,-200)(1,0)[(A(MB))C` A((MB)C) `\alpha_{A,MB,C}]{970}1a
\putmorphism(570,-200)(1,0)[` A(M(BC)) `A\ot(1_{MB}^\bullet\ltr 1_C)]{1060}1a
\putmorphism(1840,-200)(1,0)[` M(A(BC)) `1_A\rtr 1_{BC}^\bullet]{840}1a

\put(580,20){\fbox{$\kappa_{1_A, 1^\bullet_{MB},1_C}$}}
\efig}
$$
in $\Dd$. With a similar interpretation of the 1h-cells in $\coPara_\Mm(\Dd)$ via 1h-cells in $\Dd$ as in the edges of $q$ above, 
one sees that the following makes sense to be taken for a commutativity $c$ 
$$c^{MA,NB}:=
\scalebox{0.8}{
\bfig
\putmorphism(260,200)(1,0)[`A(MB)`1^\bullet_{NA}\ltr 1_{MB}]{740}1a
 \putmorphism(1040,200)(1,0)[\phantom{F(A)}`AB `1_{A}\rtr 1^\bullet_{MB}]{620}1a
 \putmorphism(320,-100)(1,0)[`AB`1^\bullet_{NA}\ot^P 1^\bullet_{MB}]{1340}1b
  \putmorphism(170,-500)(1,0)[\phantom{F(A)}`(NA)B ` 1_{NA}\rtr 1^\bullet_{MB}]{840}1b
 \putmorphism(1060,-500)(1,0)[\phantom{F(A)}`AB ` 1^\bullet_{NA}\ltr 1_B]{620}1b

\putmorphism(60,200)(0,-1)[(NA)(MB)`(NA)(MB) `=]{300}1l
\putmorphism(60,-100)(0,-1)[`(NA)(MB) `=]{400}1l
\putmorphism(1650,200)(0,-1)[\phantom{Y_2}``=]{300}1r
\putmorphism(1650,-100)(0,-1)[\phantom{Y_2}``=]{400}1r
\put(870,0){\fbox{$lax\s \ot^P$}}
\put(870,-370){\fbox{$colax\s \ot^P$}}
\efig}
$$
whereby we restrict the monoidal product $\ot^P$ to the left and the right action of $\Dd$. 
That $q$ and $c$ are indeed horizontal modifications it is proved as in (the converse of) the bicategorical \cite[Theorem 7.3]{HF2} (the axiom 
\axiomref{m.ho.-2} easily follows). This completes the proof. 
%\textcolor{rojo}{
%We add that $w,z$ from the left and right strength are defined in terms of $\beta, \nu$ from the lax action. }
\qed\end{proof}

The implication $(1\Rightarrow 4)$ is a double categorical version of the 1-categorical \cite[Corollary 5.21]{CJM2}. 

\bigskip

We developed the results on the Para construction for double categories with the interest to extend some known constructions for Para 
in computer science. For this purpose our approach was to deal with horizontal structures. As for the vertical structures, 
our results from \ssref{act Para} have their vertical analogous (similarly as in Theorems 10.17 and 10.32 of \cite{Fem3}): 
vertical strengths induce vertical actions, vertical bistrengths induce vertical premonoidality. The 1-1 correspondence between 
horizontal strengths and horizontal extensions of the canonical actions of $\Dd$ on itself from \thref{iff} does not hold in the vertical direction: 
this failure already occurs for the Kleisli double category, as we discussed in \cite[Remark 10.27]{Fem3}. Moreover, a vertical strength $t$ 
induces a vertical action (and a ``vertical'' extension) via the induced horizontal strength $\hat t$, but in the converse one 
can not recover a vertical strength out of a horizontal one. Our main \thref{main eq} does not hold in its fulness in the vertical direction. 
We do have the equivalence of vertical lax monoidality and vertically bistrong commutativity of the action $(F,\beta,\nu)$ (see \prref{F is lax monoidal}), 
but commutativity, now expressed as a vertical identity modification, does not induce globular horizontal 2-cells $(f,g)$ necessary for pure centrality in point 3. 
We do have vertical premonoidality (by \thref{Kl-prem}), but the absence of pure centrality implies that neither do we obtain vertical monoidality 
(recall \thref{premon-mon}). %***

\pagebreak

{\large {\bf Appendix}}

\bigskip

We record here the ``pseudo'' version of the ``lax'' statement of \cite[Proposition 3.3]{Fem2}. 
The data of the points 1. and 2. in the proposition comprise the definition of a {\em pseudodouble quasi-functor} $H:\Aa\times\Bb\to\Cc$.

\begin{prop} \prlabel{char df} 
A pseudodouble functor $\F\colon\Aa\to\Lax_{hop}(\Bb, \Cc)$ of double categories consists of the following: \\
1. pseudodouble functors 
$$(-,A)\colon\Bb\to\Cc\quad\text{ and}\quad (B,-)\colon\Aa\to\Cc$$ 
such that $(-,A)\vert_B=(B,-)\vert_A=(B,A)$, 
for objects $A\in\Aa, B\in\Bb$, \\
%\item 
2. 2-cells
$$
\scalebox{0.86}{
\bfig
 \putmorphism(-150,50)(1,0)[(B,A)`(B', A)`(k, A)]{600}1a
 \putmorphism(450,50)(1,0)[\phantom{A\ot B}`(B', A') `(B', K)]{680}1a
\putmorphism(-180,50)(0,-1)[\phantom{Y_2}``=]{450}1r
\putmorphism(1100,50)(0,-1)[\phantom{Y_2}``=]{450}1r
\put(350,-190){\fbox{$(k,K)$}}
 \putmorphism(-150,-400)(1,0)[(B,A)`(B,A')`(B,K)]{600}1a
 \putmorphism(450,-400)(1,0)[\phantom{A\ot B}`(B', A') `(k, A')]{680}1a
\efig}
$$

$$
\scalebox{0.86}{
\bfig
\putmorphism(-150,50)(1,0)[(B,A)`(B,A')`(B,K)]{600}1a
\putmorphism(-150,-400)(1,0)[(\tilde B, A)`(\tilde B,A') `(\tilde B,K)]{640}1a
\putmorphism(-180,50)(0,-1)[\phantom{Y_2}``(u,A)]{450}1l
\putmorphism(450,50)(0,-1)[\phantom{Y_2}``(u,A')]{450}1r
\put(0,-180){\fbox{$(u, K)$}}
\efig}
\quad
\scalebox{0.86}{
\bfig
\putmorphism(-150,50)(1,0)[(B,A)`(B',A)`(k,A)]{600}1a
\putmorphism(-150,-400)(1,0)[(B, \tilde A)`(B', \tilde A) `(k,\tilde A)]{640}1a
\putmorphism(-180,50)(0,-1)[\phantom{Y_2}``(B,U)]{450}1l
\putmorphism(450,50)(0,-1)[\phantom{Y_2}``(B',U)]{450}1r
\put(0,-180){\fbox{$(k,U)$}}
\efig}
$$

$$
\scalebox{0.86}{
\bfig
 \putmorphism(-150,500)(1,0)[(B,A)`(B,A) `=]{600}1a
\putmorphism(-180,500)(0,-1)[\phantom{Y_2}`(B, \tilde A) `(B,U)]{450}1l
\put(-20,50){\fbox{$(u,U)$}}
\putmorphism(-150,-400)(1,0)[(\tilde B, \tilde A)`(\tilde B, \tilde A) `=]{640}1a
\putmorphism(-180,50)(0,-1)[\phantom{Y_2}``(u,\tilde A)]{450}1l
\putmorphism(450,50)(0,-1)[\phantom{Y_2}``(\tilde B, U)]{450}1r
\putmorphism(450,500)(0,-1)[\phantom{Y_2}`(\tilde B, A) `(u,A)]{450}1r
\efig}
$$ 
in $\Cc$ for every 1h-cells $A\stackrel{K}{\to} A'$ and $B\stackrel{k}{\to} B'$ and 1v-cells $A\stackrel{U}{\to} \tilde A$ and 
$B\stackrel{u}{\to} \tilde B$ which satisfy: 

\noindent $\bullet$ \quad \axiom{($1_B,K$)}  
$$
\scalebox{0.86}{
\bfig
 \putmorphism(-210,420)(1,0)[(B,A)`(B,A)`=]{550}1a
\putmorphism(-210,50)(1,0)[(B,A)`(B,A) `(1_B,A)]{560}1a
\putmorphism(350,400)(1,0)[\phantom{F(A)}` (B,A') `(B,K)]{600}1a
\putmorphism(360,50)(1,0)[\phantom{F(A)}`(B,A') `(B,K)]{600}1a

\putmorphism(-170,420)(0,-1)[\phantom{Y_2}``=]{350}1l
\putmorphism(320,420)(0,-1)[\phantom{Y_2}``]{370}1r
\putmorphism(300,420)(0,-1)[\phantom{Y_2}``=]{370}0r
\put(-120,250){\fbox{$(-,A)_B$}}

\putmorphism(860,420)(0,-1)[\phantom{Y_2}``=]{350}1r
\put(450,240){\fbox{$\Id_{(B,K)}$}}
\putmorphism(-180,50)(0,-1)[\phantom{Y_2}``=]{350}1r
\putmorphism(860,50)(0,-1)[\phantom{Y_2}``=]{350}1r
 \putmorphism(-150,-300)(1,0)[(B,A)`(B,A')`(B,K)]{500}1a
 \putmorphism(350,-300)(1,0)[\phantom{A\ot B}`(B, A') `(1_B, A')]{600}1a
\put(170,-110){\fbox{$(1_B,K)$}}
\efig}
\quad=\quad
\scalebox{0.86}{
\bfig
 \putmorphism(-260,200)(1,0)[(B,A)`\phantom{F(A)} `(B,K)]{550}1a
\putmorphism(330,200)(1,0)[(B,A')`(B,A')`=]{600}1a
\putmorphism(-210,200)(0,-1)[\phantom{Y_2}``=]{370}1l
\putmorphism(320,200)(0,-1)[\phantom{Y_2}``]{370}1l
\putmorphism(340,200)(0,-1)[\phantom{Y_2}``=]{370}0l
\putmorphism(960,200)(0,-1)[\phantom{Y_2}``=]{370}1r
 \putmorphism(-260,-170)(1,0)[(B,A)`\phantom{Y_2}`(B,K)]{520}1a
 \put(420,30){\fbox{$(-,A')_{B}$}} 
\putmorphism(360,-170)(1,0)[(B,A')`(B,A') `(1_B,A')]{620}1a
\put(-160,30){\fbox{$\Id_{(B,K)}$}}
\efig}
$$

%\pagebreak

\noindent $\bullet$ \quad \axiom{($k,1_A$)} \vspace{-0,9cm}
$$ 
\scalebox{0.86}{
\bfig
 \putmorphism(-150,200)(1,0)[(B,A)`(B,A)`=]{500}1a
\putmorphism(360,200)(1,0)[\phantom{F(A)}`(B,A') `(k,A)]{500}1a
\putmorphism(830,200)(0,-1)[\phantom{Y_2}``=]{350}1r
\put(460,30){\fbox{$\Id_{(k,A)}$}}

\putmorphism(-180,200)(0,-1)[\phantom{Y_2}``=]{370}1l
\putmorphism(320,200)(0,-1)[\phantom{Y_2}``]{370}1r
\putmorphism(300,200)(0,-1)[\phantom{Y_2}``=]{370}0r
 \putmorphism(-150,-170)(1,0)[(B,A)`(B,A)`(B,1_A)]{500}1a
 \put(-140,30){\fbox{$(B,-)_A$}} %(0,-180)
\putmorphism(350,-170)(1,0)[\phantom{F(A)}`(B',A) `(k,A)]{560}1a
\efig}
\quad
=
\quad
\scalebox{0.86}{
\bfig
\putmorphism(-150,420)(1,0)[(B,A)` \phantom{Y_2}`(k,A)]{450}1a
\putmorphism(370,420)(1,0)[(B',A)` (B',A) `=]{500}1a
\putmorphism(-170,420)(0,-1)[\phantom{Y_2}``=]{350}1l
\putmorphism(350,420)(0,-1)[\phantom{Y_2}``=]{350}1l
\putmorphism(860,420)(0,-1)[\phantom{Y_2}``=]{350}1r
 \putmorphism(420,50)(1,0)[\phantom{Y_2}`(B',A)`(B',1_A)]{520}1a
\putmorphism(-150,50)(1,0)[(B,A)` (B',A) `(k,A)]{500}1a
 \put(410,250){\fbox{$(B',-)_A$}} %(0,-180)
\put(-120,250){\fbox{$\Id_{(k,A)}$}}

\putmorphism(-180,50)(0,-1)[\phantom{Y_2}``=]{350}1r
\putmorphism(860,50)(0,-1)[\phantom{Y_2}``=]{350}1r
\putmorphism(-150,-300)(1,0)[(B,A)`(B,A)`(B,1_A)]{500}1a
 \putmorphism(350,-300)(1,0)[\phantom{A\ot B}`(B', A) `(k, A)]{560}1a
\put(200,-120){\fbox{$(k,1_A)$}}
\efig}
$$

\noindent where the 2-cells $(-,A)_{B}$ and $(B,-)_A$ come from laxity of the lax double functors $(-,A)$ and $(B,-)$ 

\noindent $\bullet$ \quad  \axiom{($u,1_A$)} 
$$
\scalebox{0.86}{
\bfig
\putmorphism(-250,500)(1,0)[(B,A)`(B,A)` =]{550}1a
 \putmorphism(-250,50)(1,0)[(\tilde B,A)`(\tilde B,A)` =]{550}1a
 \putmorphism(-250,-400)(1,0)[(\tilde B,A)`(\tilde B,A)` (\tilde B,1_A)]{550}1a

\putmorphism(-280,500)(0,-1)[\phantom{Y_2}``(u,A)]{450}1l
 \putmorphism(-280,70)(0,-1)[\phantom{F(A)}` `=]{450}1l

\putmorphism(300,500)(0,-1)[\phantom{Y_2}``(u,A)]{450}1r
\putmorphism(300,70)(0,-1)[\phantom{Y_2}``=]{450}1r
\put(-150,270){\fbox{$Id^{(u,A)}$}}
\put(-190,-150){\fbox{$(\tilde B,-)_A$}}
\efig}
=
\scalebox{0.86}{
\bfig
\putmorphism(-250,500)(1,0)[(B,A)`(B,A)` =]{550}1a
 \putmorphism(-250,50)(1,0)[(B,A)`(B,A)` (B,1_A)]{550}1a
 \putmorphism(-250,-400)(1,0)[(\tilde B,A)`(\tilde B,A)` (\tilde B,1_A)]{550}1a

\putmorphism(-280,500)(0,-1)[\phantom{Y_2}``= ]{450}1l
 \putmorphism(-280,70)(0,-1)[\phantom{F(A)}` `(u,A)]{450}1l

\putmorphism(300,500)(0,-1)[\phantom{Y_2}``=]{450}1r
\putmorphism(300,70)(0,-1)[\phantom{Y_2}``(u,A)]{450}1r
\put(-160,290){\fbox{$(B,-)_A$}}
\put(-140,-150){\fbox{$(u,1_A)$}}
\efig}
$$

%(12)
\noindent $\bullet$ \quad \axiom{($1_B,U$)} 
$$\scalebox{0.86}{
\bfig
\putmorphism(-250,500)(1,0)[(B,A)`(B,A)` =]{550}1a
 \putmorphism(-250,50)(1,0)[(B,A)`(B,A)` (1_B,A)]{550}1a
 \putmorphism(-250,-400)(1,0)[(B,\tilde A)`(B,\tilde A)` (1_B,\tilde A)]{550}1a

\putmorphism(-280,500)(0,-1)[\phantom{Y_2}``= ]{450}1l
 \putmorphism(-280,70)(0,-1)[\phantom{F(A)}` `(B,U)]{450}1l

\putmorphism(300,500)(0,-1)[\phantom{Y_2}``=]{450}1r
\putmorphism(300,70)(0,-1)[\phantom{Y_2}``(B,U)]{450}1r
\put(-170,290){\fbox{$(-,A)_B$}}
\put(-150,-160){\fbox{$(1_B,U)$}}
\efig}\quad
=
\scalebox{0.86}{
\bfig
\putmorphism(-250,500)(1,0)[(B,A)`(B,A)` =]{550}1a
 \putmorphism(-250,50)(1,0)[(B,\tilde A)`(B,\tilde A)` =]{550}1a
 \putmorphism(-250,-400)(1,0)[(B,\tilde A)`(B,\tilde A)` (\tilde B,1_{\tilde A})]{550}1a

\putmorphism(-280,500)(0,-1)[\phantom{Y_2}``(B,U)]{450}1l
 \putmorphism(-280,70)(0,-1)[\phantom{F(A)}` `=]{450}1l

\putmorphism(300,500)(0,-1)[\phantom{Y_2}``(B,U)]{450}1r
\putmorphism(300,70)(0,-1)[\phantom{Y_2}``=]{450}1r
\put(-170,270){\fbox{$Id^{(B,U)}$}}
\put(-180,-160){\fbox{$(-, \tilde A)_B$}}
\efig}
$$

\noindent $\bullet$ \qquad \axiom{($1^B,K$)}  \quad $(1^B,K)=Id_{(B,K)}$ \qquad\hspace{-0,16cm}\text{and}\qquad 
$\bullet$ \quad  \axiom{($k,1^A$)} \qquad $(k,1^A)=Id_{(k,A)}$ %$\bullet$ \quad ($(k,1^A)$) \label{k,1uA} \quad $(k,1^A)=Id_{(k,A)}$  

%(22)
\noindent $\bullet$ \qquad  \axiom{($1^B,U$)} \quad $(1^B,U)=Id^{(B,U)}$ \qquad\text{and}\qquad \hspace{-0,22cm} 
$\bullet$ \quad  \axiom{($u,1^A$)} \qquad $(u,1^A)=Id^{(u,A)}$

%\item (11)  
\noindent $\bullet$ \quad  \axiom{($k'k,K$)} 
$$
\scalebox{0.78}{
\bfig
 \putmorphism(450,650)(1,0)[(B', A) `(B'', A) `(k', A)]{680}1a
 \putmorphism(1140,650)(1,0)[\phantom{A\ot B}`(B'', A') ` (B'', K)]{680}1a

 \putmorphism(-150,200)(1,0)[(B, A) `(B', A)`(k, A)]{600}1a
 \putmorphism(450,200)(1,0)[\phantom{A\ot B}`(B', A') `(B', K)]{680}1a
 \putmorphism(1130,200)(1,0)[\phantom{A\ot B}`(B'', A') ` (k', A')]{680}1a

\putmorphism(450,650)(0,-1)[\phantom{Y_2}``=]{450}1r
\putmorphism(1750,650)(0,-1)[\phantom{Y_2}``=]{450}1r
\put(1000,420){\fbox{$ (k',K)$}}

 \putmorphism(-150,-250)(1,0)[(B, A)`(B, A') `(B,K)]{640}1a
 \putmorphism(480,-250)(1,0)[\phantom{A'\ot B'}`(B', A') `(k, A')]{680}1a

\putmorphism(-180,200)(0,-1)[\phantom{Y_2}``=]{450}1l
\putmorphism(1120,200)(0,-1)[\phantom{Y_3}``=]{450}1r
\put(310,-50){\fbox{$ (k,K)$}}

 \putmorphism(1170,-250)(1,0)[\phantom{A\ot B}`(B'', A') ` (k', A')]{650}1a
\putmorphism(450,-250)(0,-1)[\phantom{Y_2}``=]{450}1r
\putmorphism(1750,-250)(0,-1)[\phantom{Y_2}``=]{450}1r

 \putmorphism(480,-700)(1,0)[(B, A') `(B'', A'') `(k'k, A')]{1320}1a
\put(920,-470){\fbox{$ (-,A')_{k'k}$}}
\efig}
=\quad
\scalebox{0.78}{
\bfig
 \putmorphism(-150,500)(1,0)[(B, A) `(B', A)`(k,A)]{580}1a
 \putmorphism(450,500)(1,0)[\phantom{(B, A)} `(B'', A) `(k',A)]{660}1a
\putmorphism(-180,500)(0,-1)[\phantom{Y_2}``=]{450}1r
\put(240,270){\fbox{$ (-,A)_{k'k}$}}

 \putmorphism(-150,50)(1,0)[(B,A)` `(k'k,A)]{1080}1a
 \putmorphism(1080,50)(1,0)[(B'',A)`(B'', A') ` (B'',K)]{680}1a

\putmorphism(1050,500)(0,-1)[\phantom{Y_2}``=]{450}1r

 \putmorphism(-150,-400)(1,0)[(B, A)`(B, A') `(B,K)]{640}1a
 \putmorphism(530,-400)(1,0)[\phantom{Y_2X}`(B'', A') `(k'k,A')]{1220}1a
\put(570,-200){\fbox{$ (k'k,K)$}}

\putmorphism(-180,50)(0,-1)[\phantom{Y_2}``=]{450}1l
\putmorphism(1700,50)(0,-1)[\phantom{Y_3}``=]{450}1r
\efig}
$$ 
where $(-,A)_{k'k}$ is the 2-cell from the laxity of $(-,A)$ 

\noindent $\bullet$ \quad  \axiom{($k,K'K$)}  
$$
\scalebox{0.78}{
\bfig
 \putmorphism(450,500)(1,0)[(B', A) `(B', A') `(B', K)]{680}1a
 \putmorphism(1140,500)(1,0)[\phantom{A\ot B}`(B', A'') ` (B', K')]{680}1a

 \putmorphism(-150,50)(1,0)[(B, A) `(B', A)`(k,A)]{600}1a
 \putmorphism(450,50)(1,0)[\phantom{A\ot B}`(B', A'') `(B', K'K)]{1350}1a

\putmorphism(450,500)(0,-1)[\phantom{Y_2}``=]{450}1r
\putmorphism(1750,500)(0,-1)[\phantom{Y_2}``=]{450}1r
\put(880,290){\fbox{$ (B',-)_{K'K}$}}

 \putmorphism(-150,-400)(1,0)[(B, A)`(B, A'') `(B, K'K)]{980}1a
 \putmorphism(780,-400)(1,0)[\phantom{A'\ot B'}`(B', A'') `(k,A'')]{980}1a

\putmorphism(-180,50)(0,-1)[\phantom{Y_2}``=]{450}1l
\putmorphism(1750,50)(0,-1)[\phantom{Y_3}``=]{450}1r
\put(560,-200){\fbox{$ (k,K'K)$}}

\efig}
=
\scalebox{0.78}{
\bfig
 \putmorphism(-150,450)(1,0)[(B,A)`(B',A)`(k,A)]{600}1a
 \putmorphism(450,450)(1,0)[\phantom{A\ot B}`(B', A') `(B',K)]{680}1a

 \putmorphism(-150,0)(1,0)[(B,A)`(B,A')`(B,K)]{600}1a
 \putmorphism(450,0)(1,0)[\phantom{A\ot B}`(B', A') `(k,A')]{680}1a
 \putmorphism(1120,0)(1,0)[\phantom{A'\ot B'}`(B', A'') `(B', K')]{660}1a

\putmorphism(-180,450)(0,-1)[\phantom{Y_2}``=]{450}1r
\putmorphism(1100,450)(0,-1)[\phantom{Y_2}``=]{450}1r
\put(350,210){\fbox{$(k,K)$}}
\put(1000,-250){\fbox{$(k,K')$}}

 \putmorphism(450,-450)(1,0)[\phantom{A''\ot B'}` (B, A'') `(B, K')]{680}1a
 \putmorphism(1100,-450)(1,0)[\phantom{A''\ot B'}`(B', A'') ` (k, A'')]{660}1a

\putmorphism(450,0)(0,-1)[\phantom{Y_2}``=]{450}1l
\putmorphism(1750,0)(0,-1)[\phantom{Y_2}``=]{450}1r
 \putmorphism(-150,-450)(1,0)[(B,A)`(B,A')`(B,K)]{600}1a
\putmorphism(-180,-450)(0,-1)[\phantom{Y_2}``=]{450}1r
\putmorphism(1100,-450)(0,-1)[\phantom{Y_2}``=]{450}1r
 \putmorphism(-150,-900)(1,0)[(B,A)`(B,A'')`(B, K'K)]{1280}1a
\put(260,-670){\fbox{$(B,-)_{K'K}$}}
\efig}
$$
where $(B,-)_{K'K}$ is the 2-cell from the laxity of $(B,-)$ 
%(21)

\noindent $\bullet$ \quad   \axiom{($u, K'K$)}  
$$
\scalebox{0.86}{
\bfig
\putmorphism(-150,500)(1,0)[(B,A)`(B,A')`(B,K)]{600}1a
 \putmorphism(470,500)(1,0)[\phantom{F(A)}`(B,A'') `(B, K')]{600}1a
 \putmorphism(-150,50)(1,0)[(B,A)`(B,A'')`(B,K'K)]{1220}1a

\putmorphism(-180,500)(0,-1)[\phantom{Y_2}``=]{450}1r
\putmorphism(1080,500)(0,-1)[\phantom{Y_2}``=]{450}1r
\put(250,290){\fbox{$(B,-)_{K'K}$}}

\putmorphism(-150,-400)(1,0)[(\tilde B,A)`(\tilde B,A'') `(\tilde B,K'K)]{1200}1a

\putmorphism(-180,50)(0,-1)[\phantom{Y_2}``(u,A)]{450}1l
\putmorphism(1080,50)(0,-1)[\phantom{Y_3}``(u,A'')]{450}1r
\put(250,-160){\fbox{$(u,K'K)$}} 
\efig}
=
\scalebox{0.86}{
\bfig
\putmorphism(-150,500)(1,0)[(B,A)`(B,A')`(B,K)]{600}1a
 \putmorphism(470,500)(1,0)[\phantom{F(A)}`(B,A'') `(B, K')]{600}1a

 \putmorphism(-150,50)(1,0)[(\tilde B, A)`(\tilde B,A')`(\tilde B,K)]{600}1a
 \putmorphism(470,50)(1,0)[\phantom{F(A)}`(B'',\tilde A) `(\tilde B,K')]{620}1a

\putmorphism(-180,500)(0,-1)[\phantom{Y_2}``(u,A)]{450}1l
\putmorphism(450,500)(0,-1)[\phantom{Y_2}``]{450}1r
\putmorphism(300,500)(0,-1)[\phantom{Y_2}``(u,A')]{450}0r
\putmorphism(1080,500)(0,-1)[\phantom{Y_2}``(u,A'')]{450}1r
\put(-40,280){\fbox{$(u,K)$}}
\put(620,280){\fbox{$(u,K')$}}

\putmorphism(-150,-400)(1,0)[(\tilde B, A)`(\tilde B,A'') `(\tilde B,K'K)]{1200}1a

\putmorphism(-180,50)(0,-1)[\phantom{Y_2}``=]{450}1l
\putmorphism(1080,50)(0,-1)[\phantom{Y_3}``=]{450}1r
\put(260,-160){\fbox{$(\tilde B,-)_{K'K}$}}

\efig}
$$

%(12) 
\noindent $\bullet$ \quad  \axiom{($k'k, U$)} 
$$
\scalebox{0.86}{
\bfig
\putmorphism(-150,500)(1,0)[(B,A)`(B',A)`(k,A)]{600}1a
 \putmorphism(450,500)(1,0)[\phantom{F(A)}`(B'',A) `(k',A)]{620}1a

 \putmorphism(-150,50)(1,0)[(B,\tilde A)`(B',\tilde A)`(k,\tilde A)]{600}1a
 \putmorphism(470,50)(1,0)[\phantom{F(A)}`(B'',\tilde A) `(k',\tilde A)]{620}1a

\putmorphism(-180,500)(0,-1)[\phantom{Y_2}``(B,U)]{450}1l
\putmorphism(450,500)(0,-1)[\phantom{Y_2}``]{450}1r
\putmorphism(240,500)(0,-1)[\phantom{Y_2}``(B',U)]{450}0r
\putmorphism(1080,500)(0,-1)[\phantom{Y_2}``(B'',U)]{450}1r
\put(-40,280){\fbox{$(k,U)$}}
\put(620,280){\fbox{$(k',U)$}}

\putmorphism(-150,-400)(1,0)[(B,\tilde A)`(B'',\tilde A) `(k'k,\tilde A)]{1200}1a

\putmorphism(-180,50)(0,-1)[\phantom{Y_2}``=]{450}1l
\putmorphism(1080,50)(0,-1)[\phantom{Y_3}``=]{450}1r
\put(260,-160){\fbox{$(-,\tilde A)_{k'k}$}}
\efig}
=
\scalebox{0.86}{
\bfig
\putmorphism(-150,500)(1,0)[(B,A)`(B',A)`(k,A)]{600}1a
 \putmorphism(450,500)(1,0)[\phantom{F(A)}`(B'',A) `(k',A)]{620}1a
 \putmorphism(-150,50)(1,0)[(B,A)`(B'',A)`(k'k, A)]{1220}1a

\putmorphism(-180,500)(0,-1)[\phantom{Y_2}``=]{450}1r
\putmorphism(1080,500)(0,-1)[\phantom{Y_2}``=]{450}1r
\put(260,290){\fbox{$(-,A)_{k'k}$}}

\putmorphism(-150,-400)(1,0)[(B,\tilde A)`(B'',\tilde A) `(k'k,\tilde A)]{1200}1a

\putmorphism(-180,50)(0,-1)[\phantom{Y_2}``(B,U)]{450}1l
\putmorphism(1080,50)(0,-1)[\phantom{Y_3}``(B'',U)]{450}1r
\put(300,-180){\fbox{$(k'k,U)$}} 
\efig}
$$%\quad\text{and}\quad 

\noindent $\bullet$ \quad   \axiom{($\frac{u}{u'}, K$)} \qquad $(\frac{u}{u'}, K)=\frac{(u,K)}{(u', K)}$  \qquad\text{and}\qquad 
$\bullet$ \quad  \axiom{($k,\frac{U}{U'}$)}  \qquad $(k,\frac{U}{U'})=\frac{(k,U)}{(k,U')}$ 
%\noindent $\bullet$ \qquad ($(k,\frac{U}{U'})$) \label{k,UU'}  \quad $(k,\frac{U}{U'})=\frac{(k,U)}{(k,U')}$ 

%(22)
\noindent $\bullet$ \quad  \axiom{($u,\frac{U}{U'}$)}
$$(u,\frac{U}{U'})=
\scalebox{0.86}{
\bfig
 \putmorphism(-150,500)(1,0)[(B,A)`(B,A) `=]{600}1a
\putmorphism(-180,500)(0,-1)[\phantom{Y_2}`(B, \tilde A) `(B,U)]{450}1l
\put(0,50){\fbox{$(u,U)$}}
\putmorphism(-150,-400)(1,0)[(\tilde B, \tilde A)`(\tilde B, \tilde A) `=]{640}1a
\putmorphism(-180,50)(0,-1)[\phantom{Y_2}``(u,\tilde A)]{450}1l
\putmorphism(450,50)(0,-1)[\phantom{Y_2}``(\tilde B, U)]{450}1r
\putmorphism(450,500)(0,-1)[\phantom{Y_2}`(\tilde B, A) `(u,A)]{450}1r
\putmorphism(-820,50)(1,0)[(B, \tilde A)``=]{520}1a
\putmorphism(-820,50)(0,-1)[\phantom{(B, \tilde A')}``(B,U')]{450}1l
\putmorphism(-820,-400)(0,-1)[(B, \tilde A')`(\tilde B, \tilde A')`(u,\tilde A')]{450}1l
\putmorphism(-820,-850)(1,0)[\phantom{(B, \tilde A)}``=]{520}1a
\putmorphism(-150,-400)(0,-1)[(\tilde B, \tilde A)`(\tilde B, \tilde A') `(\tilde B, U')]{450}1r
\put(-650,-630){\fbox{$(u,U')$}}
\efig}
$$

\noindent $\bullet$ \quad  \axiom{($\frac{u}{u'},U$)} 
$$(\frac{u}{u'},U)=
\scalebox{0.86}{
\bfig
 \putmorphism(-150,500)(1,0)[(B,A)`(B,A) `=]{600}1a
\putmorphism(-180,500)(0,-1)[\phantom{Y_2}`(B, \tilde A) `(B,U)]{450}1l
\put(0,50){\fbox{$(u,U)$}}
\putmorphism(-150,-400)(1,0)[(\tilde B, \tilde A)` `=]{500}1a
\putmorphism(-180,50)(0,-1)[\phantom{Y_2}``(u,\tilde A)]{450}1l
\putmorphism(450,50)(0,-1)[\phantom{Y_2}`(\tilde B, \tilde A)`(\tilde B, U)]{450}1r
\putmorphism(450,500)(0,-1)[\phantom{Y_2}`(\tilde B, A) `(u,A)]{450}1r
\putmorphism(450,50)(1,0)[\phantom{(B, \tilde A)}`(\tilde B, A)`=]{620}1a
\putmorphism(1070,50)(0,-1)[\phantom{(B, \tilde A')}``(u',A)]{450}1r
\putmorphism(1070,-400)(0,-1)[(\tilde B', A)`(\tilde B', \tilde A)`(\tilde B', U)]{450}1r
\putmorphism(450,-850)(1,0)[\phantom{(B, \tilde A)}``=]{520}1a
\putmorphism(450,-400)(0,-1)[\phantom{(B, \tilde A)}`(\tilde B', \tilde A) `(U',\tilde A)]{450}1l
\put(600,-630){\fbox{$(u',U)$}}
\efig}
$$

%\item (11)\quad 
\noindent $\bullet$ \quad  \axiom{$(k,K)$-l-nat}  
$$
\scalebox{0.86}{
\bfig
 \putmorphism(-150,500)(1,0)[(B,A)`(B', A)`(k, A)]{600}1a
 \putmorphism(450,500)(1,0)[\phantom{A\ot B}`(B', A') `(B', K)]{680}1a
 \putmorphism(-150,50)(1,0)[(B,A)`(B,A')`(B,K)]{600}1a
 \putmorphism(450,50)(1,0)[\phantom{A\ot B}`(B', A') `(k, A')]{680}1a

\putmorphism(-180,500)(0,-1)[\phantom{Y_2}``=]{450}1r
\putmorphism(1100,500)(0,-1)[\phantom{Y_2}``=]{450}1r
\put(350,260){\fbox{$(k,K)$}}

\putmorphism(-150,-400)(1,0)[(\tilde B,A)`(\tilde B,A') `(\tilde B, K)]{640}1a
 \putmorphism(450,-400)(1,0)[\phantom{A'\ot B'}` (\tilde B',A') `(l,A')]{680}1a

\putmorphism(-180,50)(0,-1)[\phantom{Y_2}``(u,A)]{450}1l
\putmorphism(450,50)(0,-1)[\phantom{Y_2}``]{450}1r
\putmorphism(300,50)(0,-1)[\phantom{Y_2}``(u,A')]{450}0r
\putmorphism(1100,50)(0,-1)[\phantom{Y_2}``]{450}1r
\putmorphism(1080,50)(0,-1)[\phantom{Y_2}``(v,A')]{450}0r
\put(-20,-180){\fbox{$(u,K)$}}
\put(660,-180){\fbox{$(\omega,A')$}}

\efig}
\quad=\quad
\scalebox{0.86}{
\bfig
 \putmorphism(-150,500)(1,0)[(B,A)`(B', A)`(k, A)]{600}1a
 \putmorphism(450,500)(1,0)[\phantom{A\ot B}`(B', A') `(B', K)]{680}1a

 \putmorphism(-150,50)(1,0)[(\tilde B,A)`(\tilde B',A)`(l, A)]{600}1a
 \putmorphism(450,50)(1,0)[\phantom{A\ot B}`(\tilde B',A') `(\tilde B',K)]{680}1a
%.....
\putmorphism(-180,500)(0,-1)[\phantom{Y_2}``]{450}1l
\putmorphism(-160,500)(0,-1)[\phantom{Y_2}``(u,A)]{450}0l
\putmorphism(450,500)(0,-1)[\phantom{Y_2}``]{450}1l
\putmorphism(610,500)(0,-1)[\phantom{Y_2}``(v,A)]{450}0l %470
\putmorphism(1120,500)(0,-1)[\phantom{Y_3}``(v,A')]{450}1r
\put(-40,270){\fbox{$(\omega,A)$}} %(0,-180)
\put(650,270){\fbox{$(v,K)$}}
%.....
\putmorphism(-150,-400)(1,0)[(\tilde B,A)`(\tilde B,A') `(\tilde B, K)]{640}1a
 \putmorphism(450,-400)(1,0)[\phantom{A'\ot B'}` (\tilde B',A') `(l,A')]{680}1a

\putmorphism(-180,50)(0,-1)[\phantom{Y_2}``=]{450}1l
\putmorphism(1120,50)(0,-1)[\phantom{Y_3}``=]{450}1r
\put(300,-200){\fbox{$(l,K)$}}

\efig}
$$

\noindent $\bullet$ \quad  \axiom{$(k,K)$-r-nat}  
$$
\scalebox{0.86}{
\bfig
 \putmorphism(-150,500)(1,0)[(B,A)`(B', A)`(k, A)]{600}1a
 \putmorphism(450,500)(1,0)[\phantom{A\ot B}`(B', A') `(B', K)]{680}1a
 \putmorphism(-150,50)(1,0)[(B,A)`(B,A')`(B,K)]{600}1a
 \putmorphism(450,50)(1,0)[\phantom{A\ot B}`(B', A') `(k, A')]{680}1a

\putmorphism(-180,500)(0,-1)[\phantom{Y_2}``=]{450}1r
\putmorphism(1100,500)(0,-1)[\phantom{Y_2}``=]{450}1r
\put(350,260){\fbox{$(k,K)$}}

\putmorphism(-180,50)(0,-1)[\phantom{Y_2}``(B,U)]{450}1l
\putmorphism(450,50)(0,-1)[\phantom{Y_2}``]{450}1r
\putmorphism(300,50)(0,-1)[\phantom{Y_2}``(B,V)]{450}0r
\putmorphism(1100,50)(0,-1)[\phantom{Y_2}``]{450}1r
\putmorphism(1080,50)(0,-1)[\phantom{Y_2}``(B',V)]{450}0r
\put(-20,-180){\fbox{$(B,\zeta)$}}
\put(660,-180){\fbox{$(k,V)$}}

\putmorphism(-150,-400)(1,0)[(B,\tilde A)`(B,\tilde A') `(B,L)]{640}1a
 \putmorphism(450,-400)(1,0)[\phantom{A'\ot B'}` (B',\tilde A') `(k,\tilde A')]{680}1a
\efig}
\quad=\quad
\scalebox{0.86}{
\bfig
 \putmorphism(-150,500)(1,0)[(B,A)`(B', A)`(k, A)]{600}1a
 \putmorphism(450,500)(1,0)[\phantom{A\ot B}`(B', A') `(B', K)]{680}1a

 \putmorphism(-150,50)(1,0)[(B,\tilde A)`(B',\tilde A)`(k,\tilde A)]{600}1a
 \putmorphism(450,50)(1,0)[\phantom{A\ot B}`(B',\tilde A') `(B',L)]{680}1a

\putmorphism(-180,500)(0,-1)[\phantom{Y_2}``]{450}1l
\putmorphism(-160,500)(0,-1)[\phantom{Y_2}``(B,U)]{450}0l
\putmorphism(450,500)(0,-1)[\phantom{Y_2}``]{450}1l
\putmorphism(610,500)(0,-1)[\phantom{Y_2}``(B',U)]{450}0l %470
\putmorphism(1120,500)(0,-1)[\phantom{Y_3}``(B',V)]{450}1r
\put(-40,270){\fbox{$(k,U)$}} %(0,-180)
\put(650,270){\fbox{$(B', \zeta)$}}

\putmorphism(-180,50)(0,-1)[\phantom{Y_2}``=]{450}1l
\putmorphism(1120,50)(0,-1)[\phantom{Y_3}``=]{450}1r
\put(300,-200){\fbox{$(k,L)$}}

\putmorphism(-150,-400)(1,0)[(B,\tilde A)`(B,\tilde A') `(B,L)]{640}1a
 \putmorphism(450,-400)(1,0)[\phantom{A'\ot B'}` (B',\tilde A') `(k,\tilde A')]{680}1a
\efig}
$$

%(22)\quad
\noindent $\bullet$ \quad  \axiom{$(u,U)$-l-nat} 
$$
\scalebox{0.86}{
\bfig
 \putmorphism(-150,500)(1,0)[(B,A)`(B,A) `=]{600}1a
 \putmorphism(450,500)(1,0)[(B,A)` `(k,A)]{450}1a
\putmorphism(-180,500)(0,-1)[\phantom{Y_2}`(B, \tilde A) `(B,U)]{450}1l
\put(-20,50){\fbox{$(u,U)$}}
\putmorphism(-170,-400)(1,0)[(\tilde B, \tilde A)` `=]{480}1a
\putmorphism(-180,50)(0,-1)[\phantom{Y_2}``(u,\tilde A)]{450}1l
\putmorphism(450,50)(0,-1)[\phantom{Y_2}`(\tilde B, \tilde A)`(\tilde B, U)]{450}1l
\putmorphism(450,500)(0,-1)[\phantom{Y_2}`(\tilde B, A) `(u,A)]{450}1l
\put(600,260){\fbox{$(\omega,A)$}}
\putmorphism(430,50)(1,0)[\phantom{(B, \tilde A)}``(l, A)]{500}1a
\putmorphism(1070,50)(0,-1)[\phantom{(B, A')}`(\tilde B', \tilde A)`]{450}1r
\putmorphism(1050,50)(0,-1)[``(\tilde B',U)]{450}0r
\putmorphism(1070,500)(0,-1)[(B', A)`(\tilde B', A)`]{450}1r
\putmorphism(1050,500)(0,-1)[``(v,A)]{450}0r
\putmorphism(450,-400)(1,0)[\phantom{(B, \tilde A)}``(l, \tilde A)]{500}1a
\put(600,-170){\fbox{$(l,U)$}}
\efig}
\quad=\quad
\scalebox{0.86}{
\bfig
 \putmorphism(-150,500)(1,0)[(B,A)`(B',A) `(k,A)]{600}1a
 \putmorphism(450,500)(1,0)[\phantom{(B,A)}` `=]{450}1a
\putmorphism(-180,500)(0,-1)[\phantom{Y_2}`(B, \tilde A) `]{450}1l
\putmorphism(-160,500)(0,-1)[` `(B,U)]{450}0l
\put(620,50){\fbox{$(v,U)$}}
\putmorphism(-170,-400)(1,0)[(\tilde B, \tilde A)` `(l, \tilde A)]{470}1a
\putmorphism(-180,50)(0,-1)[\phantom{Y_2}``]{450}1l
\putmorphism(-160,50)(0,-1)[\phantom{Y_2}``(u,\tilde A)]{450}0l
\putmorphism(450,50)(0,-1)[\phantom{Y_2}`(\tilde B', \tilde A)`(v,\tilde A)]{450}1r
\putmorphism(450,500)(0,-1)[\phantom{Y_2}`(B', \tilde A) `(B',U)]{450}1r
\put(0,260){\fbox{$(k,U)$}}
\putmorphism(-190,50)(1,0)[\phantom{(B, \tilde A)}``(k, \tilde A)]{500}1a
\putmorphism(1070,50)(0,-1)[\phantom{(B, A')}`(\tilde B', \tilde A)`(\tilde B',U)]{450}1r
\putmorphism(1070,500)(0,-1)[(B', A)``(v,A)]{450}1r
\putmorphism(1100,500)(0,-1)[`(\tilde B', A)`]{450}0r
\putmorphism(450,-400)(1,0)[\phantom{(B, \tilde A)}``=]{500}1b
\put(0,-170){\fbox{$(\omega, \tilde A)$}}
\efig}
$$

\noindent $\bullet$ \quad \axiom{$(u,U)$-r-nat} 
$$
\scalebox{0.86}{
\bfig
 \putmorphism(-150,500)(1,0)[(B,A)`(B,A) `=]{600}1a
 \putmorphism(550,500)(1,0)[` `(B,K)]{380}1a
\putmorphism(-180,500)(0,-1)[\phantom{Y_2}`(B, \tilde A) `(B,U)]{450}1l
\put(-20,50){\fbox{$(u,U)$}}
\putmorphism(-170,-400)(1,0)[(\tilde B, \tilde A)` `=]{480}1a
\putmorphism(-180,50)(0,-1)[\phantom{Y_2}``(u,\tilde A)]{450}1l
\putmorphism(450,50)(0,-1)[\phantom{Y_2}`(\tilde B, \tilde A)`(\tilde B, U)]{450}1l
\putmorphism(450,500)(0,-1)[\phantom{Y_2}`(\tilde B, A) `(u,A)]{450}1l
\put(620,280){\fbox{$(u,K)$}}
\putmorphism(430,50)(1,0)[\phantom{(B, \tilde A)}``(\tilde B,K)]{500}1a
\putmorphism(1070,50)(0,-1)[\phantom{(B, A')}`(\tilde B, \tilde A')`]{450}1r
\putmorphism(1090,50)(0,-1)[\phantom{(B, A')}``(\tilde B,V)]{450}0r
\putmorphism(1070,500)(0,-1)[(B, A')`(\tilde B, A')`]{450}1r
\putmorphism(1090,500)(0,-1)[``(u,A')]{450}0r
\putmorphism(450,-400)(1,0)[\phantom{(B, \tilde A)}``(\tilde B, L)]{500}1a
\put(620,-170){\fbox{$ (\tilde{B},\zeta)$ } } % ???
\efig}
\quad=\quad
\scalebox{0.86}{
\bfig
 \putmorphism(-150,500)(1,0)[(B,A)`(B,A') `(B,K)]{600}1a
 \putmorphism(450,500)(1,0)[\phantom{(B,A)}` `=]{450}1a
\putmorphism(-180,500)(0,-1)[\phantom{Y_2}`(B, \tilde A) `]{450}1l
\putmorphism(-160,500)(0,-1)[\phantom{Y_2}` `(B,U)]{450}0l
\put(620,50){\fbox{$(u,V)$}}
\putmorphism(-170,-400)(1,0)[(\tilde B, \tilde A)` `(\tilde B, L)]{500}1a
\putmorphism(-180,50)(0,-1)[\phantom{Y_2}``]{450}1l
\putmorphism(-160,50)(0,-1)[\phantom{Y_2}``(u,\tilde A)]{450}0l
\putmorphism(450,50)(0,-1)[\phantom{Y_2}`(\tilde B, \tilde A')`(u,\tilde A')]{450}1r
\putmorphism(450,500)(0,-1)[\phantom{Y_2}`(B, \tilde A') `(B,V)]{450}1r
\put(0,260){\fbox{$(B,\zeta)$}}
\putmorphism(-190,50)(1,0)[\phantom{(B, \tilde A)}``(B,L)]{500}1a
\putmorphism(1070,50)(0,-1)[\phantom{(B, A')}`(\tilde B, \tilde A')`(\tilde B,V)]{450}1r
\putmorphism(1070,500)(0,-1)[(B, A')``(u,A')]{450}1r
\putmorphism(1110,500)(0,-1)[`(\tilde B, A')`]{450}0r
\putmorphism(450,-400)(1,0)[\phantom{(B, \tilde A)}``=]{500}1b
\put(0,-170){\fbox{$(u,L)$}}
\efig}
$$
for any 2-cells 
\begin{equation} \eqlabel{omega-zeta}
\scalebox{0.86}{
\bfig
\putmorphism(-150,175)(1,0)[B` B'`k]{450}1a
\putmorphism(-150,-175)(1,0)[\tilde B`\tilde B' `l]{440}1b
\putmorphism(-170,175)(0,-1)[\phantom{Y_2}``u]{350}1l
\putmorphism(280,175)(0,-1)[\phantom{Y_2}``v]{350}1r
\put(0,-15){\fbox{$\omega$}}
\efig}
\quad\text{and}\quad
\scalebox{0.86}{
\bfig
\putmorphism(-150,175)(1,0)[A` A'`K]{450}1a
\putmorphism(-150,-175)(1,0)[\tilde A`\tilde A' `L]{440}1b
\putmorphism(-170,175)(0,-1)[\phantom{Y_2}``U]{350}1l
\putmorphism(280,175)(0,-1)[\phantom{Y_2}``V]{350}1r
\put(0,-15){\fbox{$\zeta$}}
\efig}
\end{equation}
in $\Bb$, respectively $\Aa$. 
\end{prop}

\vspace{1,3cm}

\begin{thm} \thlabel{lr central} \cite[Theorem 6.2]{Fem3} \\
Let $\Bb$ be a double category. The following are equivalent:
\begin{enumerate}
\item there is a pseudodouble functor $\F\colon\Bb\to\Pseudo_{ps}(\Bb, \Bb)$; 
\item there is a pseudodouble quasi-functor $H:\Bb\times\Bb\to\Bb$ (with families of pseudodouble functors 
$(-,A),(B,-)\colon\Bb\to\Bb$ for $A,B\in\Bb$, four families of 2-cells satisfying 20 axioms); 
\item the following hold (meaning that $\Bb$ is purely central):
\begin{enumerate} [i)]
\item $\Bb$ is binoidal with pseudodouble functors $A\ltimes -,-\rtimes B:\Bb\to\Bb$ for $A,B\in\Bb$; 
\item every 1h-cell $K$ in $\Bb$ is left central via a horizontal pseudonatural transformation $K\ltimes-$ %$(-,K)$ 
and every 1h-cell $k$ is right central in $\Bb$ via a horizontal pseudonatural transformation $-\rtimes k$, %$(k,-)$, 
and it is \vspace{-0,2cm} 
$$A\ltimes-\vert_k=-\rtimes k\vert_A \quad\text{and}\quad -\rtimes B\vert_K=K\ltimes-\vert_B;$$
\item every 1v-cell $U$ is left central in $\Bb$ via a vertical pseudonatural transformation $U\ltimes-$ %$(-,U)$ 
and every 1v-cell $u$ is right central in $\Bb$ via a vertical pseudonatural transformation $-\rtimes u$, and %$(u,-)$.
it is \vspace{-0,2cm}
$$A\ltimes-\vert_u=-\rtimes u\vert_A \quad\text{and}\quad -\rtimes B\vert_U=U\ltimes-\vert_B;$$
\item every 2-cell $\zeta$ is left central via a modification $\zeta\ltimes-$ and every 2-cell $\omega$ is right central 
via a modification $-\rtimes\omega$, and it is \vspace{-0,2cm}
$$A\ltimes-\vert_\omega=-\rtimes\omega\vert_A \quad\text{and}\quad -\rtimes B\vert_\zeta=\zeta\ltimes-\vert_B;$$
\item it is 
$$K\ltimes-\vert_k=(-\rtimes k\vert_K)^{-1}, \qquad  U\ltimes-\vert_u=(-\rtimes u\vert_U)^{-1},$$
$$K\ltimes-\vert_u=-\rtimes u\vert_K \quad\text{and}\quad U\ltimes-\vert_k=-\rtimes k\vert_U$$ 
for all 1h-cells $K,k$ and 1v-cells $U,u$ in $\Bb$.  
\end{enumerate}
In particular, with notations as above, one has: 
\begin{itemize}
\item the 2-cell component $K\ltimes-\vert_k$ of the oplax (resp. lax) structure of the horizontal transformation $K\ltimes-$ coincides with the 2-cell component at $K$ of the lax (resp. oplax) structure of the transformation $-\rtimes k$ 
(which is $(-\rtimes k\vert_K)^{-1}$ in \equref{f lr}); 
\item the 2-cell component $U\ltimes-\vert_u$ of the lax (resp. oplax) structure of the vertical transformation $U\ltimes-$ coincides with the 2-cell component at $U$ of the oplax (resp. lax) structure of the transformation $-\rtimes u$ 
(which is $(-\rtimes u\vert_U)^{-1}$). % in \equref{v lr}
\end{itemize}
\end{enumerate}
\end{thm}

\bigskip

\textbf{Ethical approval declarations:} not applicable. \\

\textbf{Funding declaration:} none. 

%\thanks{The author was supported by the Science Fund of the Republic of Serbia, Grant No.~7749891, Graphical Languages - GWORDS. }


\begin{thebibliography}{99}


\bibitem{B} I. Bakovi\'c, {\em Bigroupoid 2-torsors}, Dissertation, LMU M\"unchen: Faculty of Mathematics, Computer Science and Statistics. 
(2008)

\bibitem{CGHF}
M. Capucci, B. Gavranovi\'c, J. Hedges, E. Fjeldgren Rischel: 
Towards Foundations of Categorical Cybernetics,  In Proceedings ACT 2021, arXiv:2211.01102

\bibitem{CG} M. Capucci, B. Gavranovi\'c, {\em Actegories for the Working Amthematician},
arXiv preprint available at https://arxiv.org/abs/2203.16351. 2022.


\bibitem{CJM} M. Capucci, D. Jaz Myers, {\em Constructing triple categories of cybernetic processes}, 

\bibitem{CJM2} M. Capucci, D. Jaz Myers, {\em Contextads as Wreaths; Kleisli, Para, and Span Constructions as Wreath Products}, 
arXiv:2410.21889. 


\bibitem{ChG} E. Cheng, N. Gurski, {\em The periodic table of n-categories ii: degenerate tricategories}, 
Cahiers de Topologie et G\'eom\'etrie Diff\'erentielle Cat\'egoriques {\bf 52}/2, (2011).


\bibitem{DS} B. Day, R. Street, {\em Monoidal Bicategories and Hopf Algebroids}, Adv. Math. {\bf 129}/1 (1997), 99--157. 


\bibitem{Fem2} B. Femi\'c, {\em Bifunctor Theorem and strictification tensor product for double categories with lax double functors}, 
Theory Appl. Categ. {\bf 39} (2023), 824--873. 


\bibitem{Fem3} B. Femi\'c, {\em Kleisli and premonoidal double categories}. Preprint arXiv:2401.17494.


\bibitem{FST} B. Fong, D. Spivak, R. Tuy\'eras, {\em Backprop as Functor: A compositional perspective on supervised learning}, 34th Annual ACM/IEEE Symposium on Logic in Computer Science (LICS) 2019, pp. 1-13, 2019. (arXiv:1711.10455, LICS'19)

\bibitem{GGV} N. Gambino, R. Garner, C. Vasilakopoulou, {\em Monoidal Kleisli Bicategories and the Arithmetic Product of 
Coloured Symmetric Sequences},  arXiv:2206.06858v2. 


\bibitem{GG} R. Garner, N. Gurski, {\em The low-dimensional structures formed by tricategories}, 
Journal Mathematical Proceedings of the Cambridge Philosophical Society {\bf 146} (2009), 551--589.


\bibitem{GSh} R. Garner, M. Shulman, {\em Enriched categories as a free cocompletion}, Adv. Math. {\bf 289} (2016), 1--94. 


\bibitem{GPS} R. Gordon, A. J. Power, R. Street, {\em Coherence for tricategories}, Memoirs of
the Amer. Math. Soc. {\bf 117}/558 (1995), 19, 28.

\bibitem {MG} M. Grandis, {\em Higher Dimensional Categories: From Double to Multiple Categories}, World Scientific (2019). 


\bibitem {GP:L} M. Grandis, R. Par\'e, {\em Limits in double categories}, Cahiers  Topol. 
G\'eom. Diff. Cat\'eg. {\bf 40} (1999), 162--220. 


\bibitem{GP:Adj} M. Grandis, R. Par\'e, {\em Adjoint for double categories}, Cahiers de Topologie et G\'eom\'etrie Diff\'erentielle
Cat\'egoriques {\bf 45}/3 (2004), 193--240. 


\bibitem{HT} C. Hermida, R. Tennent, {\em Monoidal indeterminates and categories of possible worlds}, Theoretical Computer Science vol 430. 2012. Preliminary version in MFPS 2009. doi:j.tcs.2012.01.001


\bibitem{HP} M. Hyland, J. Power, {\em Pseudo-commutative monads and pseudo-closed 2-categories}, J. Pure Appl. Algebra {\bf 175} (2002), 141--185. 


\bibitem{MV} J. Moeller, Ch. Vasilakopoulou, {\em Monoidal Grothendieck Construction}, 
Theory and Applications of Categories {\bf 35}/31 (2020), 1159--1207.

\bibitem{HF} H. Paquet, P. Saville, {\em Strong pseudomonads and premonoidal bicategories}, arXiv:2304.11014 (choose v1). 


\bibitem{HF1} H. Paquet, P. Saville, {\em Effectful semantics in 2-dimensional categories: premonoidal and Freyd bicategories}, 
	Electronic Proceedings in Theoretical Computer Science {\bf 397}/3 (2023), 190--209 %https://act2023.github.io/papers/paper76.pdf. 
(Proceedings of the Sixth International Conference on Applied Category Theory 2023), 
%University of Maryland, 31 July - 4 August 2023 
%(Electronic Proceedings in Theoretical Computer Science {\bf 397}), 
Open Publishing Association, https://doi.org/10.4204/EPTCS.397.12.

\bibitem{HF2} H. Paquet, P. Saville, {\em Effectful semantics in 2-dimensional categories: strong, commutative, and concurrent pseudomonads}, 
LICS '24: Proceedings of the 39th Annual ACM/IEEE Symposium on Logic in Computer Science Article {\bf 61}, 1 -- 15, 
https://doi.org/10.1145/3661814.3662130. 

\bibitem{P} D. Pavlovi\'c, {\em Categorical logic of names and abstraction in action calculi}, Math. Structures Comput. Sci. {\bf 7}/6 (1997) 619--637. 

\bibitem{Shul} M. Shulman, {\em Constructing symmetric monoidal bicategories}, arXiv: 1004.0993

%\bibitem{Shul1} M. Shulman, {\em Comparing composites of left and right derived functors}, New York J. Math. {\bf 17} (2011), 75--125. 

\bibitem{St} M. Stay, {\em Compact Closed Bicategories}, Theories and Applications of Categories {\bf 31}/26 (2016), 755--798.

\bibitem{UV} T. Uustalu, V. Vene, {\em Comonadic Notions of Computation}, Electronic Notes in Theoretical Computer Science {\bf 203} (2008), 
263--284. 

\end{thebibliography}
\end{document}